\documentclass[11pt]{article}
\title{Measuring edge importance: a quantitative analysis of the stochastic shielding approximation for random processes on graphs}
\author{Deena R.~Schmidt$^{1,2,*}$ and Peter J.~Thomas$^{1,2,3}$\\
\texttt{dschmidt@case.edu} and \texttt{pjthomas@case.edu}\\
Departments of $^1$Mathematics, Applied Mathematics and Statistics, \\$^2$Biology, $^3$Cognitive Science\\
Case Western Reserve University\\
Cleveland, Ohio 44106} 

\usepackage{color, subfigure, bbm}
\usepackage{graphicx,amssymb,amsmath,amsfonts}
\usepackage{epsf,wrapfig}

\newcommand{\galan}{Gal\'{a}n~}
\newcommand{\ER}{Erd\"{o}s-R\'{e}nyi~}
\newcommand{\cov}{\text{Cov}}
\newcommand{\ntot}{N_{\text{tot}}}
\newcommand{\R}{\mathbb{R}}
\newcommand{\N}{\mathbb{N}}
\renewcommand{\matrix}[2]{ \left(\begin{array}{#1} #2 \end{array}\right)}
\newtheorem{theorem}{Theorem}
\newtheorem{lemma}[theorem]{Lemma}

\begin{document}
\maketitle

\begin{abstract}
Mathematical models of cellular physiological mechanisms often involve random walks on graphs representing transitions within networks of functional states.  Schmandt and \galan recently introduced a novel \emph{stochastic shielding approximation} as a fast, accurate method for generating approximate sample paths from a finite state Markov process in which only a subset of states are observable.  For example, in ion channel models, such as the Hodgkin-Huxley or other conductance based neural models, a nerve cell has a population of ion channels whose states comprise the nodes of a graph, only some of which allow a transmembrane current to pass.  The stochastic shielding approximation consists of neglecting fluctuations in the dynamics associated with edges in the graph not directly affecting the observable states.  We consider the problem of finding the \emph{optimal} complexity reducing mapping from a stochastic process on a graph to an approximate process on a smaller sample space, as determined by the choice of a particular linear measurement functional on the graph.  The partitioning of ion channel states into conducting \textit{versus} nonconducting states provides a case in point.  In addition to establishing that Schmandt and Gal\'{a}n's approximation is in fact optimal in a specific sense, we use recent results from random matrix theory to provide heuristic error estimates for the accuracy of the stochastic shielding approximation for an ensemble of random graphs.  Moreover, we provide a novel quantitative measure of the contribution of individual transitions within the reaction graph to the accuracy of the approximate process.
\end{abstract}

*Corresponding author: \texttt{dschmidt@case.edu}

%Keywords: Markov process, complexity reduction, ion channel, Hodgkin-Huxley model, networks, random graphs

\newpage
%\tableofcontents
\section{Introduction}  

Many biological systems exhibit a combination of stochastic (chance, random, noisy) and deterministic dynamics \cite{Allen2003,Berg1993,CalvettiSomersalo2013}. For example, mathematical models involving stochastic processes arise in physiology \cite{ChennubhotlaBahar2007PLoSComputBiol,GeQianQian2012PhysRep,LaingLord2010,ZhangQianQian2012PhysRep}, ecology \cite{LugoMcKane2008PRE,nisbet1982,Snyder2003Ecology}, and genetic regulatory systems \cite{Elowitz+Levien+Sigga:2002:Science,GoldingPaulssonZawilskiCox2005Cell,LuShenZongHastyWolynes2006PNAS}.  Such mathematical models often originate as intrinsically complex, high-dimensional systems with many degrees of freedom, and many sources of variability. This inherent complexity presents two related challenges.  First, the essential dynamics of such systems may be hard to discern, and model reduction based on first principles for stochastic systems on complex networks is difficult.    Second, in order to predict the behavior of such systems under normal, pathological or experimental conditions, one must usually resort to numerical simulation studies.  Even with the tremendous progress in computing power over the last decades, intrinsically high dimensional stochastic systems remain prohibitive to simulate exhaustively.  Moreover, because of their dimensionality, the results of ensembles of stochastic simulations can be challenging to interpret.  Therefore, there is demand for efficient dimension reduction methods, both to provide high quality approximate numerical solutions to the stochastic evolution equations arising in high dimensional systems, and to provide an efficient conceptual framework for interpretation of the behavior of such systems.

In \cite{SchmandtGalan2012PRL}, Schmandt and \galan  introduced a \textit{stochastic shielding approximation} as a fast, accurate method for generating sample paths from a finite state Markov process in which only a subset of states are observable.  For example, in ion channel models, such as the Hodgkin-Huxley or other conductance based neural models, a nerve cell has a population of ion channels whose configurational states comprise the nodes of a graph, only some of which allow a transmembrane current to pass.  That is, each vertex of the ion channel state graph is labeled with a scalar ``conductance", which is either zero (non-conducting) or one (conducting).  The stochastic shielding approximation consists of neglecting fluctuations  associated with edges in the graph not directly affecting the observable states.  Specifically, stochastic processes representing transitions along edges connecting identically labeled states are replaced by mean transition rates, while processes associated with transitions connecting distinguishable states are unchanged.   This approximation is an example of complexity reduction in the sense of reducing a stochastic process generated by $K$ independent processes to a process on a smaller sample space, i.e.~generated by $K' < K$ processes.  Schmandt and \galan observe that, remarkably, the variance of the observable state (the membrane conductance) is almost identical in the reduced and the unreduced system.\footnote{Cf.~\cite{SchmandtGalan2012PRL} supplemental material section 5.}  While the approximate process does not faithfully reproduce \emph{all} aspects of the full process, it reproduces those features relevant to the neurophyisologist as well as to the larger biological system in which it is embedded.

Here we consider the problem of finding the \emph{optimal} complexity reducing mapping from a stochastic process on a graph to an approximate process on a smaller sample space, as determined by the choice of a particular linear measurement functional on the graph.  The partitioning of ion channel states into conducting versus nonconducting states provides a case in point.  In this paper we establish that Schmandt and  Gal\'{a}n's approximation is in fact optimal in a specific sense.  We derive a quantitative measure of the contributions of individual edges in the graph to the accuracy of the approximation, relative to the chosen measurement functional.  This  approach allows quantitative comparison of edge importance, and sheds light on the parametric dependence of relative edge importance, for instance in a voltage-gated ion channel.  In addition, we provide heuristic error estimates for the accuracy of the stochastic shielding approximation for an ensemble of symmetric random graphs.

Motivated by \cite{SchmandtGalan2012PRL}, we consider a multidimensional Ornstein-Uhlenbeck process on a graph $\mathcal{G}(\mathcal{V},\mathcal{E})$ with $n$ nodes and $m$ edges (reactions), and a linear measurement functional $M\in\R^n$.  We show that the stochastic shielding approximation is the most accurate dimension reduction possible among those neglecting fluctuations in the same number of underlying processes.  Neglecting a set of reactions in the full stochastic process $X$ creates an approximate process $\tilde{X}$ which matches the behavior of the full process in the mean, but deviates from the full process in the fluctuations.  

Extending this idea for an ensemble of symmetric directed graphs $\mathcal{G}(\mathcal{V},\mathcal{E})$, we establish two main results.
Lemma \ref{lem:decompose-U}, our first main result, allows us to find the optimal complexity reducing mapping from a stochastic process on a graph to an approximate process on a smaller sample space, as determined by the measurement  $M$.  Neglecting the fluctuations associated with a subset $\mathcal{E}'$ of the edge set $\mathcal{E}$ defines a new process $\tilde{X}(t)$ that deviates from the full process $X(t)$ by an amount that we call the deficiency, $U(t)=\tilde{X}(t)-X(t)$.  The observed error, given $M$, is then $M^\intercal U$; its mean is zero by construction, and its variance is $R=E[(M^\intercal U)^2]$.  In Lemma \ref{lem:decompose-U} we provide an exact formula for the contribution of the $k^{\text{th}}$ edge to this error.  This formula, which arises from a spectral decomposition of the graph Laplacian associated with the full process, gives an explicit criterion for choosing the $k$ most important edges in the graph, for any $0<k<m$.  

Our second main result, Theorem \ref{thm:MainResult}, applies this criterion to networks generated from a broad class of random graph ensembles with a randomly chosen binary measurement vector $M$.  We show that the importance measures of individual edges cluster tightly around one of two values.  For moderately large graphs, these clusters correspond with very high accuracy to Schmandt and Gal\'{a}n's stochastic shielding heuristic; an extremely accurate, reduced complexity approximation is obtained by neglecting fluctuations associated with edges connecting states that are indistinguishable under the measurement $M$.  We illustrate this result with a sample from the \ER random graph ensemble in \S \ref{ssec:ER}.

The analysis of Schmandt and \galan focused on accurate, efficient approximation of Markov processes arising from ion channel models.  In \S \ref{sec:HH} we apply our analysis to processes on two graphs arising from the classical Hodgkin-Huxley system of ion channels: the 5-state model for the voltage-gated potassium channel, and the 8-state model for the voltage-gated sodium channel.  In a more general setting, the transition rates connecting adjacent states in these models are voltage-dependent.  Here we restrict attention to the stationary case, corresponding biologically to the behavior of the channels under ``voltage clamped" conditions. For both the voltage-gated potassium and voltage-gated sodium channel state graphs we show that our ranking reproduces the Schmandt-\galan stochastic shielding heuristic over all physiologically relevant  voltages. This example also demonstrates that our results apply to graphs with non-symmetric adjacency matrices, as well as to the symmetric case.

In \S \ref{sec:discuss} we discuss possible extensions of our results to examples including signal transduction networks and calcium-induced calcium release models, as well as systems with graded rather than binary measurement functionals.

\section{Model}\label{sec:model}

\subsection{Connection to the Population Process}
\label{ssec:PopulationProcess}
We develop our results in the context of stationary Ornstein-Uhlenbeck processes.  In contrast, Schmandt and \galan \cite{SchmandtGalan2012PRL} introduced stochastic shielding in the broader context of density dependent random walks on a graph from which our OU process arises as a large population approximation.
To set the stage before moving to the OU process framework, we briefly describe a population process on a graph of the type considered by Schmandt and Gal\'{a}n.  In particular, we consider a stationary stochastic process on a directed graph $\mathcal{G}=(\mathcal{V},\mathcal{E})$ where $|\mathcal{V}|=n$ and $|\mathcal{E}|=m$, the number of nodes and edges in the graph, respectively.  Each directed edge corresponds to one reaction in the system. The $k^{th}$ edge $ij(k)=(i(k),j(k)) \in \mathcal{E}$ is defined to start at node $i(k)$ and end at node $j(k)$, so that the $k^{th}$ reaction effects a transition from state $i$ to state $j$.  Following \cite{Higham2008SIREV,LeeOthmer2010JMB}, we let $\zeta_k$ be the stoichiometry vector associated with edge $ij(k)\in\mathcal{E}$.  That is, the $i^{th}$ component of $\zeta_k$ is -1, the $j^{th}$ component is 1, and all other components are zero.
\begin{eqnarray}\label{eq:stoich-vector}
\zeta_k = \matrix{r}{\zeta_k(1)\\\vdots\\\zeta_k(i)\\\vdots\\\zeta_k(j)\\\vdots\\\zeta_k(n)} &=& \matrix{r}{0\\\vdots\\-1\\\vdots\\1\\\vdots\\0}
\end{eqnarray}
Under stationary conditions, such as a population of ion channels under voltage clamp, the occupancy numbers of different states of a continuous time Markov process can be represented as the solution of the stochastic equation obtained from a random time-change representation in terms of Poisson processes \cite{AndersonKurtz2011chapter}.  
If $\alpha_{k}$ gives the instantaneous \textit{per capita} transition rate from state $i(k)$ to state $j(k)$, then the full Markov process is specified by a collection of independent standard (unit rate) Poisson processes $Y_{k}$ each representing the occurrence of $i(k)\to j(k)$ transitions as follows.  Letting $N(t)\in\N^n$ be the nonnegative integer-valued vector representing the number of individuals in each of $n$ states, we may write $N(t)$ as a sum of transitions occurring at random times specified by the collection of $Y_{k}$.  
\begin{equation}\label{eq:population-process}
N(t)=N(0)+\sum_{k\in\mathcal{E}}\zeta_k Y_k \left(\int_{0}^t\alpha_k N_{i(k)}(s)\,ds  \right)
\end{equation}
Because each transition preserves the total number of individuals (i.e.~the components of $\zeta_k$ sum to zero for each $k$), we have
$\sum_iN_i(t)=\ntot=\sum_iN_i(0)$ for all $t>0$.

In Appendix \ref{app:tau-leaping} we show that, provided $\ntot$ is sufficiently large, we can approximate the deviation of $N(t)$ from its mean $\bar{N}\in\R^n$ by a multidimensional, Gaussian, Ornstein-Uhlenbeck process $X(t)\in\R^n, X(t) \approx N(t)-\bar{N}$ which satisfies a stochastic differential equation of the form given in Equation \ref{eq:OUPforX} below.  
%
%[Quote final result of Appendix \ref{app:tau-leaping} - analogous expression for $X$.]
In particular, we show that $X(t)$ can be approximated by an SDE of the form
\begin{equation}\label{eq:SDE-complicated}
%dX(t) = \sum_{k\in\mathcal{E}} \zeta_k \left( [\bar{N}_{i(k)}+X_{i(k)}(t)]\alpha_k dt + \sqrt{\bar{N}_{i(k)}(t) \alpha_k} dW_k(t) \right).
% note that the first sum involving N cancels out to zero because of equilibrium!
dX(t) = \sum_{k\in\mathcal{E}} \zeta_k \left( X_{i(k)}(t)\alpha_k dt + \sqrt{\bar{N}_{i(k)}(t) \alpha_k} dW_k(t) \right).
\end{equation}

\subsection{Multidimensional Ornstein-Uhlenbeck Process}
\label{ssec:OUP}

To obtain our main mathematical result, we consider a multidimensional Ornstein-Uhlenbeck process $X\in\R^n$ on the directed graph $\mathcal{G}=(\mathcal{V},\mathcal{E})$ where $|\mathcal{V}|=n$ and $|\mathcal{E}|=m$. The state of the system at time $t$, $X(t)$, satisfies Equation \ref{eq:SDE-complicated}, which we write in the equivalent form 
\begin{equation}\label{eq:OUPforX}
dX=LX\,dt + B\,dW.
\end{equation}
Here $L = (A-D)^\intercal$ is the graph Laplacian ($A$ is the weighted adjacency matrix of $\mathcal{G}$ with entries $A_{ij}=\alpha_k>0$ if there is an edge from node $i(k)$ to $j(k)$ and zero otherwise, and $D$ is the diagonal matrix such that entry $D_{ii}=\sum_j A_{ij}$ is the out-degree of node $i$). $B$ is the $n \times m$ noise matrix, and $W\in\R^m$ is an $m$-dimensional Brownian motion, i.e.~each component $dW_k$ represents the increment of an independent standard Brownian motion capturing the fluctuations of the $k^{\text th}$ reaction about its mean\footnote{If $Q=(q_{ij})$ is the generator matrix of the stochastic process on the graph, with $q_{ij}=\alpha_k$ whenever $k$ is the edge leading from $i$ to $j$, then $L=Q^\intercal$.}.  Matrix $B$ decomposes into a sum over the $m$ reactions
\begin{equation}\label{eq:decompB}
B = \sum_{k=1}^m B_k
\end{equation}
such that the $k^{th}$ column of matrix $B_k = \sigma_k \zeta_k$ and all other columns of $B_k$ are zero.  

The stochastic shielding approximation for a system of the form given in Equation \ref{eq:OUPforX} amounts to preserving the mean, but neglecting the fluctuations, for the processes driving a subset of the reactions, i.e.~replacing $B$ with an alternative matrix $\tilde{B}$ obtained by replacing a subset of columns in $B$ with null vectors. The trajectories of the resulting SDE, $\tilde{X}(t)$ (see Equation \ref{eq:OUPforXss}), are approximations of the trajectories of the full system.

In order to compare different complexity reduction choices, we define the \emph{deficiency} of an approximation to be the difference between the true and approximate trajectories, $U(t)=\tilde{X}(t)-X(t)$, when projected onto the measurement functional of interest $M$.  As suggested by Schmandt and Gal\'{a}n, the \emph{stationary covariance of the projection of the deficiency on the measurement vector} provides an appropriate measure for comparing the quality of alternative reductions.  That is, we use  $R=\mbox{Cov}[M^\intercal  U]=\mbox{Cov}[M^\intercal (\tilde{X}-X)]$ as our error measure.
We focus on reductions that preserve the behavior of the system (Equation \ref{eq:OUPforX}) relative to a given linear measurement functional $M\in\R^n$.  In the case of ion channels, $M\in\{0,1\}^n$ represents the conductance of each channel state.  We consider the case of graded rather than binary measurements in \S \ref{sec:discuss}.
%Motivated by this case, we will focus on binary measurement vectors $M$.  
Whether binary or graded, the measurement vector identifies the stochastic process of interest as the projection $Y(t)=M^\intercal X(t)$.  

Formally, we consider two processes $X(t)$ (full process) and $\tilde{X}(t)$ (reduced process) defined on a common probability space  
$(\Omega,\mathcal{F}_t,P)$.  The sample space $\Omega=C[0,\infty)^n$, filtration $\mathcal{F}_t$, and Wiener measure $P$ are those associated with $m$ independent copies of the standard Brownian process.  The approximate process $\tilde{X}(t)$ has the same sample space $\Omega$ and is measurable with respect to the same filtration $\mathcal{F}_t$, but also with respect to a smaller filtration $\tilde{\mathcal{F}}_t\subset\mathcal{F}_t$ generated by the Wiener processes associated with a subset of edges of the graph. The covariance of the deficiency, then, is well-defined in terms of the underlying measure $P$ on the full probability space.

In Appendix \ref{app:OUP-covariance} we show the standard result \cite{Gardiner2009book} that the stationary covariance of the full process decomposes into a sum of the contributions from the $m$ different reaction processes:
\begin{equation}\label{eq:covariance-of-X-first-occurence}
\cov[X(t),X^\intercal (t)] = \lim_{t\to\infty}\int_0^t\sum_{k=1}^m \sigma_k^2 \exp[L(t-t')]\zeta_k\zeta_k^\intercal \exp[L^\intercal (t-t')]\,dt'.
\end{equation}
Similarly, the covariance function for the projection $Y(t)=M^\intercal X(t)$ also decomposes into a sum as above.  
Because the eigenvector corresponding to the leading (0) eigenvalue of $L$ has constant components, it is orthogonal to $\zeta_k$ for each $k$.  Therefore the corresponding eigenspace is contained in the kernel of the matrix $B_kB_k^\intercal$, for each $k$, which guarantees that the limit on the RHS of Equation \ref{eq:covariance-of-X-first-occurence} remains finite.

Neglecting a set of reactions $\mathcal{E}'\subset\mathcal{E}$ creates an approximate processes, $\tilde{X}(t)$, which matches the behavior of the full process in the mean, but deviates from the full process in the fluctuations.  This reduced process satisfies the following SDE
\begin{equation}\label{eq:OUPforXss}
d\tilde{X}=L\tilde{X}\,dt + \tilde{B}\,dW
\end{equation}
where $\tilde{B} = \sum_{k\in\mathcal{E}\backslash\mathcal{E}'} B_k$
sums over the edges we keep.  Given the linear measurement functional $M\in\R^n$ above, we define the approximate projection $\tilde{Y}(t)=M^\intercal \tilde{X}(t)$.  Note that in the case of an ion channel system, $M$ is binary so $Y$ and $\tilde{Y}$ just pull out the observable (i.e., conducting) states of each system.  

Neglecting a subset of reactions also introduces an error in the representation of the measurement $Y(t)$ versus $\tilde{Y}(t)$ due to the difference between $X(t)$ and $\tilde{X}(t)$.  Recall that $U(t) = \tilde{X}(t)-X(t)$ is the deficiency of the reduced model compared to the full model. Then $\tilde{Y}(t)-Y(t) = M^\intercal U(t)$, and $U(t)$ satisfies the SDE 
\begin{equation}\label{eq:sde-u-of-t}
dU=LU\,dt+(\tilde{B}-B)\,dW.
\end{equation}
It is important to note that the noise source $dW$ that appears in Equations \ref{eq:OUPforX} and \ref{eq:OUPforXss} refers to the \emph{same} noise process $W$ in both cases.  The deficiency of the approximation relative to the full process is given by taking the limit of the mean squared error (MSE) of $\tilde{Y}-Y$ (equivalent to the stationary covariance of $\tilde{Y}-Y$), which, as shown in the proof of Lemma \ref{lem:decompose-U}, is an expression of the sum over all neglected reactions.
%%%%%%%%%%%%%%%%%%%%%%%%%%%%
% Lemma 1
%%%%%%%%%%%%%%%%%%%%%%%%%%%%
\begin{lemma}\label{lem:decompose-U}
For a connected graph with a symmetric Laplacian $L$, let $X$ and $\tilde{X}$ be the full and reduced processes defined by Equations \ref{eq:OUPforX} and \ref{eq:OUPforXss}, respectively, and let $M\in\R^n$.  Let $\mathcal{E}'\subset \mathcal{E}$ be the subset of edges neglected in the definition of $\tilde{X}$.  Let $L$ be diagonalizable with eigenpairs  $\{(\lambda_i,v_i)\}_{i=1}^n$  listed with eigenvalues $\lambda_i$ in order of decreasing real part and $||v_i||_2=1$.  Then the stationary covariance of the discrepancy $\tilde{Y}-Y=M^\intercal(\tilde{X}-X)$ satisfies
\begin{equation}
R[\mathcal{E}']=\lim_{t\to\infty}\cov(\tilde{Y}-Y)=\sum_{k\in\mathcal{E}'}R_k\end{equation}
where
\begin{equation}\label{eq:error_Rk}
R_k=\sigma_k^2 \sum_{i=2}^n\sum_{j=2}^n
\left(\frac{-1}{\lambda_i+\lambda_j}\right) (M^\intercal v_i) (v_i^\intercal \zeta_k) (\zeta_k^\intercal v_j) (v_j^\intercal M).\end{equation}
\end{lemma}
We can rank the error terms $R_k$ in descending order, thereby ordering the corresponding reactions in terms of their ``importance".  The most important reaction is the one with the largest value of $R_k$; if neglected, it would  introduce the largest error.  See Appendix \ref{subsec:proof-lemma1} for the proof of Lemma \ref{lem:decompose-U}.
%
%By computing $R_k$ for all $k$, we can evaluate the contribution of each reaction to the accuracy of the whole process and prioritize which terms to keep in the approximation. 
Note that an individual term in the sum (\ref{eq:error_Rk}) will be zero if either $\zeta_k \perp v_i$ or if $M \perp v_i$ for a given eigenvector $v_i$.  Typically, however, these vectors will not be orthogonal.  Therefore, it is of interest to know how the values of $R_k$ are distributed for different examples: graphs of actual ion channel states such as those in the classical Hodgkin-Huxley model, and more generally, ensembles of random graphs.  In \S \ref{sec:HH}, we compute the distribution of $R_k$ for the graphs of the potassium and sodium channel states in the Hodgkin-Huxley model.  In \S \ref{sec:RandomEnsemble}, we consider an ensemble of random graphs such as the \ER ensemble with randomly assigned binary measurement vector $M$ and prove our main result which is a statement about the expected value of $R_k$.  

For a random graph ensemble, the eigenvectors of the graph Laplacian are distributed randomly on the unit sphere \cite{KnowlesYin2011ProbTheorRelFields,TaoVu2012RandMat}.  Hence, they are unlikely to be exactly orthogonal to either $\zeta_k$ or $M$.  Given a series of assumptions (see \S \ref{ssec:assumptions}) that are true for naturally occurring random ensembles such as the symmetric Gaussian and \ER ensembles, we state our main result.
%%%%%%%%%%%%%%%%%%%%%%%%%%%%
% Theorem 2
%%%%%%%%%%%%%%%%%%%%%%%%%%%%
\begin{theorem}\label{thm:MainResult}
Given an ensemble of symmetric directed irreducible graphs $\mathcal{G}(\mathcal{V},\mathcal{E})$ with $n$ nodes satisfying assumptions A0-A5 (see \S \ref{ssec:assumptions}), a binary measurement vector $M \in\{0,1\}^n$, and a stoichiometry vector $\zeta_k$ corresponding to the $k^{th}$ reaction, the mean squared error $R_k$ resulting from neglecting the $k^{th}$ reaction has expected value
\begin{equation}
E[R_k |M] = \frac{\sigma^2_k|M^\intercal \zeta_k|}{n\,C} +  O(n^{-q}), \textit{ as }n\to\infty, \textit{ for some } q>1
\end{equation}
\emph{where the constant $C$ depends on the mean edge weight.}
\end{theorem}
This result says that there are two kinds of connections between nodes in the graph: those with a mean $R_k$ value of order $1/n$ and those with a smaller mean value of order less than $n^{-q}$, where $q>1$ is driven by the fourth moment of the eigenvector components (see assumption 4a in \S \ref{ssec:assumptions} for details).  For the case of the Gaussian ensemble, $q=2$.  Empirically, for the \ER random graph ensemble, $q \approx 5/3$ (see discussion in \S \ref{ssec:ER} and also Figure \ref{fig:covariance_all}).  The first type of reaction connects differently labeled nodes in terms of their measurement under $M$; those are the important reactions in the graph.  The other type connects identically labeled nodes; those are unimportant and, hence, can be neglected under the stochastic shielding approximation without sacrificing much accuracy.  The proof of Theorem \ref{thm:MainResult} is given in \S \ref{ssec:ThmProof}.  Before discussing more complicated examples, we illustrate the decomposition of the full process into approximate subprocesses for a simple 3-state example in the next subsection.

\subsection{3-State Example}\label{ssec:3state}

\begin{figure}[ht!] 
\centering
\includegraphics[width=3in]{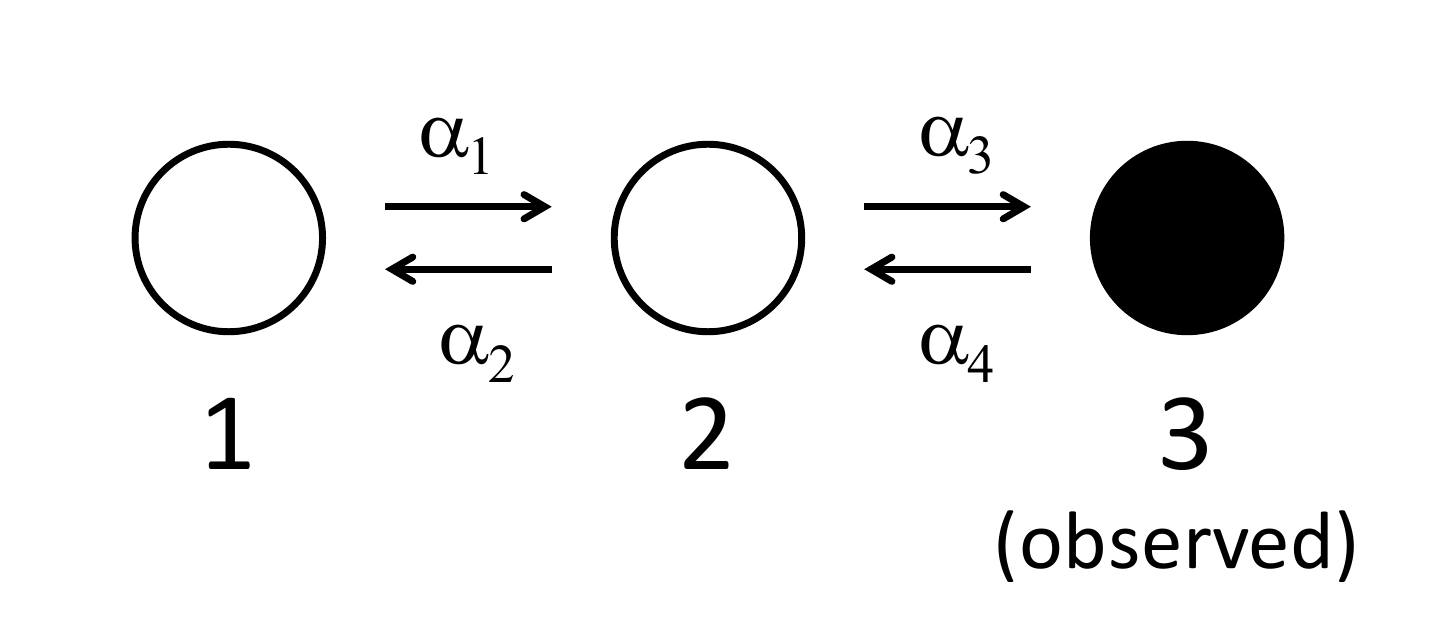}
\caption{Graph with 3 nodes and 4 reactions (edges) such that a transition from state $i(k)$ to state $j(k)$ happens at rate $\alpha_k$. For this example, we assume that only state 3 is observed.  This is the system shown in Figure 1 of Schmandt and \galan \cite{SchmandtGalan2012PRL}.}
\label{fig:3node-graph}
\end{figure}

We illustrate Schmandt and Gal\'{a}n's \cite{SchmandtGalan2012PRL} stochastic shielding heuristic with the following simple example they considered.  Figure \ref{fig:3node-graph} shows a 3-state chain which has adjacency matrix entries $A_{ij}=\alpha_k=1$ if there is an edge from $i(k)$ to $j(k)$ and zero otherwise.  State 3 is designated as the only observable state.  We think of this as the conducting state in an ion channel model.
%%%%%%%%%%%%%%%%%%%%%%%%%%%%%%%
% Figure 1
%%%%%%%%%%%%%%%%%%%%%%%%%%%%%%%
Table \ref{tab:3-node-stoich} illustrates the notation introduced in Equation \ref{eq:stoich-vector} for this case.  
%%%%%%%%%%%%%%%%%%%%%%%%%%%%%%%
% Table 1
%%%%%%%%%%%%%%%%%%%%%%%%%%%%%%%
\begin{table}[htbp]
   \centering
   $\begin{array}{|c|ccc|c|} % Column formatting, @{} suppresses leading/trailing space
   \hline
k&i(k)&\to&j(k)&M^\intercal \zeta_k\\ 
\hline
1&1&\to&2&0\\
2&2&\to&1&0\\
3&2&\to&3&+1\\
4&3&\to&2&-1\\
\hline
   \end{array}$
   \caption{Indexing of nodes and edges for the 3-state process, cf.~Equation \ref{eq:stoich-vector} and Figure \ref{fig:3node-graph}.  The first column gives the reaction number, the middle column gives the direction of the reaction, and the last column gives the contribution of the reaction to the measurement $Y=M^\intercal X$.}
   \label{tab:3-node-stoich}
\end{table}

In this case, we suppose $\sigma_k=1$ in the noise matrix $B$ and use the linear measurement functional $M=[0,0,1]^\intercal $ to pull out the third component of $X(t)$, yielding the projection $Y(t) = M^\intercal X(t)=X_3(t)$.  The vector $X(t)=(X_1(t),X_2(t),X_3(t))^\intercal $ gives the occupancy of the system states at time $t$ and satisfies the constant coefficient SDE given in Equation \ref{eq:OUPforX} with 
\begin{eqnarray}\label{eq:3state-system}
L&=& (A-D)^\intercal = \matrix{rrrr}
{-1 & 1 & 0\\
1 & -2 & 1\\
0 & 1 & -1}\\
B&=&\matrix{cccc}{\sigma_1\zeta_1&\sigma_2\zeta_2&\sigma_3\zeta_3&\sigma_4\zeta_4} = \matrix{rrrr}
{-1 & 1 & 0 & 0\\
1 & -1 & -1 & 1\\
0 & 0 &1 & -1} \\
W(t)&=&\matrix{c}{W_1(t)\\W_2(t)\\W_3(t)\\W_4(t)}
\end{eqnarray}
where the $W_k(t)$ are independent and identically distributed standard Brownian motions, and 
\begin{equation}
A=\matrix{rcr}
{0 & 1 & 0 \\
1 & 0 & 1\\
0 & 1 & 0}, \,\,\,
D=\matrix{rcr}
{1 & 0 & 0 \\
0 & 2 & 0\\
0 & 0 & 1}.
\end{equation}
Since we're assuming $\sigma_k=1$ for all $k$, the $k^{th}$ column of $B$ is exactly the stoichiometry vector associated with the $k^{th}$ reaction, and in particular, $B_k B_k^\intercal  = \zeta_k \zeta_k^\intercal $.  

The full process $X(t)$ has four stochastic transitions and a reduced process $\tilde{X}(t)$ is defined by keeping a subset of the four stochastic transitions.  We use the notation $\tilde{X}=X_{(i,j,k)}$ to explicitly define which columns of the full noise matrix $B$ are neglected in the approximate process, i.e. which stochastic transitions are neglected.  We are interested in the accuracy of the approximation of the trajectory itself.

\begin{figure}[htbp] 
	\centering
  \includegraphics[width=6in]{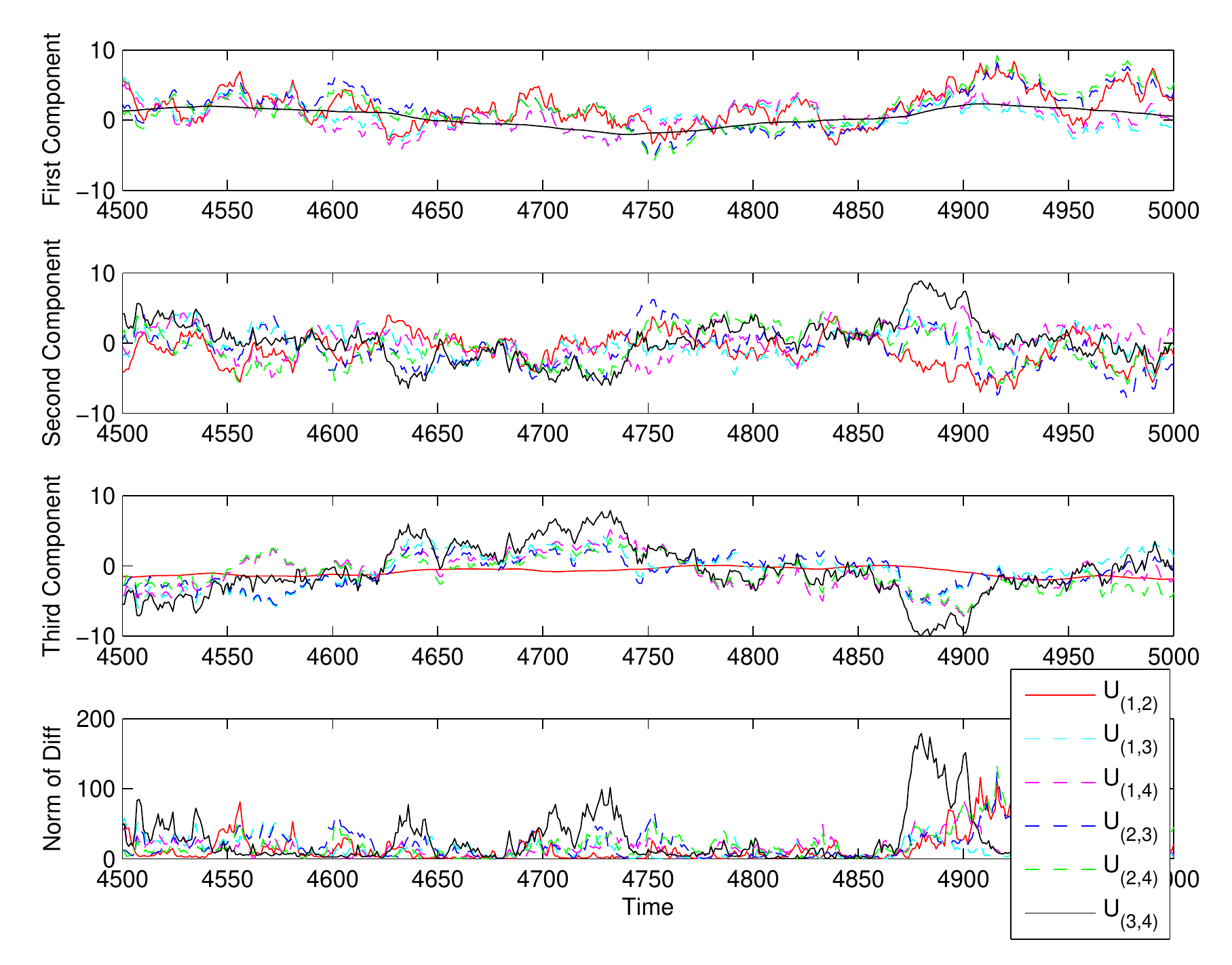} 
  \caption{Comparison of the deficiency $U_{(i,j)}(t)=X_{(i,j)}(t)-X(t)$ projected onto each component of the system of trajectories of an OU process on $\R^3$.  Top panel: $U_{(3,4)}$ is essentially zero which shows that reduced process $X_{(3,4)}$ is optimal for preserving the accuracy of the first component.  Second panel: no reduced process is optimal for preserving the accuracy of the second component.  Third panel: $U_{(1,2)}$ is essentially zero which shows that $X_{(1,2)}$ is optimal for preserving the accuracy of the third component (the conducting state in our 3-state example).  Bottom panel: squared norm of the deficiency $||U_{(i,j)}||^2 =||X_{(i,j)}-X||^2$.}
  \label{fig:three_state_trajectories}
\end{figure}

Figure \ref{fig:three_state_trajectories} illustrates the deficiency $U_{(i,j)}(t)=X_{(i,j)}(t)-X(t)$ between the full process and all possible two noise source reductions $X_{(i,j)}$ on the 3-state chain, as projected onto each of the three components in the system.
%%%%%%%%%%%%%%%%%%%%%%%%%%%%%%%
% Figure 2
%%%%%%%%%%%%%%%%%%%%%%%%%%%%%%%
The ``optimal complexity reduction" is not well defined in general because it is underspecified. For example, asking to reduce the norm of the deficiency $U$ while eliminating two of the four noise sources gives no preference between the six possible reductions.  Asking for the best reduction to preserve a specific component may give an answer: to preserve the trajectory as projected onto the $1^{st}$ component, keep the two noise sources that directly affect it (transitions between edges 1 and 2); for the $3^{rd}$ component, keep the other two (transitions between edges 3 and 4); for the $2^{nd}$ component there is no preference since it is affected directly by all transitions.  This gives an intuitive explanation of stochastic shielding consistent with Schmandt and Gal\'{a}n's explanation. 

If we fix a point in the underlying sample space (a choice of four Poisson processes $Y_k(t)$ in the system $N(t)$ or a choice of four white noise processes $dW_k(t)$ in the system $X(t)$) and then choose to neglect the fluctuations in two of the four, i.e. by replacing $Y_k(t)$ with $E[Y_k(t)]$ or $dW_k(t)$ with zero, respectively, then the question is: which choice leads to the most accurate representation of the process as seen by the measurement?  
%The error measure defined in the previous subsection refers to a pathwise comparison of the original and the reduced process.
   
By Lemma \ref{lem:decompose-U}, we have the following expression for the edge importance terms $R_k$:
\begin{equation}\label{eq:error_Rk_3state}
R_k = \sum_{i=2}^3\sum_{j=2}^3
\left(\frac{-1}{\lambda_i+\lambda_j}\right) (M^\intercal v_i) (v_i^\intercal \zeta_k) (\zeta_k^\intercal v_j) (v_j^\intercal M).\\
\end{equation}
Evaluating this expression for the measurement functional $M=[0,0,1]^\intercal$ yields
\begin{align*}
R_1 &= R_2 = 0.0417 \\
R_3 &= R_4 = 0.2917 
\end{align*}
Table \ref{tab:errors_3stateSDE} shows the stationary covariance of the discrepancy $M^\intercal U_{(i,j,k)}=M^\intercal (X_{(i,j,k)}-X)$ for all possible reduced processes $X_{(i,j,k)}$.  For instance, $X_{(1,2)}$ is the reduced process that neglects fluctuations in reactions 1 and 2 and the stationary covariance of $M^\intercal U_{(1,2)}$ is $R_1+R_2=0.0833$.  Note that $X_{(1,2)}$ is the optimal reduced process in terms of the Schmandt and \galan stochastic shielding approximation (among all approximations neglecting exactly two edges) for the 3-state chain.  
%%%%%%%%%%%%%%%%%%%%%%%%%%%%%%%
% Table 2
%%%%%%%%%%%%%%%%%%%%%%%%%%%%%%%
\begin{table}[ht]
\centering
\begin{tabular}{ l c r }
	\hline
  $M^\intercal U_{(i,j,k)}$ & $\sum R_{k'}$ & Value \\ \hline
  $M^\intercal U_{(1)}$ & $R_1$ & 0.0417 \\
  $M^\intercal U_{(2)}$ & $R_2$ & 0.0417 \\
  $M^\intercal U_{(3)}$ & $R_3$ & 0.2917 \\
  $M^\intercal U_{(4)}$ & $R_4$ & 0.2917 \\
  $M^\intercal U_{(1,2)}$ & $R_1+R_2$ & 0.0833* \\
  $M^\intercal U_{(3,4)}$ & $R_3+R_4$ & 0.583 \\
  $M^\intercal U_{(1,3)}$ & $R_1+R_3$ & 0.3333 \\
  $M^\intercal U_{(1,4)}$ & $R_1+R_4$ & 0.3333 \\
  $M^\intercal U_{(2,3)}$ & $R_3+R_3$ & 0.3333 \\
  $M^\intercal U_{(2,4)}$ & $R_2+R_4$ & 0.3333 \\
  $M^\intercal U_{(1,2,3)}$ & $R_1+R_2+R_3$ & 0.375 \\
  $M^\intercal U_{(1,2,4)}$ & $R_1+R_2+R_4$ & 0.375 \\
  $M^\intercal U_{(1,3,4)}$ & $R_1+R_3+R_4$ & 0.625 \\
  $M^\intercal U_{(2,3,4)}$ & $R_2+R_3+R_4$ & 0.625 \\ \hline
\end{tabular}
\caption{Table of discrepancies $M^\intercal U_{(i,j,k)}=M^\intercal (X_{(i,j,k)}-X)$ for the 3-state Markov process. The discrepancy $M^\intercal U_{(1,2)}$ (marked by $*$) corresponds to reduced process $X_{(1,2)}$ projected onto the third component, which is the optimal two-edge-neglecting approximation of $X$ for this example, in agreement with Schmandt and \galan \cite{SchmandtGalan2012PRL}.}
\label{tab:errors_3stateSDE}
\end{table}

\begin{figure}[ht!] 
\centering
\includegraphics[width=4in]{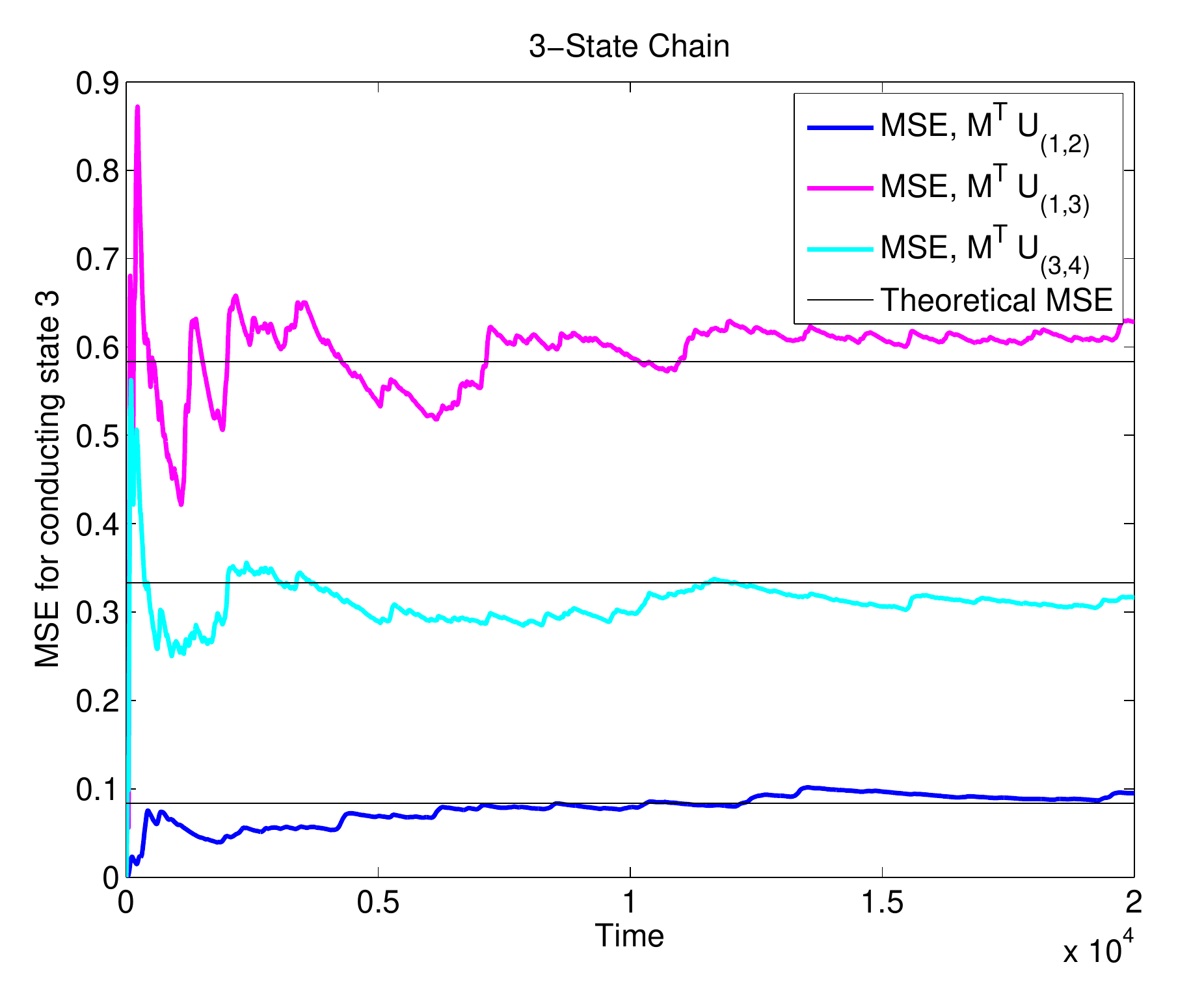}
\caption{Comparison of mean squared errors of $M^\intercal U_{(i,j)}$ in the 3-state chain, i.e. the projection of $U_{(i,j)}$ onto the third component. The theoretical MSE values are computed by summing the appropriate edge importance values $R_k$. $M^\intercal U_{(1,2)}$ has the smallest MSE out of the three classes of 2-noise source reduced process, showing that, as observed by Schmandt and Gal\'{a}n, $X_{(1,2)}$ is optimal at preserving the accuracy of the full process with respect to the third component of the system.}
\label{fig:three_state_chain_MSE}
\end{figure}

Figure \ref{fig:three_state_chain_MSE} shows the mean squared error as a function of time for $M^\intercal U_{(i,j)}(t)$ corresponding to the three classes of reduced processes $X_{(i,j)}(t)$ on the 3-state chain (i.e., the classes are $X_{(1,2)}$, $X_{(3,4)}$, and \{$X_{(1,3)}, X_{(1,4)}, X_{(2,3)}, X_{(2,4)}$\}, corresponding to the three different $M^\intercal U_{(i,j)}(t)$ values shown in Table \ref{tab:errors_3stateSDE} above).  The error function is shown with the theoretical MSE ($\sum_{k\in\mathcal{E}'} R_k$) for each case.  
%%%%%%%%%%%%%%%%%%%%%%%%%%%%%%%
% Figure 3
%%%%%%%%%%%%%%%%%%%%%%%%%%%%%%%
Therefore, since $M\zeta_1 = M\zeta_2 = 0$, $M\zeta_3 = 1$, and $M\zeta_4 = -1$ we confirm the claim made by Schmandt and \galan \cite{SchmandtGalan2012PRL} that reactions 3 and 4 are important whereas reactions 1 and 2 are unimportant in terms of stochastic shielding for this 3-state example.

%\newpage
\section{Analysis of Stochastic Shielding for a Random Graph Ensemble}
\label{sec:RandomEnsemble}

For any particular Ornstein-Uhlenbeck process on a graph, Lemma \ref{lem:decompose-U} provides the edge importance values $R_k$ (Equation \ref{eq:error_Rk}), which may be used to compute explicitly the contribution to the deficiency made by neglecting any particular reaction, relative to a given measurement vector $M$.  In order to make general observations about the stochastic shielding approximation, we now consider an ensemble of random graphs.  The proof of our main result (Theorem \ref{thm:MainResult}, restated below) will rely on properties of the joint distribution of components of eigenvectors of $L$, the graph Laplacian.  Previously, we used $i$ and $j$ to refer to the source and destination nodes in a reaction. In this section, we will adapt the notation so that edge $k$ is a reaction from node $l_-$ to node $l_+$, denoted by $l_{\pm}(k)\in\mathcal{E}$ (see Equation \ref{eq:zeta_k_2}).  In this section, $i$ and $j$ will instead index eigenvectors of $L$.  
\begin{eqnarray}\label{eq:zeta_k_2}
\zeta_k = \matrix{r}{\zeta_k(1)\\\vdots\\\zeta_k(l_-)\\\vdots\\\zeta_k(l_+)\\\vdots\\\zeta_k(n)} &=& \matrix{r}{0\\\vdots\\-1\\\vdots\\1\\\vdots\\0}
\end{eqnarray}
Because our methods combine heuristic numerical evidence with probabilistic calculations, we use ``$\approx$'' to represent ``heuristic equality".  Where precise order estimates are available, we use ``$O$'' notation.  For the reader's convenience, we restate Theorem \ref{thm:MainResult}.

%%%%%%%%%%%%%%%%%%%%%%%%%%%%
% Restate Theorem
%%%%%%%%%%%%%%%%%%%%%%%%%%%%
{\bf Theorem 2} \emph{Given an ensemble of symmetric directed irreducible graphs $\mathcal{G}(\mathcal{V},\mathcal{E})$ with $n$ nodes satisfying assumptions A0-A5 (see \S \ref{ssec:assumptions}), a binary measurement vector $M \in\{0,1\}^n$, and a stoichiometry vector $\zeta_k$ corresponding to the $k^{th}$ reaction, the mean squared error $R_k$ resulting from neglecting the $k^{th}$ reaction has expected value}
\begin{equation}
E[R_k |M] = \frac{\sigma^2_k|M^\intercal \zeta_k|}{n\,C} +  O(n^{-q}), \textit{ as }n\to\infty, \textit{ for some } q>1
\end{equation}
\emph{where the constant $C$ depends on the mean edge weight.}

In other words, since
\begin{equation}
  |M^\intercal \zeta_k| =\begin{cases}
    1, & \text{\small if reaction $k$ connects nodes with different $M$ values}\\
    0, & \text{\small if reaction $k$ connects nodes with the same $M$ value}
  \end{cases}
\end{equation}
reactions connecting nodes with identical values of $M$ have a small contribution to the error, so these reactions can be neglected under the stochastic shielding approximation.  This result relies on a list of assumptions which are described in detail below.  The proof of this theorem requires Lemma \ref{lem:random} which is stated after the assumptions and proved in Appendix \ref{ssec:proof_lemma2}.

%%%%%%%%%%%%%%%%%%%%%%%%%%%%%%%%%%%%%%%%%%%%%%%%%%%%%%
% Assumptions
%%%%%%%%%%%%%%%%%%%%%%%%%%%%%%%%%%%%%%%%%%%%%%%%%%%%%%
\subsection{Assumptions on the Random Graph Ensemble}
\label{ssec:assumptions}
We state a sequence of assumptions on the random graph ensemble needed to establish our main result.  Each assumption is reasonable for a broad class of graphs of interest, for reasons articulated in the Remarks following each assumption.  In several instances we impose on our random graph ensemble, as assumptions, properties that are known to hold for broad classes of random matrices, such as the Wigner ensemble \cite{KnowlesYin2011ProbTheorRelFields,TaoVu2012RandMat}.  The ensemble we consider is not equivalent to a generalized Wigner ensemble.  Nevertheless, for the reasons detailed below, it appears reasonable, that certain aspects of the eigenvector and eigenvalue distribution may be similar in the two cases.

\begin{figure}[ht!] 
\centering
	\includegraphics[width=4in]{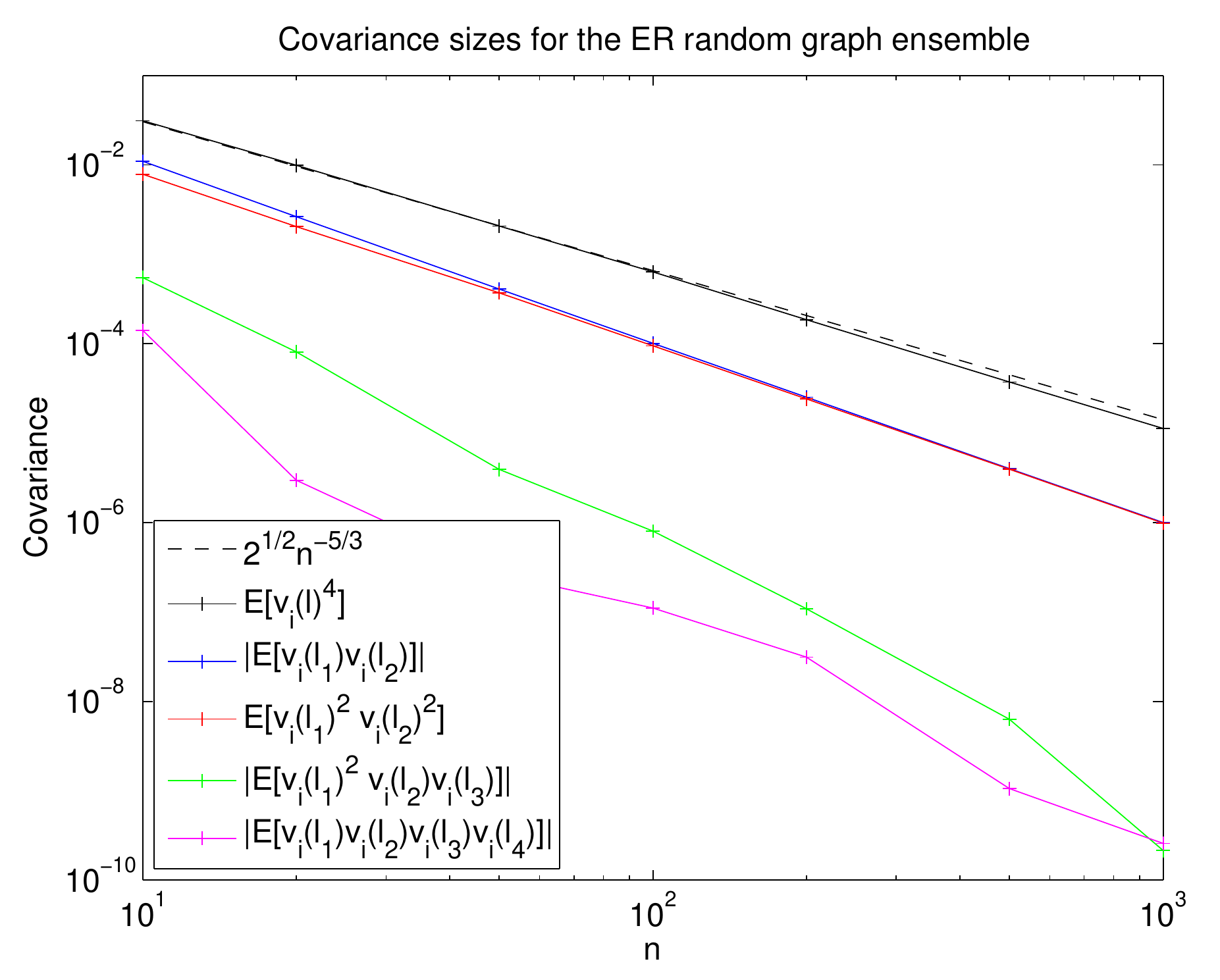}
\caption{Pairwise and fourth order covariance sizes of the eigenvector components of the graph Laplacian for the \ER random graph ensemble. To evaluate the fourth moment and the mixed moments listed in the legend, we computed the average value over $\ge100$ independent samples for each value of $n$.
%To evaluate the fourth moment for each value of $n$, we took a single sample and averaged over all $l_2$-normalized eigenvector components.  To evaluate mixed moments, we computed the average value over $100$ independent samples.  
Empirically, the expected value of $v_i(l)^4$ is approximately $\sqrt{2} n^{-5/3}$ (black); the dashed line is $\sqrt{2} n^{-5/3}$.  The absolute value of the expectation of $v_i(l_1)v_j(l_2)$ is $n^{-2}$ if $i=j$ (blue) and essentially $0$ if $i\ne j$ (data not shown; the average value was $10^{-19}$ or smaller).  The expectation of $v_i(l_1)^2 v_i(l_2)^2$ is approximately $n^{-2}$ (red).  The absolute value of the expectation of $v_i(l_1)^2 v_i(l_2) v_i(l_3)$ and $v_i(l_1) v_i(l_2) v_i(l_3) v_i(l_4)$ are both of order $n^{-3}$ (green and magenta).  This is numerical evidence for assumptions A3-A5 above.}
\label{fig:covariance_all}
\end{figure}

We consider an ensemble of symmetric directed irreducible graphs $\mathcal{G}=(\mathcal{V},\mathcal{E})$ with $|\mathcal{V}|=n$.  Let $\zeta_k$ be the stoichiometry vector corresponding to the $k^{th}$ reaction (Equation \ref{eq:zeta_k_2}) and let $(\lambda_i,v_i)$ denote the eigenpairs of the graph Laplacian $L=(A-D)^\intercal$ listed with eigenvalues in descending order.  We assume that the eigenvector components are $l_2$-normalized with mean 0 and variance $1/n$.  Suppose the measurement vector $M\in\{0,1\}^n$ contains at least one zero and at least one unit entry.\footnote{If $M$ has the same value for all nodes, the output is constant and the error is identically zero.  The expression in Theorem \ref{thm:MainResult} holds trivially so we ignore this case.}  We assume the following: 

\begin{enumerate}
\item [A0.] (Following \cite{DingJiang2010AnnalsApplProb}). Let $a_{ij} \ge 0$, the entries of the adjacency matrix, be random variables defined on a common probability space, with $\{a_{ij},1\le i < j\le n\}$ independent and identically distributed, 
with $a_{ij}=a_{ji}, E[a_{ij}]=\mu_A, V[a_{ij}]=\sigma_A^2>0$ for all $1\le i < j \le n$, and 
$\sup_{1\le i<j\le n} E\left|\left(a_{ij}-\mu_A\right)/\sigma_A\right|^\kappa<\infty$ for some $\kappa>0$.

\item [A1a.] The graph is drawn from a random ensemble with the property that the eigenvalues $\lambda_i$ and eigenvectors $v_i$ of the associated graph Laplacian are nearly independent.  That is, for any $i,j,k,l\in\{1,\dots,n \}$ and arbitrary measurable functions $f:\R^2\to\R$ and $g:\R^n \times \R^n \to\R$
\begin{equation}
E[f(\lambda_i,\lambda_j)g(v_k,v_l)]=E[f(\lambda_i,\lambda_j)]E[g(v_k,v_l)]+O\left(\frac{1}{n^4}\right), \text{ as }n\to\infty.
\end{equation}
\item [] \emph{Remark 1a:} A1 holds for the symmetric Gaussian ensemble as well as for the more general Wigner ensemble \cite{KnowlesYin2011ProbTheorRelFields,TaoVu2012RandMat}.  Indeed for these ensembles the eigenvalues and eigenvectors are independent.  The weaker assumption, that they are at most weakly correlated, appears reasonable for e.g.~the ensemble of graph Laplacians obtained from the symmetric \ER random graph ensemble.
\item [A1b.] The graph is drawn from a random ensemble with the property that the joint (eigenvalue, eigenvector) distribution is nearly invariant under permutation of eigenvectors.  That is, we assume that
\begin{enumerate}
\item [i.] for any pairs $i\not=j$ and $k\not=l$, 
\begin{equation}
E[f(\lambda_i,\lambda_j)g(v_i,v_j)]=E[f(\lambda_i,\lambda_j)g(v_k,v_l)]+O\left(\frac{1}{n^4}\right), \text{ as }n\to\infty
\end{equation}
\item [ii.] and for any $i$ and any $k$,
\begin{equation}
E[f(\lambda_i,\lambda_i)g(v_i,v_i)]=E[f(\lambda_i,\lambda_i)g(v_k,v_k)]+O\left(\frac{1}{n^4}\right), \text{ as }n\to\infty
\end{equation}
\end{enumerate}
\item [] \emph{Remark 1b:} The symmetric Gaussian and Wigner ensembles are fully invariant under permutation of eigenvectors, and the weaker assumption of near invariance appears reasonable for the \ER ensemble.  In particular, the pair $\left( \frac{-1}{\lambda_i+\lambda_j} \right), \left( M^\intercal v_i v_i^\intercal \zeta_k \zeta_k^\intercal v_j v_j^\intercal M \right)$ appearing in the definition of $R_k$ (Lemma \ref{lem:decompose-U}) are assumed to be approximately uncorrelated.  This assumption is reasonable by virtue of the approximate rotational symmetry of the eigenvector distribution under our choice of random graph model, which we expect to be close (heuristically) to the eigenvector distribution of the symmetric Gaussian ensemble \cite{KnowlesYin2011ProbTheorRelFields,TaoVu2012RandMat}.
\item [A2.] $E[v_i(l)]=0$ for any $i,l\in\{1,\dots,n\}$ where $v_i(l)$ denotes the $l^{th}$ component of the $i^{th}$ eigenvector.  
\item [] \emph{Remark 2a:} Note that $E[v_i(l)^2]=1/n$ by the $l_2$-normalization of the eigenvectors because $||v||_2 = \sqrt{\sum_{l=1}^n v(l)^2}=1$ for each eigenvector $v$. This normalization leaves a 2-fold ambiguity in the choice of eigenvector $v$.  Since $+v$ and $-v$ both have $||v||_2=1$, we choose randomly between them so that the first non-zero component is positive with probability 1/2.\footnote{In contrast, Tao and Vu \cite{TaoVu2012RandMat} always choose the first non-zero component to be positive to remove this ambiguity.}  
\item [] \emph{Remark 2b:} By the symmetry of our random graph ensemble under the symmetric group acting on the change of labels, A2 holds not just for the Gaussian and Wigner ensembles, but for any reasonable symmetric ensemble.  In particular, it holds for the symmetric \ER random graph ensemble. 
\item [A3.] For any $i,j\in\{2,\dots,n\}$ and $l,l'\in\{1,\dots,n\}$,
	\begin{enumerate}
	\item [a.] $E[v_i(l)v_j(l')]=O(n^{-3})$ as $n\to\infty$, for $i\ne j$.
	\item [b.] $E[v_i(l)v_i(l')]=O(n^{-2})$ as $n\to\infty$, for $l\ne l'$.
	\end{enumerate}
\item [] \emph{Remark 3:} Figure \ref{fig:covariance_all} provides numerical evidence for the plausibility of assumption A3 in the \ER case.  As described in the figure, the empirical expectation of $v_i(l)v_i(l')$ scales as $O(n^{-2})$ for $10\le n \le 1000$; over this range the empirical expectation of $v_i(l)v_j(l'), i\not=j$, is within machine error ($\le 10^{-19}$) of zero.
\item [A4.] For any $i\in\{2,\dots,n\}$ and $l,l'\in\{1,\dots,n\}$,
	\begin{enumerate}
	\item [a.] $E[v_i(l)^4]=O(n^{-q})$ as $n\to\infty$, for some $q > 1$.
	\item [b.] $E[v_i(l)^2 v_i(l')^2]=O(n^{-2})$ as $n\to\infty$, for $l\ne l'$.
	\end{enumerate}
\item [] \emph{Remark 4:} A4a holds for the Gaussian case for $q=2$.  For the \ER case, empirically we see that A4a holds for $q\approx 5/3$ as shown in Figure \ref{fig:covariance_all}.  Specifically, empirical evidence suggests that $E[v_i(l)^4] \approx \sqrt{2} n^{-5/3}$ in this case.
\item [A5.] Suppose that $p_1$, $p_2$, $p_3$, and $p_4$ are nonnegative integers with $\sum_{m=1}^4 p_m=4$, at least three of which are nonzero.  Then for any $i\in\{2,\dots,n\}$ and for any distinct components $\{l_1,l_2,l_3,l_4\}$
\begin{equation}\label{eq:assumption5}
E[(v_i(l_1))^{p_1}(v_i(l_2))^{p_2}(v_i(l_3))^{p_3}(v_i(l_4))^{p_4}]=O(n^{-3}) \text{ as }n\to\infty.
\end{equation}
\item [] \emph{Remark 5:} The reason for this assumption will become clear in the proof of Theorem \ref{thm:MainResult}.  It is similar in spirit to the four moment theorem for eigenvector components of a Wigner or Gaussian random matrix, different versions have been established by Tao and Vu \cite{TaoVu2012RandMat} and Knowles and Yin \cite{KnowlesYin2011ProbTheorRelFields}.  Figure \ref{fig:covariance_all} provides numerical evidence for the plausibility of A5 in the \ER case.
\end{enumerate}

%%%%%%%%%%%%%%%%%%%%%%%%%%%%
% Lemma 2
%%%%%%%%%%%%%%%%%%%%%%%%%%%%
\begin{lemma}\label{lem:random}
If assumptions A0-A5 hold, then as $n\to\infty$,
\begin{enumerate}
\item [A.] $E[M^\intercal v_i v_i^\intercal \zeta_k]=E\left[\sum_{l\in 1_M} v_i(l) (v_i(l_+)-v_i(l_-))\right]=\frac{1}{n} M^\intercal \zeta_k + O(n^{-2})$.
\item [B.] $E[M^\intercal v_i v_i^\intercal \zeta_k]^2 = E\left[\sum_{l\in 1_M} v_i(l) (v_i(l_+)-v_i(l_-))\right]^2 = \frac{1}{n^2}|M^\intercal\zeta_k| + O(n^{-4})$. 
\item [C.] $E[(M^\intercal v_i v_i^\intercal \zeta_k)^2] = E\left[\left(\sum_{l\in 1_M} v_i(l)\right)^2 (v_i(l_+)-v_i(l_-))^2\right] = O(n^{-q})$ for some $q>1$.
\end{enumerate}
\end{lemma}
Note that the exponent $q>1$ in part C is governed by the fourth moment of the eigenvector components of the graph Laplacian (see Assumption 4a). The proof of Lemma \ref{lem:random} is given in Appendix \ref{ssec:proof_lemma2}.

%%%%%%%%%%%%%%%%%%%%%%%%%%%%%%%
% Figure 4
%%%%%%%%%%%%%%%%%%%%%%%%%%%%%%%

%%%%%%%%%%%%%%%%%%%%%%%%%%%%%%%%%%%%%%%%%%%%%%%%%%%%%%
% Proof of Theorem
%%%%%%%%%%%%%%%%%%%%%%%%%%%%%%%%%%%%%%%%%%%%%%%%%%%%%%
\subsection{Proof of Main Theorem}
\label{ssec:ThmProof}

Suppose assumptions A0-A5 hold.  By Lemma \ref{lem:decompose-U}, $R_k$ denotes the contribution of the $k^{th}$ reaction to the deficiency of the approximate process.  Given the measurement vector $M$, we have (exactly)
\begin{eqnarray}\label{eq:Rk_exact}
E[R_k |M] &=& E\left[ \sigma_k^2 \sum_{i=2}^n\sum_{j=2}^n \left(\frac{-1}{\lambda_i+\lambda_j}\right) (M^\intercal v_i v_i^\intercal \zeta_k)(\zeta_k^\intercal v_j v_j^\intercal M) \right]
\end{eqnarray}
This expectation is taken over the space of symmetric directed irreducible graphs $\mathcal{G}(\mathcal{V},\mathcal{E})$ where edge $k$ is chosen at random from the set of ${n \choose 2}$ possible bidirectional edges.  If $l_{\pm}(k) \notin\mathcal{E}$, then $E[R_k |M]=0$.

If the graph Laplacian were drawn from a symmetric Gaussian ensemble (or Wigner ensemble, see \cite{KnowlesYin2011ProbTheorRelFields,TaoVu2012RandMat}), then the eigenvalues and the eigenvectors would be independent. For other ensembles we impose the weaker condition of \emph{near independence} (assumption A1a), which in this case means that for each $i\ge 2$ and $j\ge 2$, we assume
\begin{align*}
&E\left[\left(\frac{-1}{\lambda_i+\lambda_j}\right) (M^\intercal v_i v_i^\intercal \zeta_k)(\zeta_k^\intercal v_j v_j^\intercal M) \right]=\\
&
E\left[\left(\frac{-1}{\lambda_i+\lambda_j}\right)\right] E\left[ (M^\intercal v_i v_i^\intercal \zeta_k) (\zeta_k^\intercal v_j v_j^\intercal M) \right]+O\left(\frac{1}{n^4} \right), \text{ as }n\to\infty.
\end{align*}

Under assumption A1b, the joint distribution of eigenvalues and eigenvectors is approximately separable into the product of two measures, one for the eigenvalues and a second for the eigenvectors.  In this case the expectation $E\left[\left(\frac{-1}{\lambda_i+\lambda_j}\right)\right]$ in the sum (\ref{eq:Rk_exact}) can be replaced by its average, 
%
% Moreover, because of the near independence of the eigenvalues and eigenvectors and symmetry with respect to permutation of indices, we may replace the prefactor $1/(\lambda_i+\lambda_j)$ with its average over all terms in the sum (\ref{eq:Rk_exact}), 
\begin{equation}\label{eq:define-S}
S=\frac{1}{(n-1)^2}\sum_{i=2}^n\sum_{j=2}^n\frac{-1}{\lambda_i+\lambda_j},
\end{equation}
to obtain
\begin{eqnarray}
E[R_k |M] &=& \sigma_k^2 \, E[S]  \, E\left[\sum_{i=2}^n\sum_{j=2}^n(M^\intercal v_i v_i^\intercal \zeta_k)(\zeta_k^\intercal v_j v_j^\intercal M) \right]+O\left(\frac{1}{n^2}\right).
\end{eqnarray}

As shown in \cite{DingJiang2010AnnalsApplProb}, the empirical eigenvalue distribution for the graph Laplacian $L$,
$$\tilde{F}_n(x)=\frac{1}{n}\sum_{i=1}^n I\left\{ \frac{\lambda_i+n\mu_A}
{\sqrt{n}\sigma_A}\le x \right\},$$ converges weakly (with probability one) as $n\to\infty$ to the free convolution $\gamma$ of the semicircle law,
$\rho_{\text{sc}}(x)=\frac{1}{2\pi}\sqrt{4-x^2}I(|x|\le 2),$
with the standard Gaussian, $g(x)=\exp[-x^2/2]/\sqrt{2\pi}.$
The measure $\gamma$ becomes concentrated around $\lambda_i\approx -n\mu_A$ as $n$ grows.  In particular, most terms in the sum (\ref{eq:define-S}) concentrate around $1/(2n\mu_A)$, as $n\to\infty$. Therefore, setting $C=2\mu_A$, we have
$E[S]\to1/(nC)$, as $n\to\infty$, yielding in the limit
\begin{eqnarray}
E[R_k |M] &=& \frac{\sigma_k^2}{nC} E\left[\sum_{i=2}^n\sum_{j=2}^n(M^\intercal v_i v_i^\intercal \zeta_k)(\zeta_k^\intercal v_j v_j^\intercal M) \right]+O\left(\frac{1}{n^2}\right).
\end{eqnarray}

\begin{figure}[ht!] 
\centering
\includegraphics[width=4in]{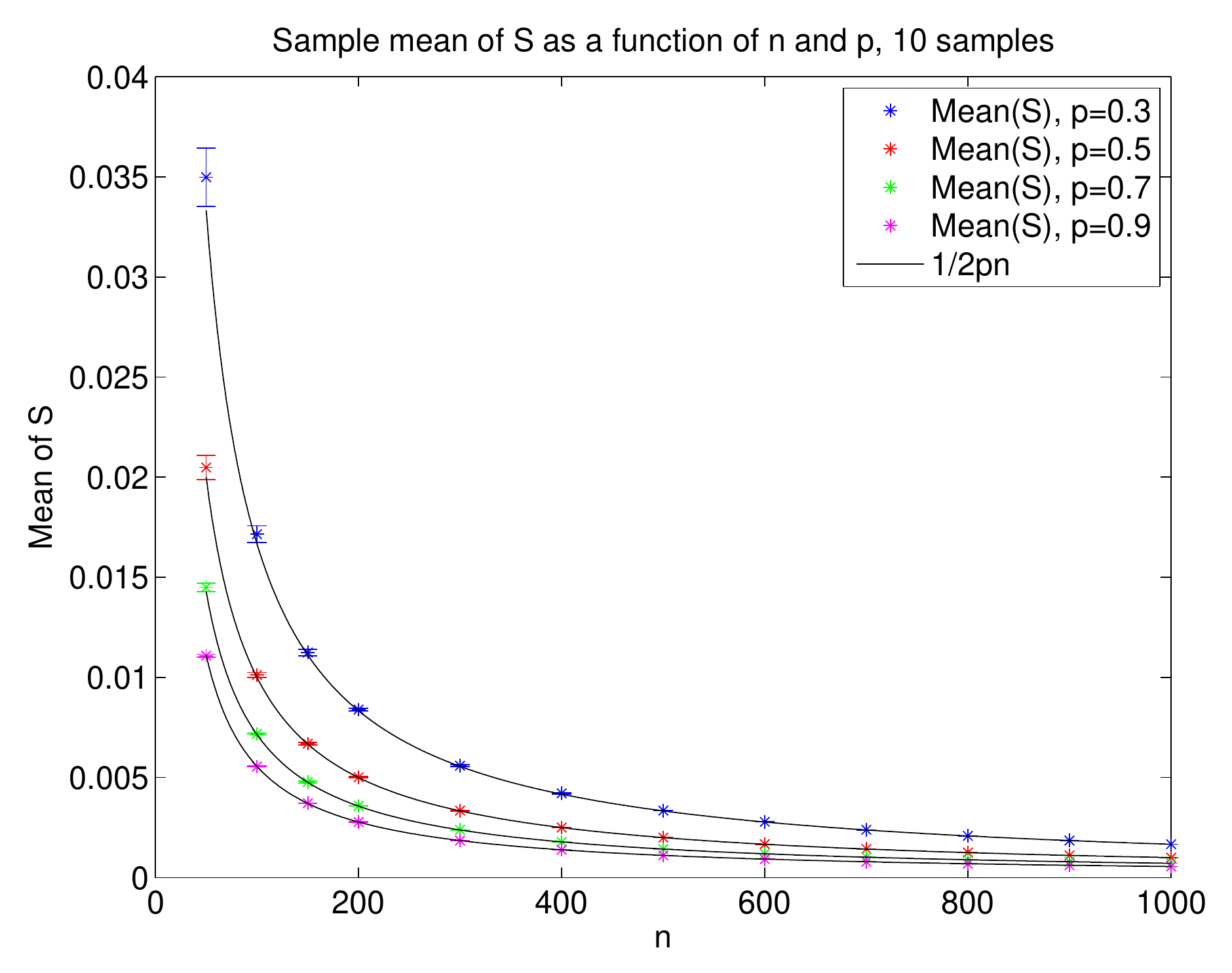}
\caption{Numerical evidence showing that the mean of $S = \frac{1}{(n-1)^2} \sum_{i=2}^n\sum_{j=2}^n \left(\frac{-1}{\lambda_i+\lambda_j}\right)$ is well approximated by $1/(2pn)$, for the \ER ensemble with $p\ge 0.3$.  For a given value of $p$, the colored asterisks show the sample mean of $S$ as a function of $n$ over 10 realizations (with error bars showing the standard deviation) and the black curve is $1/(2pn)$.}
\label{fig:meanS}
\end{figure}

For the \ER ensemble with $n$ nodes and edge probability $p$, we have $E[S]\to1/(nC)$ for $C=2p$.  Figure \ref{fig:meanS} shows that the sample mean of $S$ over 10 realizations (i.e.~10 different \ER random graph configurations with the same parameters) rapidly approaches $1/(2pn)$, as $n$ increases, for values of $p$ ranging from 0.3 to 0.9. As the factor of $1/n$ is common across all $k$,  it does not affect the stochastic shielding argument.  

To prove Theorem \ref{thm:MainResult}, we will show that
\begin{equation}
E\left[ \sum_{i=2}^n\sum_{j=2}^n (M^\intercal v_i v_i^\intercal \zeta_k) (\zeta_k^\intercal v_j v_j^\intercal M) \right] = \left\{
\begin{array}{rr}
1+O(n^{1-q}),& \left|M^\intercal \zeta_k\right|=1\\
O(n^{1-q}),& \left|M^\intercal \zeta_k\right|=0
\end{array}
\right., \text{ as }n\to\infty
\end{equation}
for some $q>1$, corresponding to the parameter $q$ appearing in Assumption A4.  This dichotomy is the basis for neglecting the edges $k$ such that $M^\intercal \zeta_k=0$, as in the stochastic shielding approximation.  To do this we show the following:
\begin{align}
&E\left[ \sum_{i=2}^n\sum_{j=2}^n (M^\intercal v_i v_i^\intercal \zeta_k) (\zeta_k^\intercal v_j v_j^\intercal M) \right]\\
&=\sum_{i=2}^n\sum_{j\ne i}E\left[ (M^\intercal v_i v_i^\intercal \zeta_k)(M^\intercal v_j v_j^\intercal \zeta_k)\right] + \sum_{i=2}^n E[(M^\intercal v_i v_i^\intercal \zeta_k)^2]\label{eq:Rk_not_factored}
\\
&=\sum_{i=2}^n\sum_{j\ne i} E[M^\intercal v_i v_i^\intercal \zeta_k] E[M^\intercal v_j v_j^\intercal \zeta_k] + \sum_{i=2}^n E[(M^\intercal v_i v_i^\intercal \zeta_k)^2] + O\left(\frac{1}{n}\right), \mbox{ as } n\to\infty \label{eq:Rk_factored}
%v_i^\intercal \zeta_k)^2] + O(n^{1-q}), \mbox{ as } n\to\infty \label{eq:Rk_factored}
\end{align}
where the first term is 
\begin{equation}\label{eq:Rk_first_term}
\sum_{i=2}^n\sum_{j\ne i} E[M^\intercal v_i v_i^\intercal \zeta_k] E[M^\intercal v_j v_j^\intercal \zeta_k] =\left|M^\intercal \zeta_k\right|+O\left(\frac{1}{n}\right), \mbox{ as } n\to\infty
\end{equation}
and the second term is 
\begin{equation}\label{eq:Rk_second_term}
\sum_{i=2}^n E[(M^\intercal v_i v_i^\intercal \zeta_k)^2] = O(n^{1-q}), \mbox{ as } n\to\infty.
\end{equation}

Starting with the first term in Equation \ref{eq:Rk_not_factored}, it follows from assumption A3a that, as $n\to\infty$,
\begin{equation}
E\left[ (M^\intercal v_i v_i^\intercal \zeta_k)(M^\intercal v_j v_j^\intercal \zeta_k)\right] = E[ M^\intercal v_i v_i^\intercal \zeta_k] E[M^\intercal v_jv_j^\intercal \zeta_k] + O\left(\frac{1}{n^3}\right)
\end{equation}
which means
\begin{equation}
\sum_{i=2}^n\sum_{j\not=i}E\left[ (M^\intercal v_i v_i^\intercal \zeta_k)(M^\intercal v_j v_j^\intercal \zeta_k)\right] = 
\sum_{i=2}^n\sum_{j\not=i}E[ M^\intercal v_i v_i^\intercal \zeta_k] E[ M^\intercal v_j v_j^\intercal \zeta_k] + O\left(\frac{1}{n}\right).
\end{equation}  

We can expand the left hand side of Equation \ref{eq:Rk_first_term} by using the definitions $M^\intercal v_i = \sum_{l\in 1_M} v_i(l)$ and $v_i^\intercal \zeta_k = v_i(l_+)-v_i(l_-)$ which yield 
\begin{eqnarray}
& & \sum_{i=2}^n\sum_{j\ne i} E[M^\intercal v_i v_i^\intercal \zeta_k] E[M^\intercal v_j v_j^\intercal \zeta_k]\\ 
&=& (n-1)(n-2) E[M^\intercal v_i v_i^\intercal \zeta_k] E[M^\intercal v_j v_j^\intercal \zeta_k]\\
&=& (n-1)(n-2) E\left[\sum_{l\in 1_M} v_i(l) (v_i(l_+)-v_i(l_-))\right]^2. \label{eq:continue}
\end{eqnarray}
By Lemma \ref{lem:random} part B, we have that $E\left[\sum_{l\in 1_M} v_i(l) (v_i(l_+)-v_i(l_-))\right]^2=\frac{1}{n^2}|M^\intercal \zeta_k| + O(n^{-4})$, as $n\to\infty$.  Continuing Equation \ref{eq:continue} above we have
\begin{eqnarray}
&=& (n-1)(n-2)\left[ \frac{1}{n^2}|M^\intercal \zeta_k| + O(n^{-4}) \right]\\
&=&  |M^\intercal \zeta_k| + O(n^{-1})
\end{eqnarray}
as $n\to\infty$, which establishes the first term (Equation \ref{eq:Rk_first_term}).

We now focus on the second term in Equation \ref{eq:Rk_factored}.  In Lemma \ref{lem:random} part C, we establish that as $n\to\infty$ 
\begin{equation}
E[(M^\intercal v_i v_i^\intercal \zeta_k)^2] = E\left[\left(\sum_{l\in 1_M} v_i(l) \right)^2 (v_i(l_+)-v_i(l_-))^2 \right ]=O(n^{-q}).
\end{equation}
Hence, $(n-1) E\left[\left(\sum_{l\in 1_M} v_i(l) \right)^2 (v_i(l_+)-v_i(l_-))^2 \right ] = O(n^{1-q})$
as $n\to\infty$ which establishes the second term (Equation \ref{eq:Rk_second_term}).  Therefore, we have established Theorem \ref{thm:MainResult}.

%%%%%%%%%%%%%%%%%%%%%%%%%%%%%%%
% Figure 5
%%%%%%%%%%%%%%%%%%%%%%%%%%%%%%%

\subsection{Symmetric \ER Random Graph Ensemble} 
\label{ssec:ER}

\begin{figure}[ht!] 
\centering
\includegraphics[width=4in]{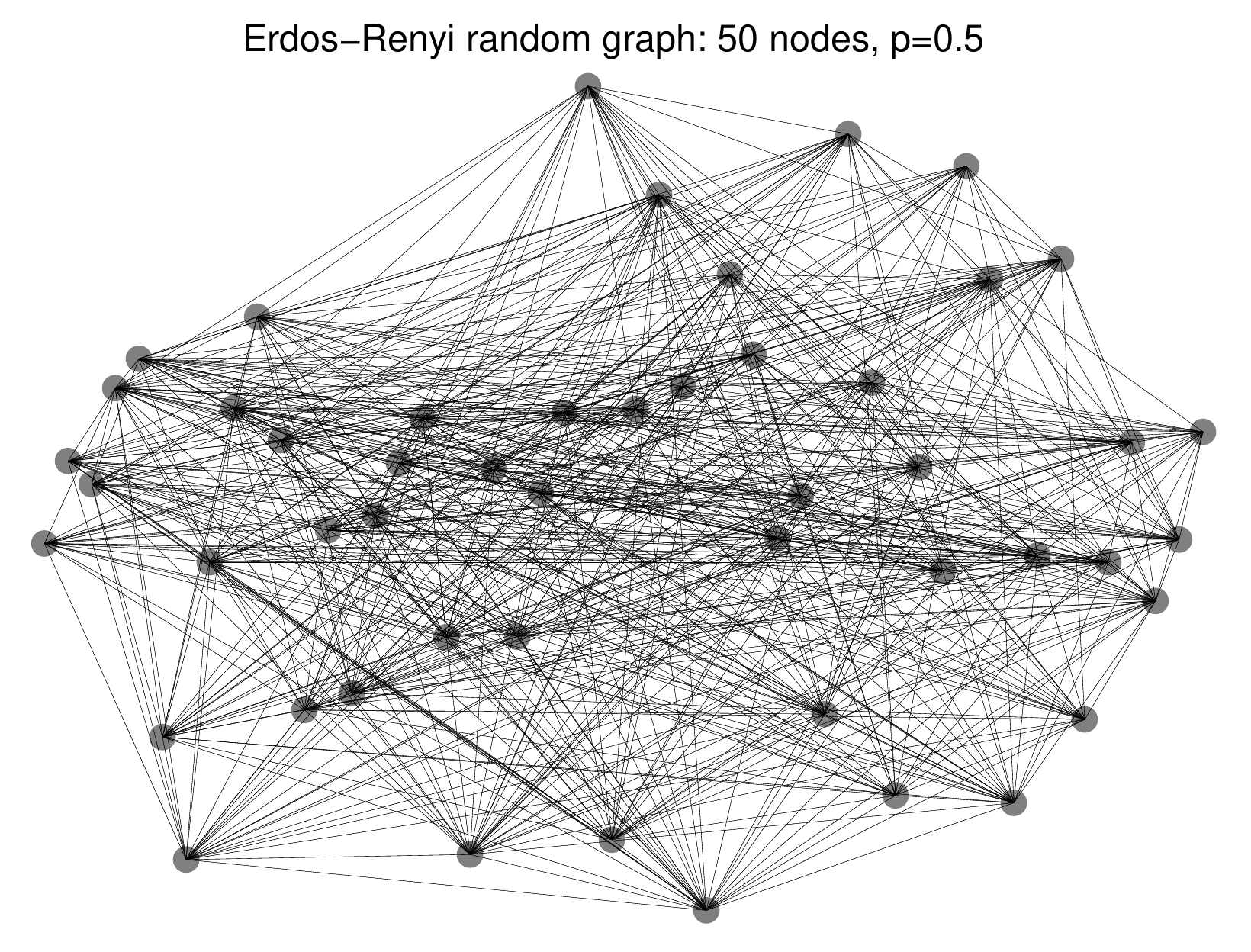}
\caption{Realization of an \ER random graph with $n=50$ nodes and edge probability $p=0.5$.}
\label{fig:ERgraph_50}
\end{figure}

Many varieties of random graphs have been used to describe biological systems \cite{AlbertBarabasi2002RMP,Durrett2007}. Here, we restrict attention to an ensemble of symmetric irreducible \ER random graphs $\mathcal{G}(n,p)$ on $n$ nodes, for which each of $(n^2-n)/2$ possible bidirectional edges occurs independently with probability $p$ \cite{ErdosRenyi1959,ErdosRenyi1960}.  The exclusion of reducible graphs does not represent a significant restriction because the \ER ensemble generates irreducible (connected) graphs with high probability for large $n$ and $p > \frac{ln(n)}{n}$ \cite{ErdosRenyi1959}.  In particular, consider an \ER random graph ensemble for $n=50$ and $p=0.5$.  See Figure \ref{fig:ERgraph_50} for a network drawn from this ensemble.  Take $A$ to be the unweighted adjacency matrix ($\alpha_k\in\{0,1\}$) and let $\sigma_k=1$ for all reactions $k$ so that the $k^{th}$ column of noise matrix $B$ is exactly the stoichiometry vector for reaction $k$.  Specifying any measurement vector $M \in\{0,1\}^{50}$ induces a partition of edges into ``important" (type 0-1) or ``unimportant" (types 0-0 or 1-1) classes.  Let $\mathcal{E}_I$ be the set of important edges and $\mathcal{E}_U$ be the set of unimportant edges.  Clearly, $\mathcal{E}=\mathcal{E}_I \cup \mathcal{E}_U$.  In the following example, we consider a vector $M$ such that half the entries are 1 and other half are 0.  

%%%%%%%%%%%%%%%%%%%%%%%%%%%%%%%
% Figure 6
%%%%%%%%%%%%%%%%%%%%%%%%%%%%%%%

Theorem \ref{thm:MainResult} says that if the matrix of eigenvector components of the \ER graph Laplacian is sufficiently similar to a random matrix drawn from the Gaussian ensemble (in terms of assumptions A0-A5) then one would expect the partitioning of the $R_k$ into two clusters.  One cluster, containing the important edges, will be centered at $1/n$.  A second cluster, containing the unimportant edges, will have smaller $R_k$ values ($O(n^{-q})$ where $q>1$ is governed by the fourth moment, see assumption A4a in \S \ref{ssec:assumptions}).  To the extent to which this similarity to the Gaussian ensemble holds, our calculation of $R_k$ involves projecting the measurement vector $M$ and the vectors $\zeta_k$ onto randomly chosen subspaces of $\R^n$. 

\begin{figure}[ht] 
\centering
	\includegraphics[width=4in]{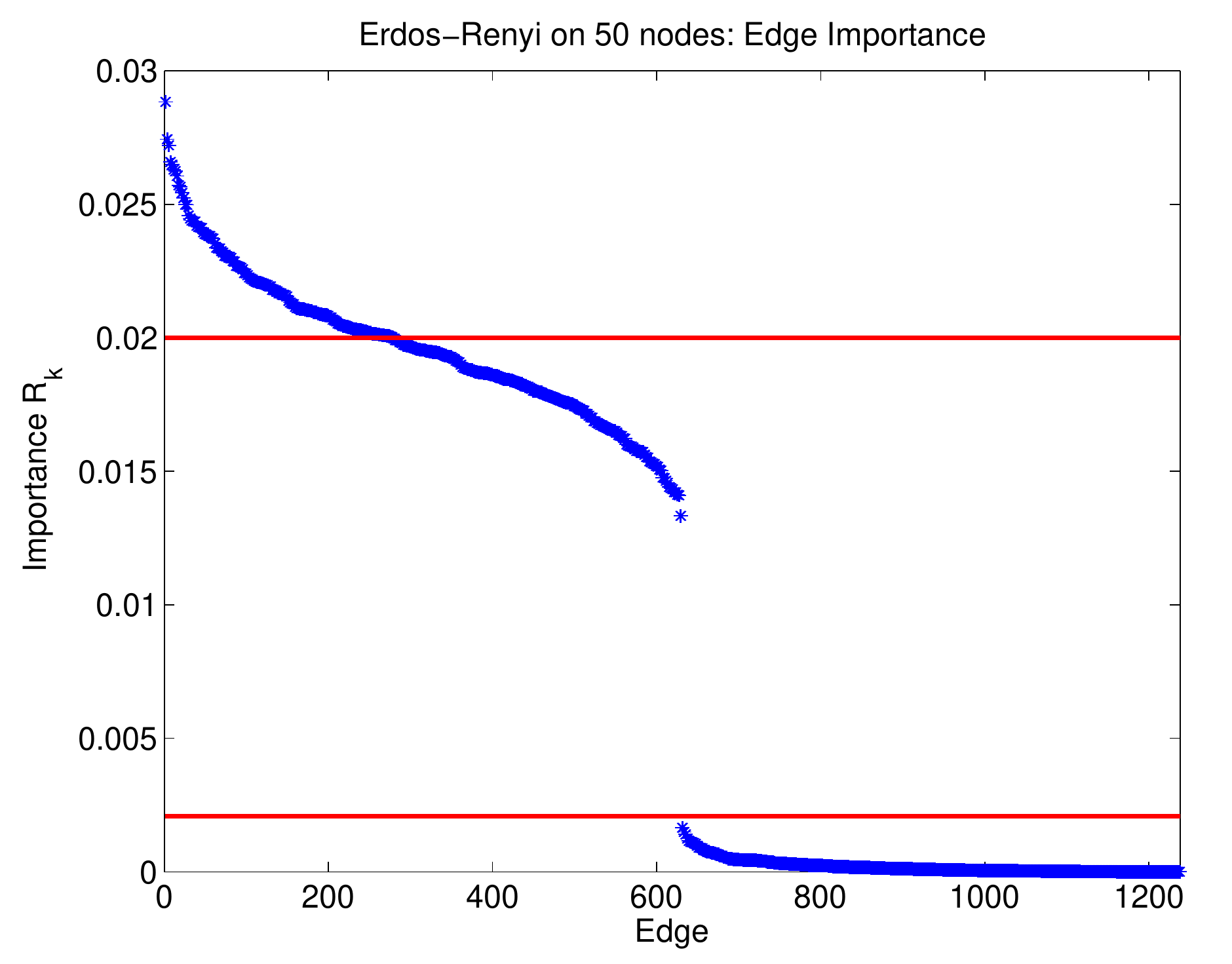}
\caption{Edge importance values $R_k$ plotted in descending order for the process on an \ER random graph with 50 nodes, edge probability 0.5, and measurement functional $M$ such that half the nodes are labeled 1 and the other half are 0.  There is a clear separation between the important edges (type 0-1) and unimportant edges (types 0-0 and 1-1).  The cluster of important edges has a mean $R_k$ value of $1/50=0.02$ whereas the  unimportant cluster lies below the line at $\sqrt{2}n^{-5/3} \approx 0.0021$.}
\label{fig:ER_edge_importance50}
\end{figure}

As shown in Figure \ref{fig:covariance_all}, assumptions A0-A5 appear to be satisfied for the symmetric \ER random graph ensemble.  In particular, the fourth moment of the eigenvector components (assumption A4a) appears to hold empirically for $q \approx 5/3$; in particular, we find that, empirically, $E[v_i(l)^4] \approx \sqrt{2}n^{-5/3}$.  This behavior suggests that the unimportant edges should have a mean $R_k$ value $\lesssim \sqrt{2}n^{-5/3}$.  Setting $n=50$, for example, we would expect one cluster of $R_k$ values centered at $1/50=0.02$ for $k\in\mathcal{E}_I$ and another cluster close to $\sqrt{2} \cdot 50^{-5/3}\approx 0.0021$ for $k\in\mathcal{E}_U$.  Figure \ref{fig:ER_edge_importance50} shows the rank order of edge importance values $R_k$ corresponding to the $m$ reactions in the \ER random graph.  The top cluster is centered at $0.02$ (upper horizontal red line) and the bottom cluster is bounded above by $0.0021$ (lower horizontal red line) consistent with Theorem \ref{thm:MainResult} for the \ER random graph ensemble with 50 nodes and edge probability $p=0.5$.  Since the measurement functional $M$ is binary, we see a significant gap between the two clusters, as expected.  If the components of $M$ are graded, i.e.~drawn uniformly from the unit interval, then this curve appears to be smooth (see discussion in \S \ref{sec:discuss}).

%%%%%%%%%%%%%%%%%%%%%%%%%%%%%%%
% Figure 7
%%%%%%%%%%%%%%%%%%%%%%%%%%%%%%% 

\begin{figure}[ht] 
\centering
	\subfigure{
	\includegraphics[width=3in]{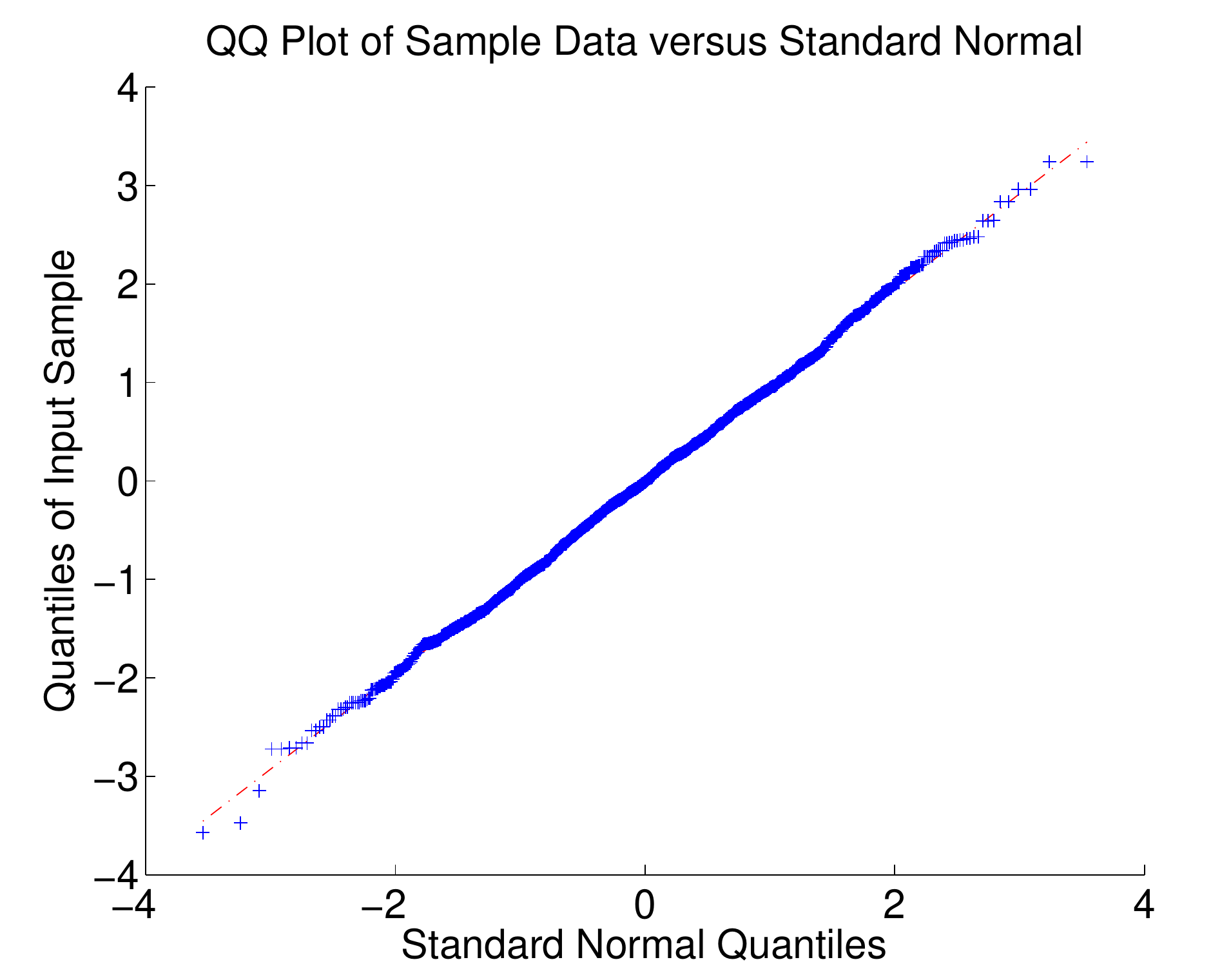}
	}
	\subfigure{
	\includegraphics[width=3in]{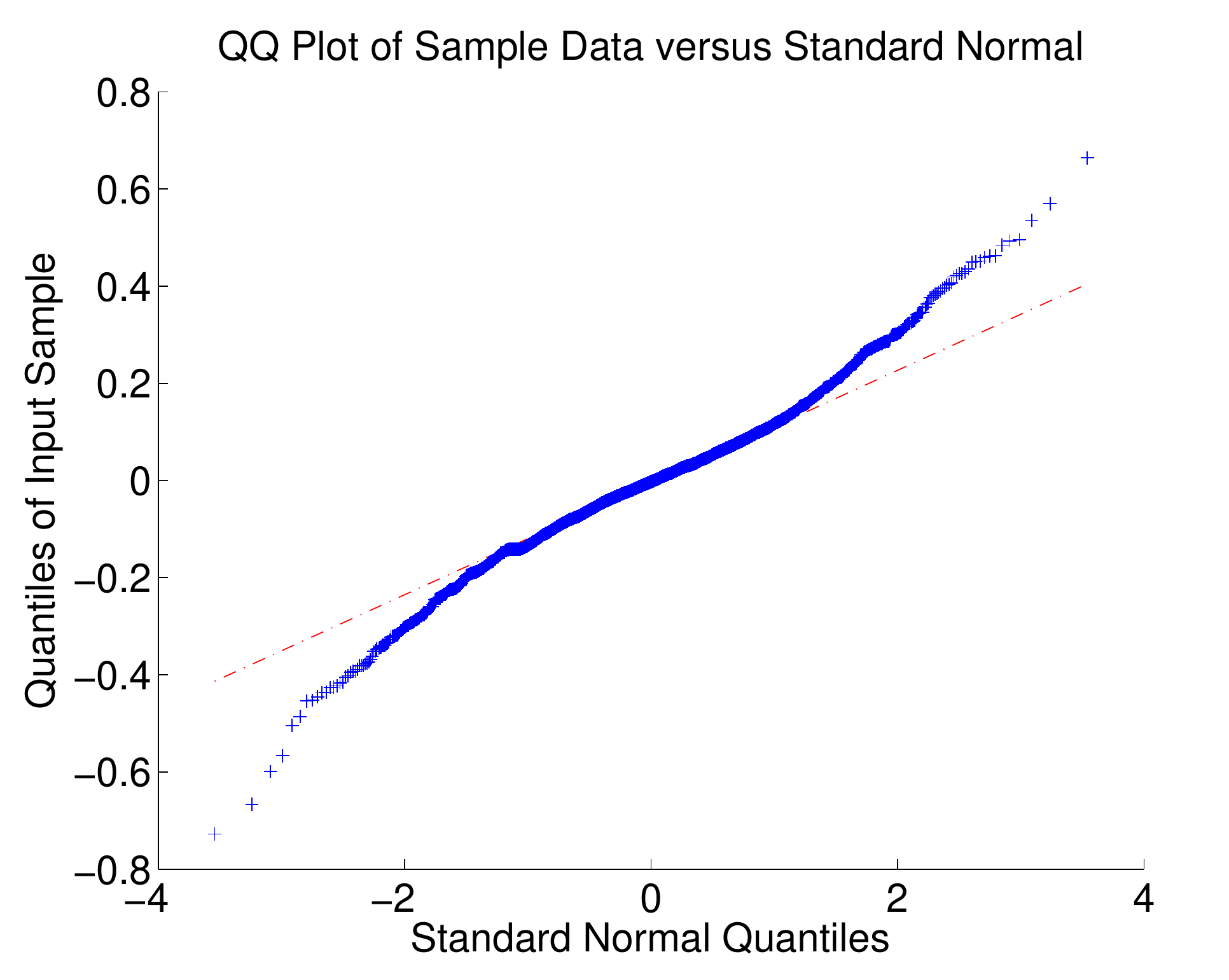}
	}
\caption{Numerical evidence illustrating that the eigenvector components of the graph Laplacian for the symmetric \ER random graph ensemble are close to Gaussian distributed (to one standard deviation).  Left: quantile-quantile plot for a Gaussian random matrix with $\mathcal{N}(0,1/50)$ entries.  Right: quantile-quantile plot of eigenvector components for the \ER case with $n=50$ nodes and edge probability $p=0.5$.}
\label{fig:qqplots}
\end{figure}

Figure \ref{fig:qqplots} illustrates the distribution of eigenvector components of the \ER graph Laplacian in comparison with a Gaussian random matrix (i.e., each entry has mean 0 and variance $1/n$).  The quantile-quantile plots show good agreement within one standard deviation and begin to deviate in the second standard deviation.  This is consistent with the observation that the fourth moment in the \ER case deviates from the Gaussian case ($q\approx 5/3$ for \ER and $q=2$ for Gaussian).  
%Hence, the eigenvector components of the symmetric \ER graph Laplacian $L$ are approximately Gaussian distributed with mean 0 and variance $1/n$.
Nevertheless, Theorem \ref{thm:MainResult} predicts that there will be two clusters of $R_k$ values as described above and shown in Figure \ref{fig:ER_edge_importance50} for the \ER case with $n=50$ and $p=0.5$.

%%%%%%%%%%%%%%%%%%%%%%%%%%%%%%%
% Figure 8
%%%%%%%%%%%%%%%%%%%%%%%%%%%%%%% 

%\newpage
\section{Application: Stochastic Shielding of Hodgkin-Huxley Channels under Voltage Clamp}
\label{sec:HH}

Hodgkin and Huxley's (HH) model for the generation and propagation of action potentials along the giant axon of the squid \textit{Loligo} lies at the foundations of modern neuroscience \cite{Catterall2012JNeurosci,jp:Hodgkin+Huxley:1952d}.   In the classic HH model, action potentials are generated through the interaction of a leak current and two voltage-gated ionic currents, carried by a sodium ion specific channel and a potassium ion specific channel.    The potassium channel comprises four identical subunits that open and close independently with voltage-dependent rates.  The channel carries a current when all four subunits are in the open state.  At the molecular level, a single channel can be represented as a continuous time Markov jump process on a chain of five states, the fifth of which has nonzero conductance.  Of the eight transitions connecting states along this chain, only the last two connect states with different conductances, therefore the stochastic shielding approximation would preserve the fluctuations of these transitions and not the other six.  

The sodium channel involves two types of subunits, an activation subunit (``$m$'')  present in three identical copies, and an inactivation subunit (``$h$'') present in a single copy.\footnote{Modern measurements of purified sodium channel preparations suggest the presence of four activation gates \cite{Catterall2012JPhysiol}; for consistency with common usage we will restrict attention here to the classical $n^4$ potassium channel and $m^3h$ sodium channel formulations of the model.}  The resulting graph has eight distinct states connected by twenty different transitions, each occurring with a voltage-dependent rate \cite{GoldwynShea-Brown2011PLoSCB,SchneidmanFreedmanSegev1998NeuralComput,SkaugenWalloe1979ActaPhysiolScand}.   Four of these twenty transitions connect states with differing conductance values (zero \textit{versus} nonzero); the fluctuations of the remaining sixteen transitions are ignored under the stochastic shielding approximation.  

Schmandt and \galan compared simulations of a system comprising 5,000 individual potassium channels and 25,000 individual sodium channels, both with and without the stochastic shielding approximation.  It is possible to construct an exact simulation scheme, analogous to Gillespie's stochastic simulation algorithm \cite{Gillespie1977}, that takes into account the nonstationarity of the transition rates (propensities) arising from their voltage dependence \cite{ClayDeFelice1983BiophysJ}.  However, Schmandt and \galan used a discrete time approximation to this process.  Appendix \ref{app:stoch-shield-GS} discusses Schmandt and Gal\'{a}n's approach in more detail.  Here we apply our analysis to evaluate the edge importance $R_k$ of each transition in the graph for the classic HH potassium and sodium channels, respectively.  Rather than consider the case of time-varying transition rates, we restrict attention to the ``voltage clamped" case.  If the membrane potential is experimentally held constant for a given cell, the \textit{per capita} transition rates remain constant and the fluctuating ion channel population forms a stationary Markov process.  In particular, our analysis approximates this stationary population process with a linear multidimensional Ornstein-Uhlenbeck process (see Appendix \ref{app:tau-leaping}); this approximation is reasonable given the large numbers of individual channels considered in Schmandt and Gal\'{a}n's simulations.

In general, the ion channel state graphs for the potassium and sodium channels in the HH model have graph Laplacians $L$ that are not symmetric.  Therefore, we need to modify our definition of the edge importance $R_k$ (Equation \ref{eq:error_Rk}) in order to apply our results.
When $L$ is not symmetric, we will assume that $L$ is nevertheless diagonalizable, i.e.~that there are eigenvalues $\lambda_i$ and a biorthogonal system of vectors $v_i, w_i$ (right and left eigenvectors) satisfying 
\begin{eqnarray}\nonumber
Lv_i&=&\lambda_i v_i\\
w_i^\intercal L&=& \lambda_i w_i^\intercal\\
w_i^\intercal v_j&=&\delta_{ij}.\nonumber
\end{eqnarray}
In this case the decomposition of $L$ becomes $L=\sum_i\lambda_iv_iw_i^\intercal$, and the definition of $R_k$ is modified as follows: 
\begin{equation}
R_k=\sigma_k^2 \sum_{i=2}^n\sum_{j=2}^n
\left(\frac{-1}{\lambda_i+\lambda_j}\right) (M^\intercal v_i) (w_i^\intercal \zeta_k) (\zeta_k^\intercal w_j) (v_j^\intercal M).
\end{equation}

\subsection{Hodgkin-Huxley Potassium Channel}
\label{ssec:HHK}

\begin{figure}[ht!] 
\centering
	\includegraphics[width=4in]{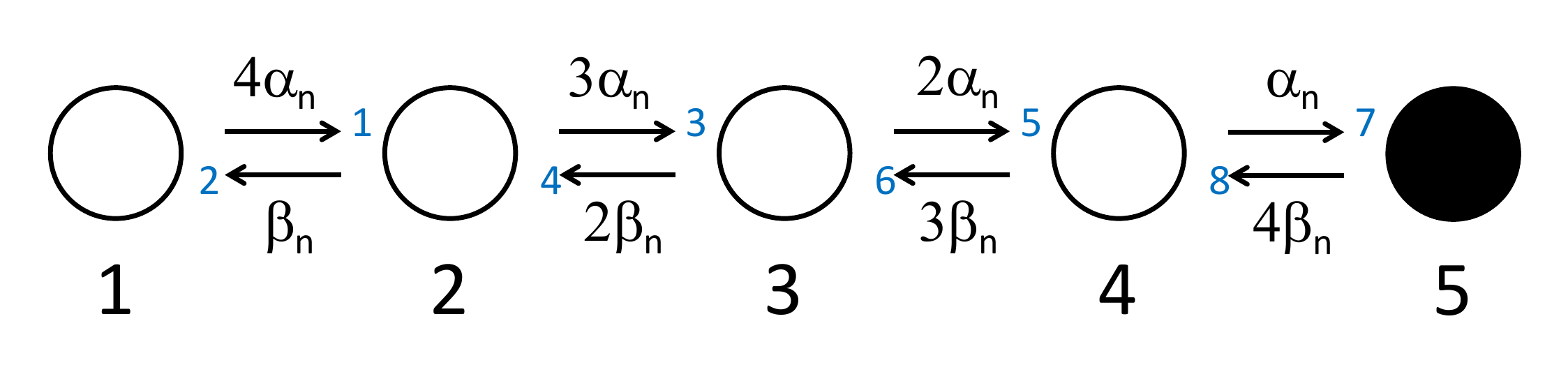}
\caption{An illustration of the Hodgkin-Huxley potassium channel state graph.  This is a 5-state chain where state 5 is the conducting state.  The eight reactions are labeled in blue and are used to define the edge importance values $R_k$ in the figures below.  The reaction rates $\alpha_n$ and $\beta_n$ are voltage-dependent as defined by Equations \ref{eq:HH-K-voltage-alpha}-\ref{eq:HH-K-voltage-beta}.}
\label{fig:HH-K-channel}
\end{figure}

The potassium channel state graph in the Hodgkin-Huxley model is a 5-state chain with one conducting state.  Following the tau-leaping construction (Appendix \ref{app:tau-leaping}) we consider a stationary OU process $X(t)\in\R^5$, with linear measurement functional $M=[0,0,0,0,1]^\intercal$.  See Figure \ref{fig:HH-K-channel} for an illustration of this channel. 
%%%%%%%%%%%%%%%%%%%%%%%%%%%%%%%
% Figure 9
%%%%%%%%%%%%%%%%%%%%%%%%%%%%%%% 
The corresponding (weighted) adjacency matrix $A$ is 
\begin{eqnarray}\label{eq:HH-K-adjacency}
A&=&\matrix{ccccc}
{0 & 4\alpha_n(V) & 0 & 0 & 0\\
\beta_n(V) & 0 & 3\alpha_n(V) & 0 & 0\\
0 & 2\beta_n(V) & 0 & 2\alpha_n(V) & 0\\
0 & 0 & 3\beta_n(V) & 0 & \alpha_n(V)\\
0 & 0 & 0 & 4\beta_n(V) & 0}
\end{eqnarray}
which is evidently not symmetric.  The voltage-dependent transition rates are given by:
\begin{eqnarray}
\alpha_n(V) &=& \frac{0.01 (V+55)}{1-e^{(-0.1(V+55))}}\label{eq:HH-K-voltage-alpha}\\
\beta_n(V) &=& 0.125 e^{\frac{-(V+65)}{80}}.\label{eq:HH-K-voltage-beta}
\end{eqnarray} 
Then the graph Laplacian $L=(A-D)^\intercal$ is voltage-dependent and is given by
\begin{eqnarray*}%\label{eq:HH-K-Laplacian}
L&=&\matrix{ccccc}
{-4\alpha_n(V) & \beta_n(V) & 0 & 0 & 0\\
4\alpha_n(V) & -(\beta_n(V)+3\alpha_n(V)) & 2\beta_n(V) & 0 & 0\\
0 & 3\alpha_n(V) & -2(\beta_n(V)+\alpha_n(V)) & 3\beta_n(V) & 0\\
0 & 0 & 2\alpha_n(V) & -(3\beta_n(V)+\alpha_n(V)) & 4\beta_n(V)\\
0 & 0 & 0 & \alpha_n(V) & -4\beta_n(V)}
\end{eqnarray*}
since the entries in the diagonal matrix $D$ are the weighted out-degrees of each node for a given voltage $V$, i.e. $D_{ii}(V) = \sum_{j=1}^5 A_{ij}(V)$.  The noise matrix $B$ is also voltage-dependent.  Recall that the $k^{th}$ column of $B$ corresponds to the $k^{th}$ reaction, and this can be written as $\sigma_k(V)\zeta_k$.  If $r_k$ is the \textit{per capita} rate of reaction $k$ (transition from node $i(k)$ to $j(k)$), then $\sigma_k(V)=\sqrt{r_k(V) \bar{N}_i(V)}$ where $\bar{N}_i(V)$ is the average number of channels at state $i$ at equilibrium for voltage $V$. 
Hence, $B$ is given by 
\begin{equation}\label{eq:generalB}
B = \left( \sqrt{r_1(V)\bar{N}_{i(1)}(V)}\zeta_1, \dots, \sqrt{r_k(V)\bar{N}_{i(k)}(V)}\zeta_k, \dots, \sqrt{r_m(V)\bar{N}_{i(m)}(V)}\zeta_m \right).
\end{equation}

\begin{figure}[ht!] 
\centering
	\includegraphics[width=4in]{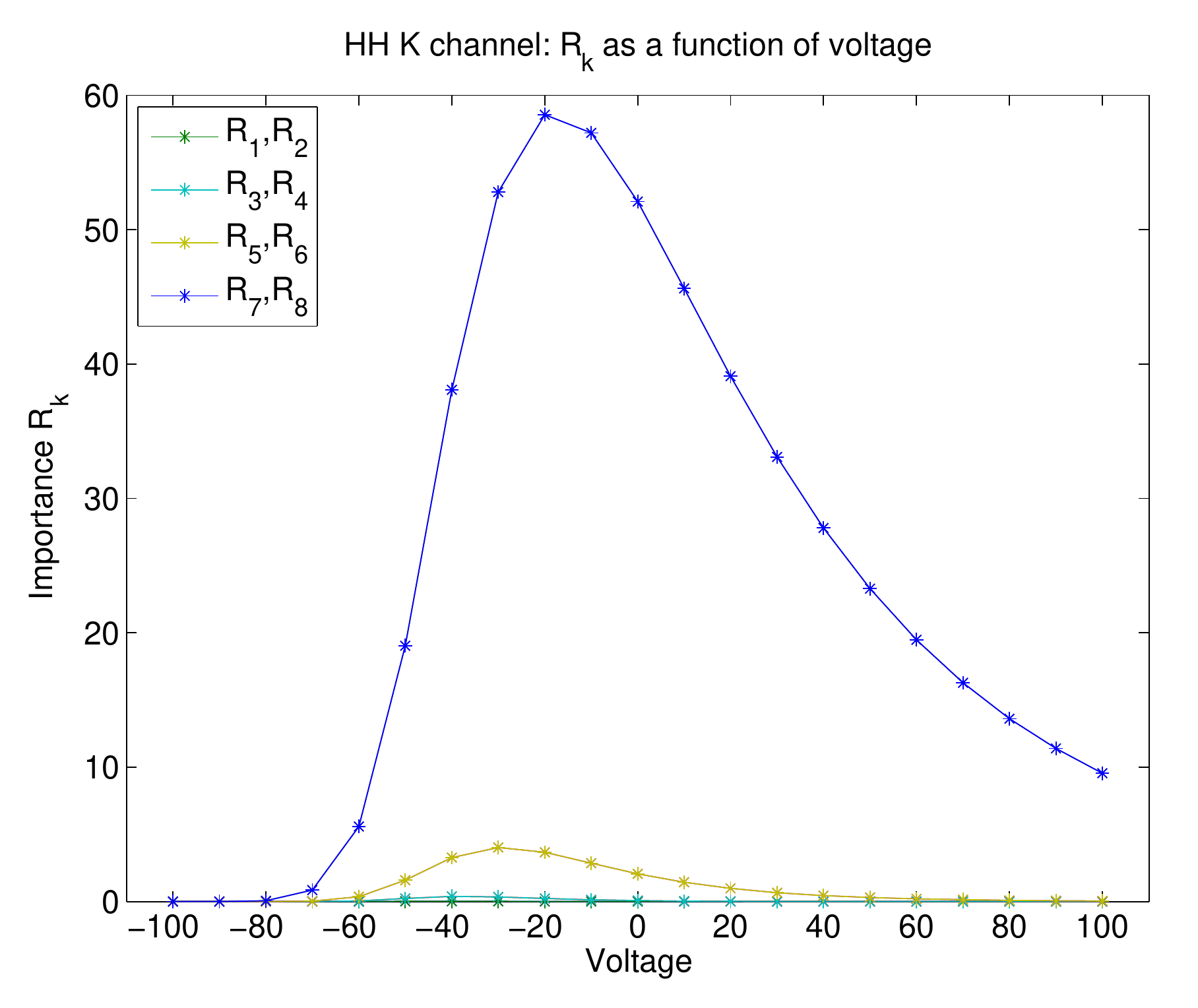}
\caption{Hodgkin-Huxley potassium channel.  This figure shows edge importance $R_k$ as a function of voltage in the range $[-100,100]$ mV for each reaction $k\in\{1,\dots,8\}$.  The blue curve corresponds to edges 7 and 8 ($R_7=R_8$), the transitions between state 4 and conducting state 5, which is the largest $R_k$ value in the voltage range above.  If neglected, these two reactions would have the highest contribution to the error.} 
\label{fig:HH-Kchannel-Rk-voltage}
\end{figure}

Figure \ref{fig:HH-Kchannel-Rk-voltage} shows the edge importance $R_k$ as a function of voltage for each reaction $k\in\{1,\dots,8\}$ in the potassium channel state graph.  Note that since the process is at steady state, and respects detailed balance, the mean flux due to the two reactions connecting the same pair of nodes will be equal and opposite.  Thus, in this case, $R_1=R_2$, $R_3=R_4$, $R_5=R_6$, and $R_7=R_8$.  The blue curve ($R_7=R_8$) corresponds to edges 7 and 8, the transitions between state 4 and conducting state 5, and has the largest edge importance value in the voltage range $[-100,100]$ mV.  This says that if either or both of these reactions are neglected, they would have the highest contribution to the error.

%%%%%%%%%%%%%%%%%%%%%%%%%%%%%%%
% Figure 10
%%%%%%%%%%%%%%%%%%%%%%%%%%%%%%% 

Physically, it is the current rather than the state occupancy that holds the greatest interest. The current through a population of potassium channels with net conductance $g$ is $I=g(V-V_k)$; here $V_k=-77$ mV is the potassium reversal potential, and the conductance $g=g^oN_o$ is the product of the unitary or single channel conductance $g^o$ with the total number of channels in the open state, $N_o$.  The variance of the current is therefore $(g^o(V-V_k))^2$ times the variance of the occupancy number, meaning that near the reversal potential, the current can have low variance even if the channel state has high variance.  For convenience we set $g^o=1$, which amounts to a change of nominal units for measuring the conductance.

\begin{figure}[ht!] 
\centering
	\subfigure{
	\includegraphics[width=3.05in]{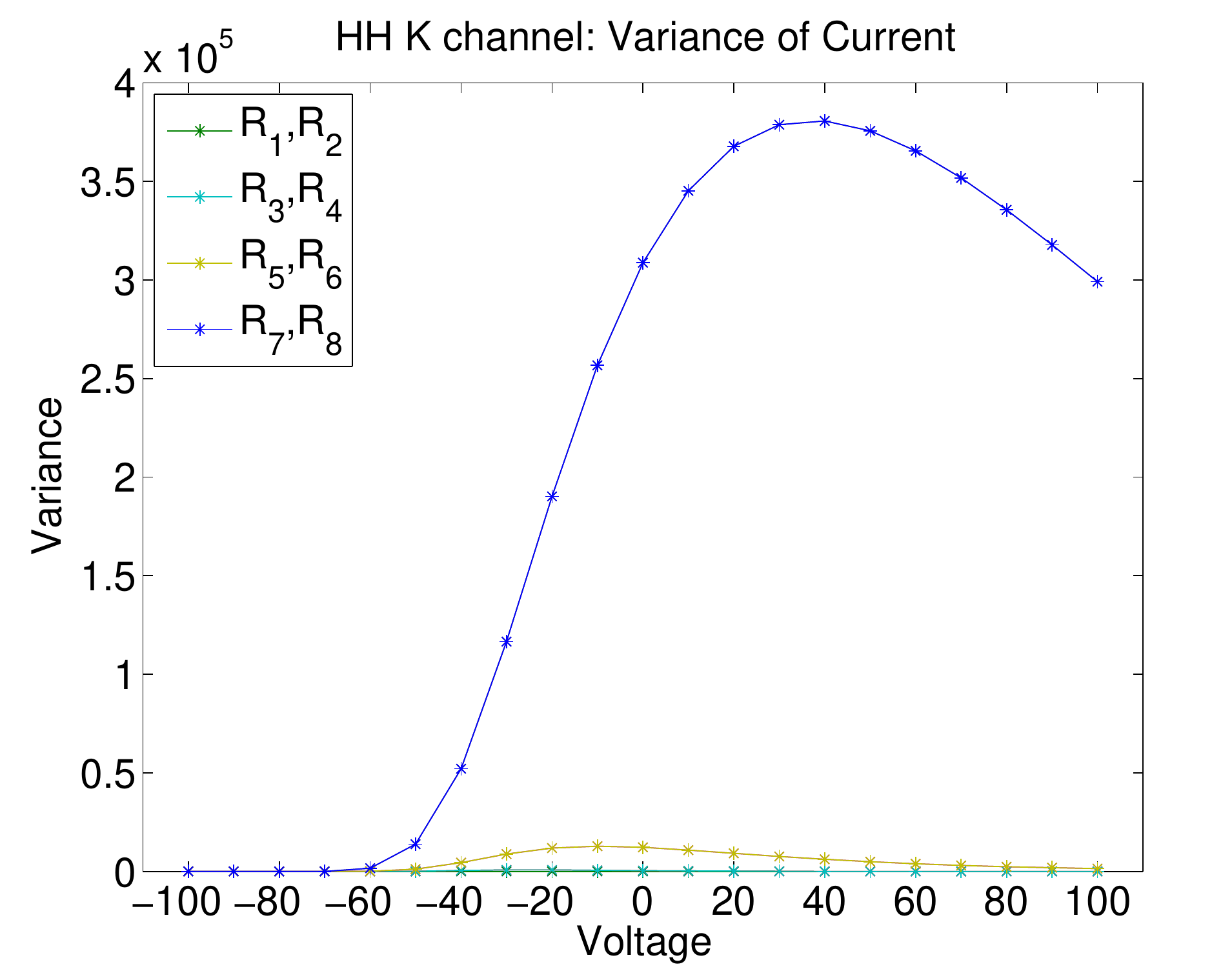}
	}
	\subfigure{
	\includegraphics[width=3in]{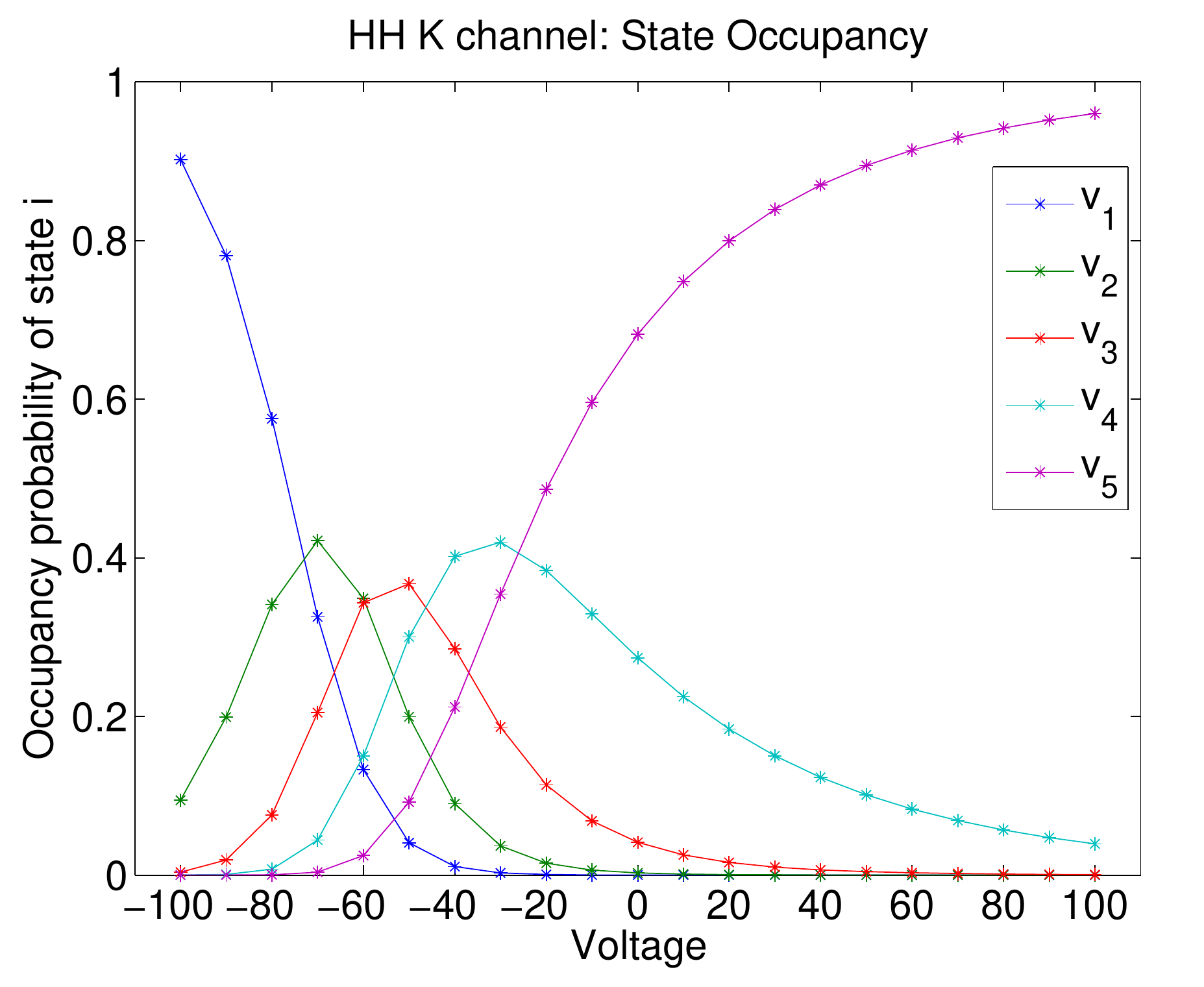}
	}
\caption{Hodgkin-Huxley potassium channel. Left: variance of the current $R_k * (V - V_k)^2$ as a function of voltage $V$ for each reaction $k$ where $V_k=-77$ mV is the reversal potential for the potassium channel.  The blue curve $R_7=R_8$ has the highest variance. Right: leading eigenvector components (normalized so that the components sum to 1) as a function of voltage.}
\label{fig:HH-Kchannel-variance-current-state-occupancy}
\end{figure}

Figure \ref{fig:HH-Kchannel-variance-current-state-occupancy} shows the variance of the nominal current, $R_k * (V - V_k)^2$ as a function of voltage $V$ for each reaction $k$ for the potassium channel.  In addition to having the highest edge importance curve, the blue curve $R_7=R_8$ also has the highest variance (left panel).  The right panel shows the probability of being in each state as a function of voltage.
%%%%%%%%%%%%%%%%%%%%%%%%%%%%%%%
% Figure 11
%%%%%%%%%%%%%%%%%%%%%%%%%%%%%%% 

\subsection{Hodgkin-Huxley Sodium Channel}
\label{ssec:HHNa}

\begin{figure}[ht!] 
\centering
	\includegraphics[width=4in]{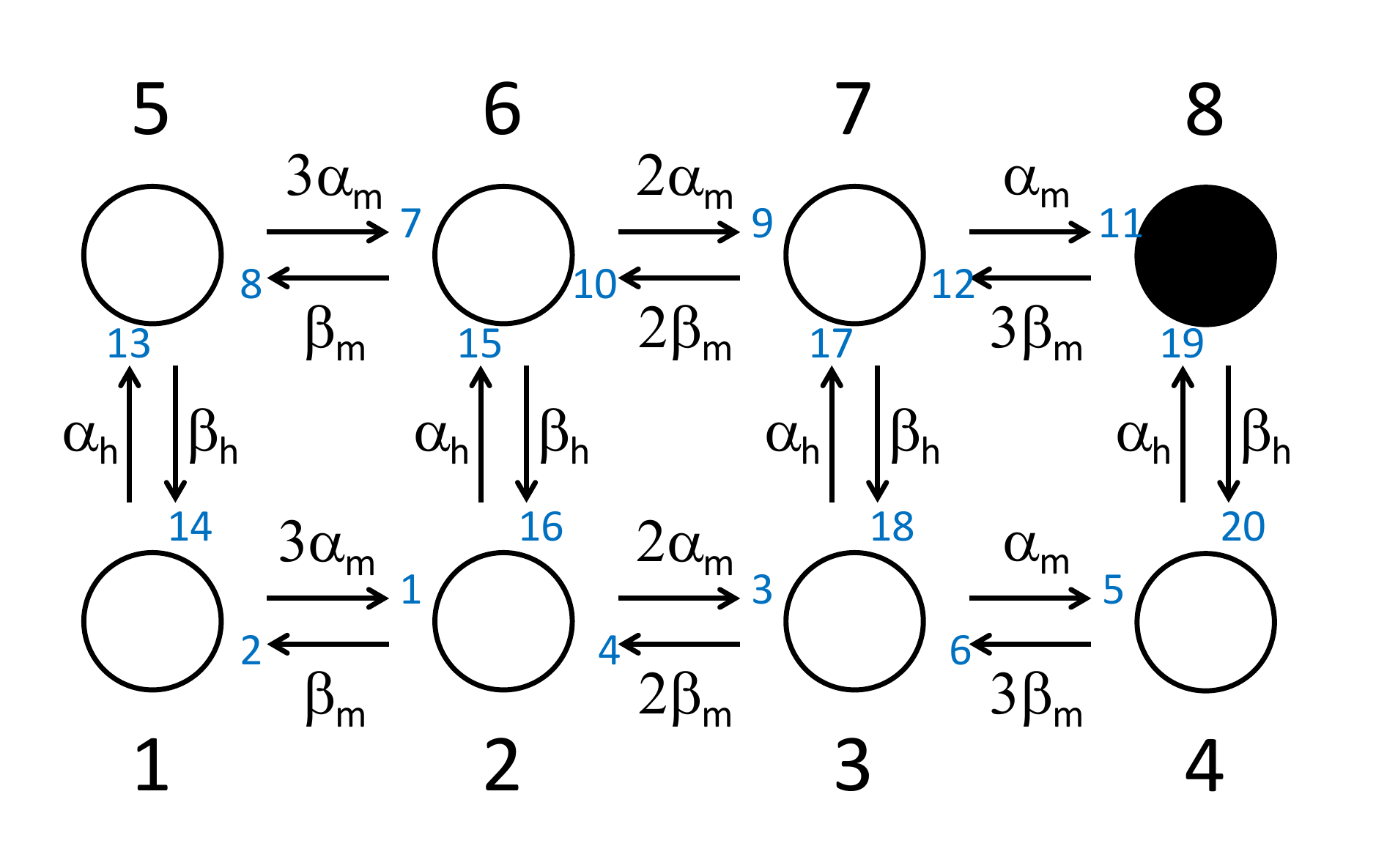}
\caption{An illustration of the Hodgkin-Huxley sodium channel.  This channel has eight states, where state 8 is the conducting state, and twenty reactions.  The reactions are labeled in blue and are used to define the edge importance values $R_k$ in the figures below.  The reaction rates $\alpha_m$, $\alpha_h$, $\beta_m$, and $\beta_h$ are voltage-dependent, defined in Equations \ref{eq:HH-Na-voltage-alpham}-\ref{eq:HH-Na-voltage-betah}.}
\label{fig:HH-Na-channel}
\end{figure}

The sodium channel state graph in the Hodgkin-Huxley model consists of two linked 4-state chains, for a total of eight states, including one conducting state, and twenty reactions.  Again following the tau-leaping construction (Appendix \ref{app:tau-leaping}) we consider a stationary OU process $X(t)\in\R^8$, with linear measurement functional $M=[0,0,0,0,0,0,0,1]^\intercal$.  See Figure \ref{fig:HH-Na-channel} for an illustration. 

%%%%%%%%%%%%%%%%%%%%%%%%%%%%%%%
% Figure 12
%%%%%%%%%%%%%%%%%%%%%%%%%%%%%%%  
The adjacency matrix in this case is

\begin{eqnarray}\label{eq:HH-Na-adjacency}
A&=&\matrix{cccccccc}
{0 & 3\alpha_m(V) & 0 & 0 & \alpha_h(V) & 0 & 0 & 0\\
\beta_m(V) & 0 & 2\alpha_m(V) & 0 & 0 & \alpha_h(V) & 0 & 0\\
0 & 2\beta_m(V) & 0 & \alpha_m(V) & 0 & 0 & \alpha_h(V) & 0\\
0 & 0 & 3\beta_m(V) & 0 & 0 & 0 & 0 & \alpha_h(V)\\
\beta_h(V) & 0 & 0 & 0 & 0 & 3\alpha_m(V) & 0 & 0\\
0 & \beta_h(V) & 0 & 0 & \beta_m(V) & 0 & 2\alpha_m(V) & 0\\
0 & 0 & \beta_h(V) & 0 & 0 & 2\beta_m(V) & 0 & \alpha_m(V)\\
0 & 0 & 0 & \beta_h(V) & 0 & 0 & 3\beta_m(V) & 0}
\end{eqnarray}
where the voltage-dependent entries are defined by
\begin{eqnarray}
\alpha_m(V) &=& \frac{0.1(V+40)}{1-e^{\frac{-(V+40)}{10}}},\label{eq:HH-Na-voltage-alpham} \quad
\beta_m(V) = 4 e^{\frac{-(V+65)}{18}}\\
\alpha_h(V) &=& 0.07 e^{\frac{-(V+65)}{20}}, \quad
\beta_h(V) = \frac{1}{1+e^{\frac{-(V+35)}{10}}}\label{eq:HH-Na-voltage-betah}.
\end{eqnarray}
The graph Laplacian $L=(A-D)^\intercal$ is
\begin{eqnarray*}%\label{eq:HH-Na-Laplacian}
L&=&\matrix{cccccccc}
{-D_{11}(V) & \beta_m(V) & 0 & 0 & \beta_h(V) & 0 & 0 & 0\\
3\alpha_m(V) & -D_{22}(V) & 2\beta_m(V) & 0 & 0 & \beta_h(V) & 0 & 0\\
0 & 2\alpha_m(V) & -D_{33}(V) & 3\beta_m(V) & 0 & 0 & \beta_h(V) & 0\\
0 & 0 & \alpha_m(V) & -D_{44}(V) & 0 & 0 & 0 & \beta_h(V)\\
\alpha_h(V) & 0 & 0 & 0 & -D_{55}(V) & \beta_m(V) & 0 & 0\\
0 & \alpha_h(V) & 0 & 0 & 3\alpha_m(V) & -D_{66}(V) & 2\beta_m(V) & 0\\
0 & 0 & \alpha_h(V) & 0 & 0 & 2\alpha_m(V) & -D_{77}(V) & 3\beta_m(V)\\
0 & 0 & 0 & \alpha_h(V) & 0 & 0 & \alpha_m(V) & -D_{88}(V)}
\end{eqnarray*}
where $D_{ii}(V) = \sum_{j=1}^8 A_{ij}(V)$ from the adjacency matrix above (Equation \ref{eq:HH-Na-adjacency}).  The noise matrix $B$ is also voltage-dependent and is given by the general expression in Equation \ref{eq:generalB}.

\begin{figure}[ht!] 
\centering
	\includegraphics[width=4in]{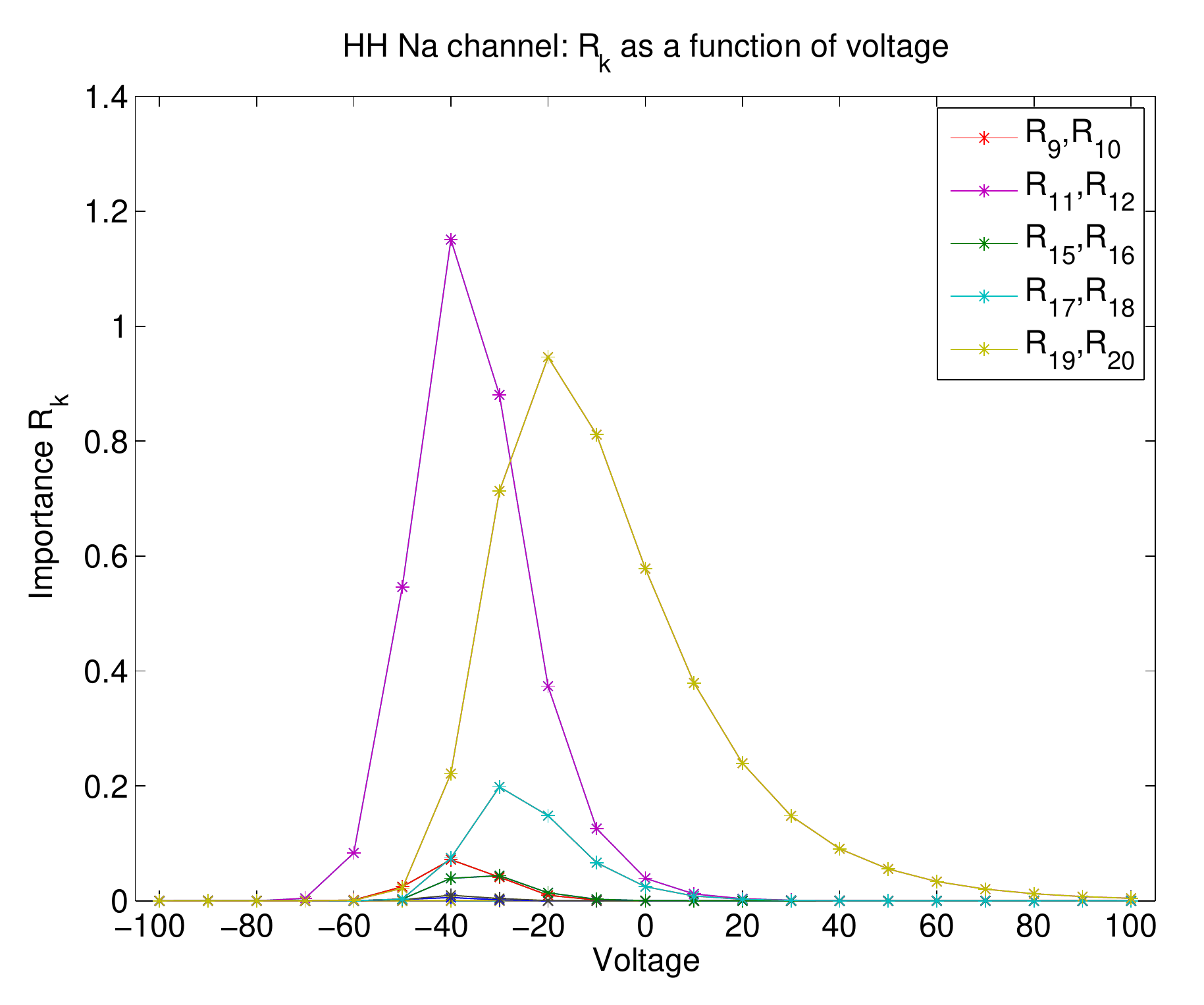}
\caption{Hodgkin-Huxley sodium channel.  This figure shows edge importance $R_k$ as a function of voltage in the range $[-100,100]$ mV for each reaction $k\in\{1,\dots,20\}$.  The magenta curve corresponds to edges 11 and 12 and the yellow curve corresponds to edges 19 and 20 (transitions between the conducting state 8 and its two nearest neighbors, states 7 and 4, respectively).  Note that $R_{11}=R_{12}$ (magenta) has the highest edge importance in the voltage range $[-100,-25]$ mV and $R_{19}=R_{20}$ (yellow) has the highest value in the range $[-25,100]$ mV.}
\label{fig:HH-Na-channel-Rk-voltage}
\end{figure}

Figure \ref{fig:HH-Na-channel-Rk-voltage} shows the edge importance $R_k$ as a function of voltage for each reaction $k\in\{1,\dots,20\}$ for the sodium channel state graph.  The sodium channel also satisfies detailed balance, so each pair of complementary reactions $k_i, k_{i+1}$ connecting the same pair of nodes will have equal edge importance values $R_{k_i}=R_{k_{i+1}}$.  The magenta curve corresponds to edges 11 and 12 and the yellow curve corresponds to edges 19 and 20, which are the transitions between state 7 and conducting state 8, and the transitions between state 4 and conducting state 8, respectively.  Note that $R_{11}=R_{12} > R_k$ (magenta) for all other reactions $k$ in the voltage range $[-100,-25]$ mV and then it switches so that $R_{19}=R_{20} > R_k$ (yellow) for all other reactions $k$ in the range $[-25,100]$ mV.  This means that if any of these four reactions are neglected, they would have the highest contribution to the error.

%%%%%%%%%%%%%%%%%%%%%%%%%%%%%%%
% Figure 13
%%%%%%%%%%%%%%%%%%%%%%%%%%%%%%% 

\begin{figure}[ht!] 
\centering
	\subfigure{
	\includegraphics[width=3.07in]{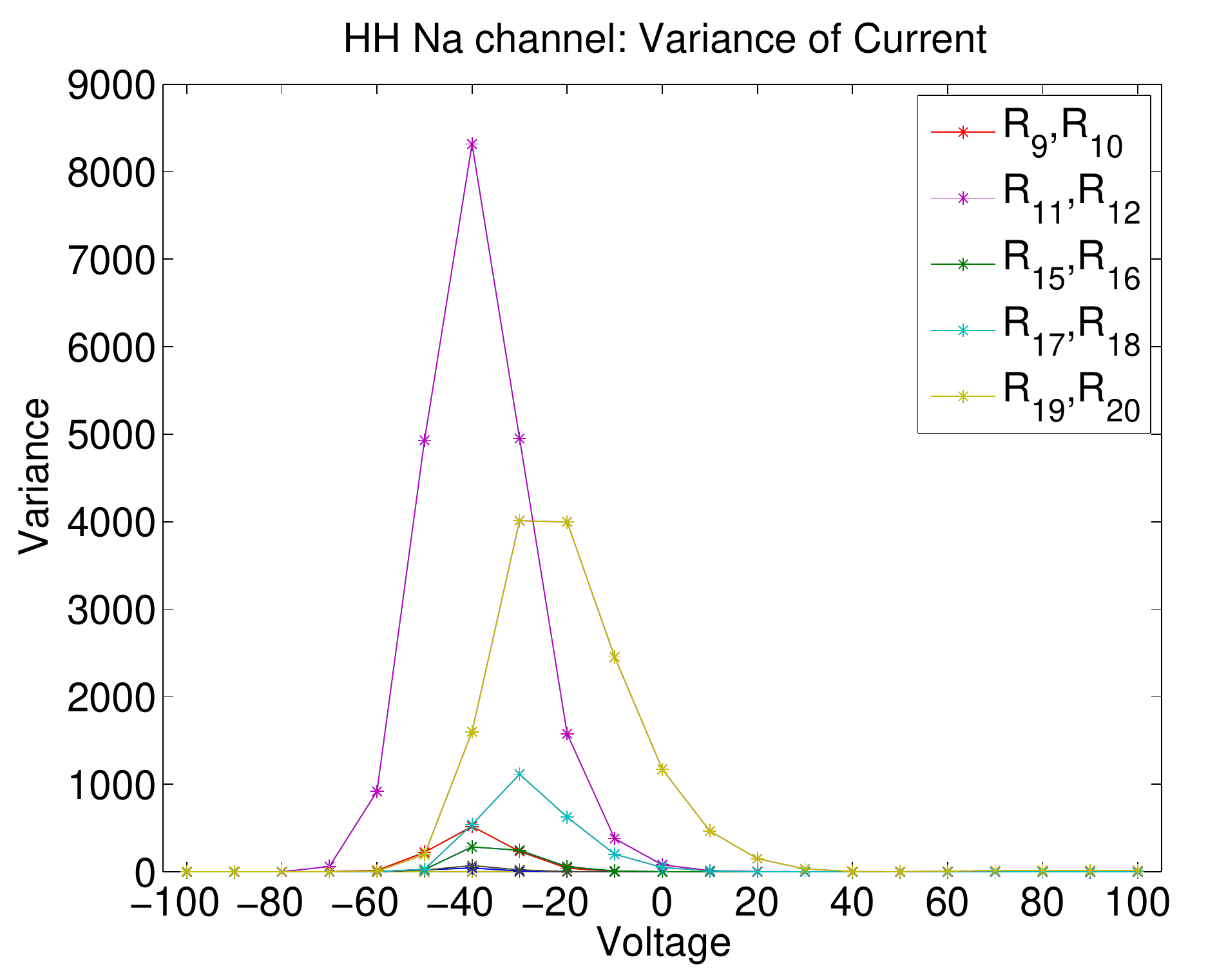}
	}
	\subfigure{
	\includegraphics[width=3in]{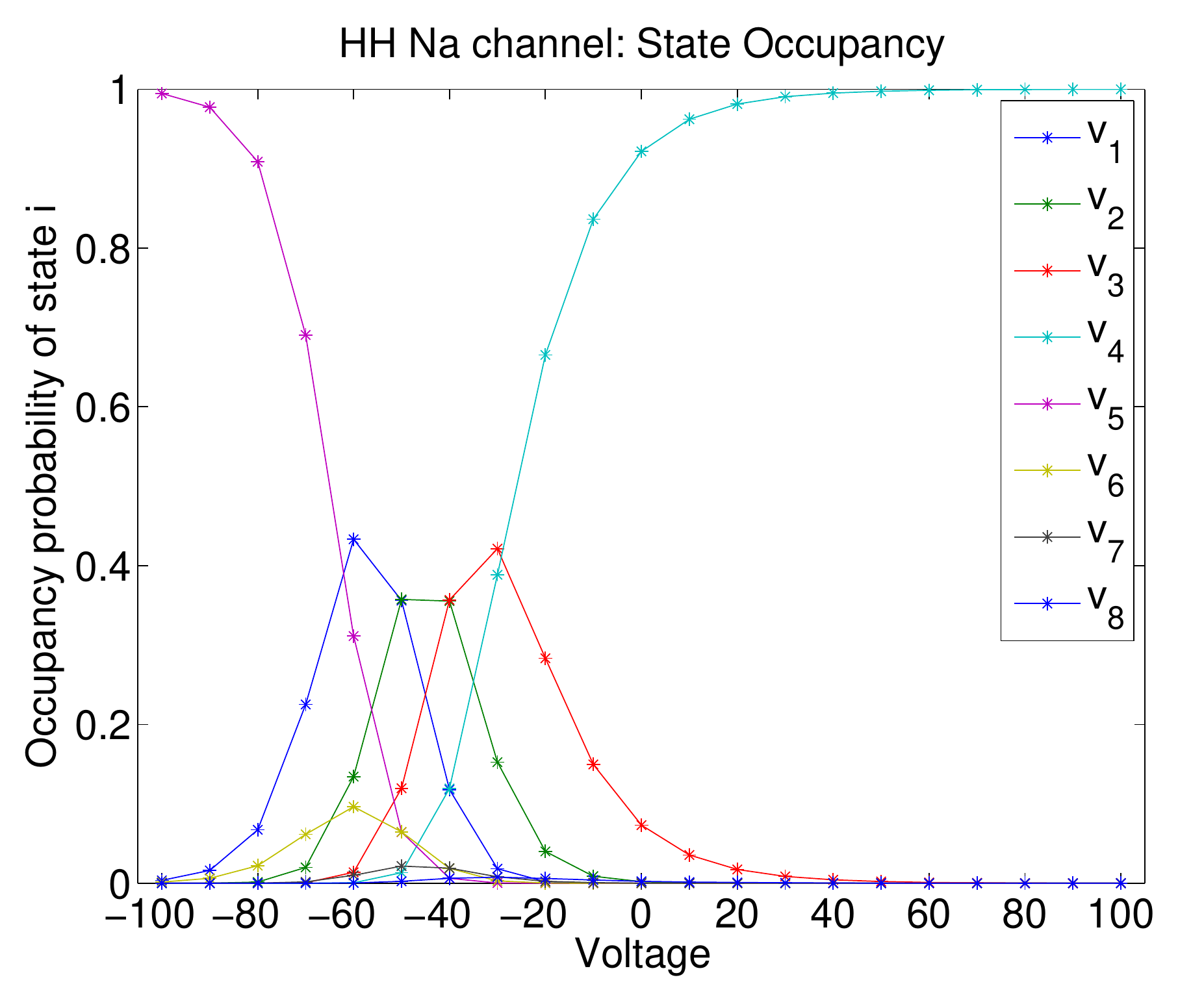}
	}
\caption{Hodgkin-Huxley sodium channel. Left: variance of the current $R_k * (V - V_k)^2$ as a function of voltage $V$ for each reaction $k$ where $V_k=45$ mV is the reversal potential for the sodium channel.  The magenta curve ($R_{11}=R_{12}$, corresponding to the transitions between state 7 and conducting state 8) has the largest variance in the voltage range $[-100,-25]$ mV and the yellow curve ($R_{19}=R_{20}$, corresponding to the transitions between state 4 and conducting state 8) has the largest variance in the voltage range $[-25,100]$ mV.  Right: leading eigenvector components (normalized so that the components sum to 1) as a function of voltage. This shows the probability of being in each state and how that changes with voltage.}
\label{fig:HH-Na-channel-variance-current-state-occupancy}
\end{figure}

Figure \ref{fig:HH-Na-channel-variance-current-state-occupancy} shows the variance of the nominal current $R_k * (V - V_k)^2$ as a function of voltage $V$ for each reaction $k$ where $V_k=45$ mV is the reversal potential for the sodium channel. Again, we choose units for conductance such that the unitary channel conductance equals 1. As before, we see that the edges with the highest edge importance have the largest variance (left panel).  The switch between the dominant curves (magenta vs yellow) agrees with the switch in Figure \ref{fig:HH-Na-channel-Rk-voltage} which occurs at -25 mV.  The right panel in Figure \ref{fig:HH-Na-channel-variance-current-state-occupancy} shows the probability of being in each state and how that changes with voltage.

%%%%%%%%%%%%%%%%%%%%%%%%%%%%%%%
% Figure 14
%%%%%%%%%%%%%%%%%%%%%%%%%%%%%%% 

In summary, our analysis fully supports the accuracy of Schmandt and Gal\'{a}n's stochastic shielding algorithm for the Hodgkin-Huxley system, at least for the voltage clamped case that we consider.  More significantly, our analysis allows one to calculate the relative importance of \emph{each} transition in a network of first order reactions, allowing a new quantitative basis for reduction of complexity of stochastic network models.  In the case of a simple chain of states such as the Hodgkin-Huxley potassium channel, the rank ordering of transitions by importance $R_k$ is the same for all voltages. As shown in Figure \ref{fig:HH-Na-channel-Rk-voltage}, however, for more complicated gating schemes, such as the Hodgkin-Huxley sodium channel, the rank ordering of transitions by importance can differ at different voltages.  

For instance, the most important transition at subthreshold voltages ($V\lesssim -40$ mV) is the transition connecting the $[m=(1,1,0),h=1]$ state (state 7 in Figure \ref{fig:HH-Na-channel}) to the $[m=(1,1,1),h=1]$ state (state 8, the conducting state). This transition corresponds biophysically to the nonconducting-to-conducting transition that occurs via \emph{activation} or \emph{deactivation} \cite{jp:Hodgkin+Huxley:1952d}, that is, the opening (or closing) of the last of three $m$-activation gates in the ion channel.  
It is significant that this transition is the most ``important" for subthreshold voltages, because the activation transition is typically the last subthreshold event during spike generation.  

On the other hand, at suprathreshold voltages the most important transition is that connecting the $[m=(1,1,1),h=1]$ state (state 8) with the $[m=(1,1,1),h=0]$ state (state 4).  Biophysically, this transition corresponds to \emph{inactivation} and \emph{deinactivation}, or the closing (and opening) of the $h$-inactivation gate.  During action potential generation this transition plays an essential role in terminating the voltage spike upstroke, and it is significant that it should be most ``important" at suprathreshold voltages.

For more general channel schemes, and more elaborate stochastic processes in general, the identification of the relative quantitative importance of different transitions or edges to the observable behavior of the system is a powerful new tool for principled complexity reduction.

%\newpage
\section{Discussion}\label{sec:discuss}

In the ongoing race between growth of empirical data sets and growth of available computing power, conceptual understanding of complex dynamical systems can get left behind.  Finding efficient lower-dimensional representations of high-dimensional systems, that accurately capture relevant aspects of system behavior, not only takes better advantage of computational resources, but can provide insights into the essential components of a system.  Hence, there has been a significant effort in recent years to develop principled complexity reduction techniques for naturally occurring complex networks.  

Schmandt and \galan \cite{SchmandtGalan2012PRL} developed a method for efficient simulation of stochastic ion channel gating in the membrane of a neuron.  The random gating of ion channels provides an important class of biological processes which are naturally represented as Markov chains on graphs \cite{ClayDeFelice1983BiophysJ,SkaugenWalloe1979ActaPhysiolScand}.  The graphs in this case arise from the different configurations of ion channel subunits or ``gates". Typically each state carries one of two functional labels: open or closed.  This coarse-grained representation of the ion channel corresponds to a linear measurement functional, in the sense that current flowing through open channels can be measured experimentally, and individual ion channels typically exhibit binary all-or-none conductance.  Schmandt and \galan implemented a novel form of coarse graining technique that ignores fluctuations between indistinguishable transitions (open-to-open or closed-to-closed) while preserving fluctuations between distinguishable states.  In order to gain a deeper understanding of why their ``stochastic shielding approximation'' works so well, we analyzed it in the context of a multidimensional Ornstein-Uhlenbeck process on a variety of networks.  First, we showed that this form of model reduction can be represented as a mapping from a many-dimensional sample space to a lower-dimensional sample space, rather than as a mapping from a many-node network to a few-node network, and that one can formulate the problem as a search for the optimal such mapping.  Second, we showed that for the specific 3-state example presented in Schmandt and Gal\'{a}n's paper, their approximation is indeed optimal in a specific sense.  Third, we obtained a theoretical result showing that stochastic shielding works for an ensemble of random graphs with arbitrarily chosen binary measurement vectors, analogous to the identification of nodes as conducting versus non-conducting in ion channel models.  Finally, we evaluated the stochastic shielding approach for the graph representing the ion channel states of the classical Hodgkin-Huxley model, and showed that this approach is optimal for a wide range of fixed voltages under ``voltage clamped'' conditions.

{\bf Relationship between different levels of modeling}\\ 
The underlying description of Schmandt and Gal\'{a}n's model \cite{SchmandtGalan2012PRL} is given by the population process described in \S \ref{ssec:PopulationProcess}, a more general framework than the Ornstein-Uhlenbeck process that we study.  The OU process connects to the population process via a tau-leaping approximation, as described in Appendix \ref{app:tau-leaping}.  The tau-leaping method involves two key assumptions.  First, assuming that the transition propensities $\alpha_{ij(k)}$ do not change dramatically in an interval of length $\tau$, we can approximate the number of transitions in each interval by a collection of independent Poisson processes.  This approach is closely related to the framework of Schmandt and Gal\'{a}n, except that they use a binomial distribution instead of a multinomial distribution (see Appendix \ref{app:stoch-shield-GS}).  Second, if the expected number of occurrences of each reaction is sufficiently large (i.e.~10s or 100s) in time  $\tau$, then it is reasonable to use a Gaussian approximation to the Poisson process.  The resulting model comprises the standard chemical Langevin formulation, in which  the size of the fluctuations associated with each transition is state dependent.  These two constraints are can always be satisfied by taking a sufficiently large number of individuals in the population.
The Ornstein-Uhlenbeck process is obtained by linearizing about the mean field steady state distribution of the tau-leaping model (see Appendix \ref{app:tau-leaping}).  The intensity of the noise terms is determined by the mean steady state occupancy of each state, resulting in a linear OU process.  Although our analysis is limited to the OU process version of the system, it is reasonable to expect that stochastic shielding will apply more broadly.  For example, in the full population process one can decompose the fluxes in the model into a sum of a mean component and a mean zero fluctuating component.  In this case, stochastic shielding amounts to setting the fluctuating component to zero while preserving the mean for those transitions connecting observationally equivalent states.

%[ add comment about nonlinear systems, spiking, limit cycles, and other transitions might turn out to be important also given some global measure of behavior]
Limiting the investigation to voltage clamped conditions facilitated a more thorough mathematical analysis of the stochastic shielding approximation, but also restricted the biological applicability of the results.  By approximating the population process with a closely related Ornstein-Uhlenbeck process  we effectively linearized the system about a fixed point given by the mean field behavior.  Therefore our analysis does not address important nonlinear dynamical behaviors  arising in many physical and biological systems, such as noise driven transport between multiple quasiequilibria, fluctuation induced spiking in excitable systems (including noise induced spiking in nerve cells), or limit cycle oscillations (including regular spiking in nerve cells).  On the one hand, we anticipate that transitions in a state graph corresponding to directly observable state changes, such as between conducting and nonconducting ion channel states, will remain ``important" under more general measures accounting for global, nonlinear behaviors.  On the other hand, it is certainly possible that additional transitions may also become important with respect to more general measures, if the linear measurement vectors employed here fail to capture their contribution to global dynamics.

{\bf Broader applications}\\ 
The stochastic shielding approximation can be directly applied to various biological networks, not just ion channel models.  For instance, Lu et al \cite{LuShenZongHastyWolynes2006PNAS} describe a signal transduction network in which the phosphorylation and transport events are arranged with a ladder topology.  The two sides of the ladder denote molecules in the nucleus and in the cytoplasm, respectively.  On each side, there are $M+1$ species having different levels of phosphorylation (see Figure 1 of \cite{LuShenZongHastyWolynes2006PNAS} for an illustration).  This is a more elaborate Markov process than a simple ion channel state model, but it can still be described with a binary measurement vector.  The readout is 1 if the system is both in the nucleus and in a specific phosphorylated state, and 0 otherwise.  The application of stochastic shielding to such a system is quite natural.  

Another broad class of examples includes calcium-induced calcium release Markov models.  Nguyen, Mathias and Smith \cite{NguyenMathiasSmith2005BullMathBio} studied a stochastic automata network description of instantaneously coupled intracellular calcium channels which they derived from Markov models of single channel gating that include calcium activation, inactivation, or both.  This high dimensional system involves a large number of functional transitions; the transition probabilities of one channel depend on the local calcium concentration which is typically influenced in turn by the state of other channels in the population.  Such models can easily become very high dimensional.  For example, DeRemigio et al \cite{DeRemigioEtAlSmith2008} considered a discrete state continuous time Markov model of coupled calcium channels, taking explicit channel position in to account, which yields up to 1.6 million distinct states. For systems of such complexity, any reduction of the complexity of the stochastic process by stochastic shielding will likely be advantageous, both for simulation and for analysis.
%Other similar systems might also involve multiplicative noise processes which is an important large case of systems that go beyond our framework thus far.

We have focused here on discrete state ion channel models with binary measurement vectors.  However, it is possible that some ion channels may have a richer than binary readout structure.  For example, Catterall and colleagues \cite{Catterall2012JPhysiol} provide structural evidence that activation of a bacterial sodium channel may possess multiple non-equivalent conducting states, raising the possibility that conductance could be graded rather than binary.  As another example in which could lead to graded measurement vectors, adaptive evolution can be represented on a random walk on a graph representing genomic variants connected by possible mutation routes \cite{KauffmanLevin1987JTB,McCandlish2011Evolution}.  
While the stochastic process representing the evolution of a human pathogen such as influenza may have an enormous number of degrees of freedom \cite{KryazhimskiyBazykinPlotkinDushoff2007ProcBiolSci,NelsonHolmes2007NatRevGenet}, the dynamics of interest may comprise a smaller number of dimensions, such as a strain's virulence or fitness, which may naturally be graded rather than discrete quantities.

\begin{figure}[htbp] 
   \centering
  \includegraphics[width=3in]{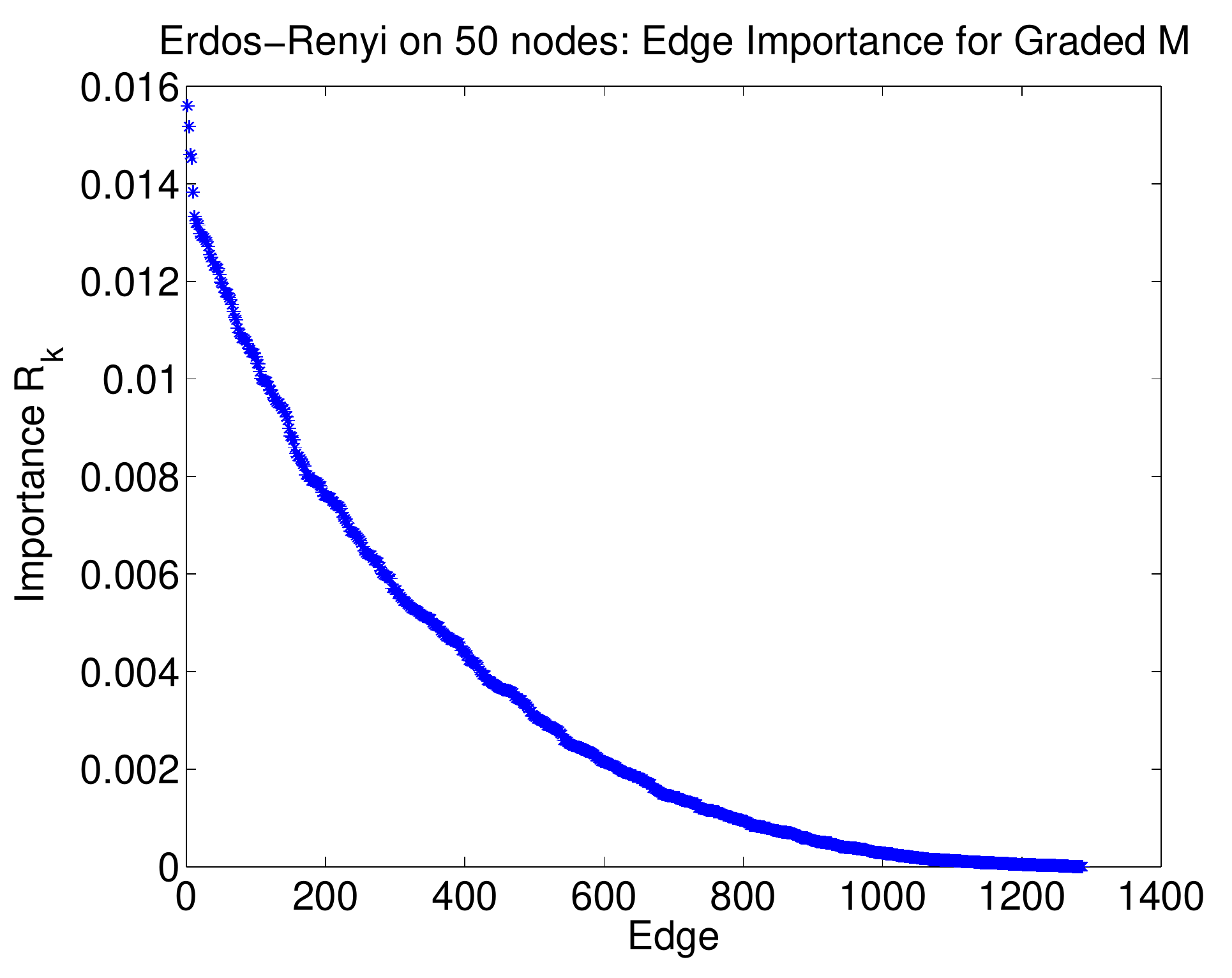}    \includegraphics[width=3in]{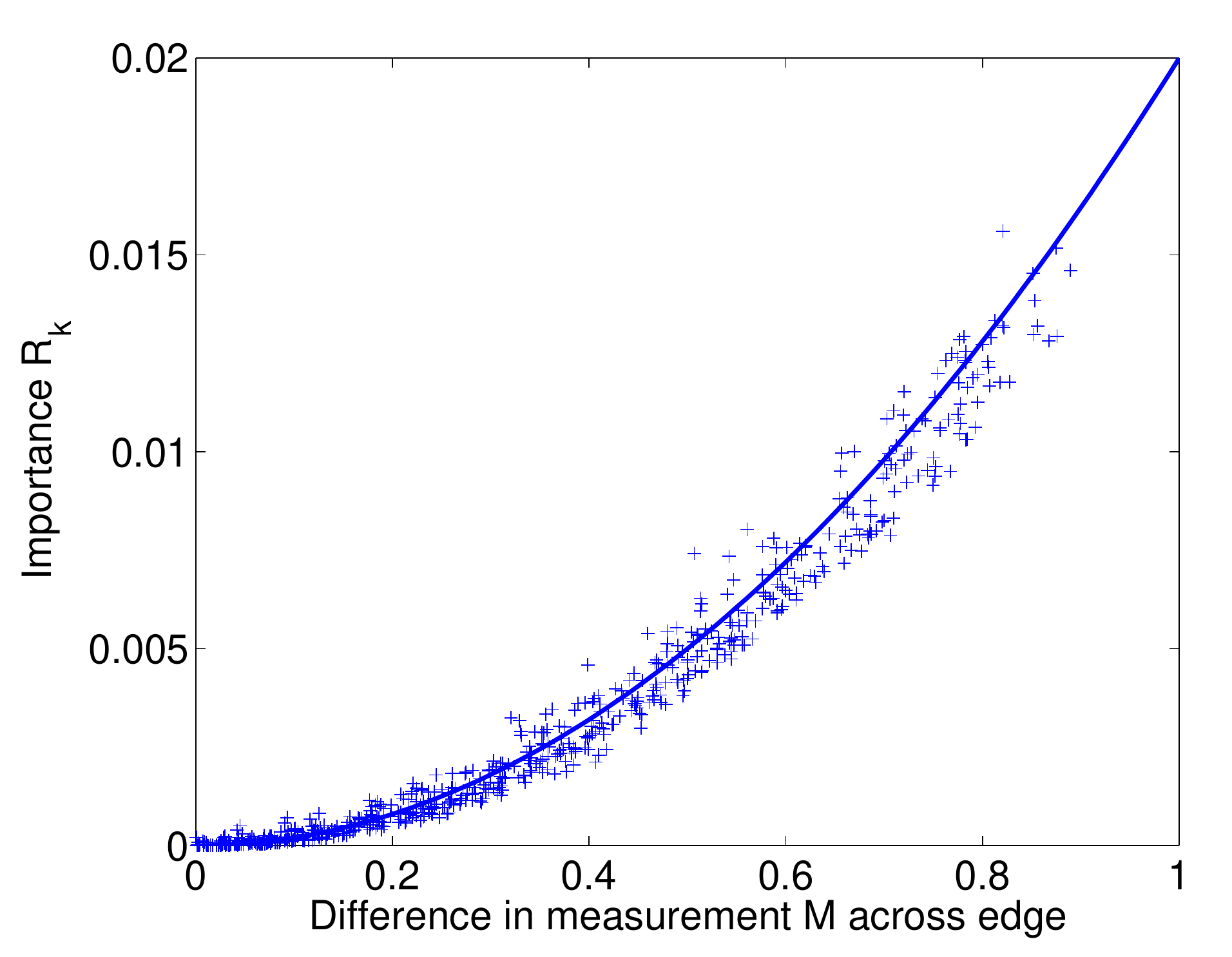} 
   \caption{Edge importance distribution for graded measurement vector $M$. The effect of neglecting the fluctuations associated with the $k^{\text{th}}$ edge in an \ER network with $n=50$ nodes and edge probability $p=0.5$, as a function of the difference in measurement $M$ at the two ends of the edge, $M^\intercal\zeta_k$.  In this example, the components of $M$ were assigned from the uniform distribution on $[0,1]$, independently of the presence or absence of edges in the graph.  Left: Rank order plot of edge importance $R_k$.  Compare to Figure \ref{fig:ER_edge_importance50}; note the absence of a clear gap distinguishing ``important'' from ``unimportant'' edges.  Right: Horizontal axis, $x=|M^\intercal\zeta_k|$.  Vertical axis, $R_k$. The superimposed curve shows the quadratic $y \approx x^2/n$, for $n=50$.}
   \label{fig:graded}
\end{figure}

Stochastic shielding in a modified form would still apply even if the measurement functional were graded continuously.  As an example, consider an \ER random graph on $n$ nodes with edge probability $p$, with graded measurement vector $M\in[0,1]^n$ instead of binary $M\in\{0,1\}^n$.  The left panel of Figure \ref{fig:graded} shows the edge importance distribution for the case $n=50$ and $p=0.5$ where the components of $M$ are chosen uniformly at random from the unit interval.  The right panel of Figure \ref{fig:graded} illustrates the difference in measurement between nodes connected by edge $k$, $x=|M^\intercal\zeta_k|$, versus the edge importance $R_k$, and shows good agreement with the curve $y \approx x^2/n$ for the case $n=50$. 

This empirical result (Figure \ref{fig:graded}, right panel) suggests the following generalization of Theorem \ref{thm:MainResult}:
\begin{equation}\label{eq:graded}
E[R_k|M]=\frac{\sigma_k^2 (M^\intercal\zeta_k)^2}{n\,C}+O(n^{-q}),\mbox{ as }n\to\infty,
\end{equation}
for some $q>1$ (e.g., $q=2$ for the Gaussian unitary ensemble, and $q\approx 5/3$ for the \ER ensemble, empirically).  
In the case of a binary measurement vector, $M \in\{0,1\}$, this formula would revert to the result given in Theorem \ref{thm:MainResult}.  
%%%%%%%%%%%%%%%%%%%%%%%%%%%%%%%
% Figure 15
%%%%%%%%%%%%%%%%%%%%%%%%%%%%%%% 
A rigorous derivation of Equation \ref{eq:graded} is beyond the scope of the present paper.

The behavior of stochastic processes arising in first-order reaction networks has been explored in broad generality by Othmer and colleagues \cite{CadgilLeeOthmer2005BullMathBiol}.  They used a spectral approach to analyze a general system of first-order reaction networks, and studied the effect of changes in the network topology on the distribution of the number of reactant molecules, as well as the difference between conversion and catalytic networks with the same topology.  Exploring sample space reductions conditioned on a linear measurement functional for such general classes of networks would be of interest.

{\bf Different levels of model simplification}\\
Model simplification is an important goal for Markov chain models in many scientific contexts, and complexity reduction has been pursued through a corresponding variety of approaches.  Newman and others have extensively developed techniques based on community structure, aggregating or lumping nodes together based on topological considerations \cite{Newman2006PNAS,NewmanBook2010}.  When applied to a stochastic process on a graph, the aggregation of $N\gg n$ to $n$ nodes is equivalent to a projection of the original process onto a subspace in which the process components on the aggregated fine-grained nodes are averaged.  
In most cases, the resulting coarsened process is no longer Markov, although in some cases exact dimension reduction to a lower dimensional Markov processes can be accomplished \cite{AgarwalaChielThomas2012JTB,Buchholz1994JApplPr,WeiKuo1969IECFundamentals}.
Other aggregation schemes, such as spectral coarse graining \cite{GfellerThesis2007,GfellerDeLosRios2007PRL,GfellerDeLosRios2008PRL}, have been proposed based on the spectral properties of the graph Laplacian.  Approaches based on topological or abstract spectral properties do not necessarily take into account functional properties of the system to be simplified.  Because stochastic shielding simplifies the representation of a stochastic process \emph{taking into account the function of the system}, namely by distinguishing conducting versus nonconducting ion channel states, it may provide insights not afforded by graph aggregation based on modularity or graph spectra.

As another example of simplification based on functional properties, Bruno, Yang and Pearson \cite{BrunoYangPearson2005PNAS} used independent open-closed transitions to describe a canonical form that can express all possible reaction schemes for binary ion channels.  

Not all prior approaches to simplification of random processes on graphs proceed by aggregating nodes.  For instance, Pearson and colleagues \cite{UllahBrunoPearson2012JTB} proposed model simplification by the elimination of nodes with low equilibrium occupancy probability using time scale separation arguments. The reduced system has fewer parameters, and the dynamics of the reduced system are identical to those of the original system except on very fast time scales.  Other simplifications based on graph sparsification have been proposed by Koutis and colleagues \cite{KoutisLevinPeng2012}.

In this paper we have investigated a novel form of simplification of stochastic processes on graphs.  Stochastic shielding is based on replacing a high-dimensional stochastic process defined on a graph with a lower-dimensional process on the \emph{same} graph, rather than replacing a complex network with a simpler one.  Specifically, we consider mappings from the original process to an approximate process defined on a significantly smaller \emph{sample space}.
In one sense, we can think of the full and a reduced system as two systems with partially shared stochastic input, and partially independent stochastic input of different magnitudes (magnitude zero, in one case).  Structurally, this situation is analogous to the kind of mixed common-noise and independent-noise scenario studied in the context of neuronal synchronization \cite{RochaDoironShea-BrownJosic2007Nature,JosicShea-BrownDoironDeLaRocha2009NeuralComp,Shea-BrownJosicDeLaRochaDoiron2008PRL}.  In another sense, stochastic shielding can be seen as a different kind of projection, \textit{versus} that induced by lumping or pruning nodes.  The latter methods simplify the graph, whereas stochastic shielding leaves the graph unchanged and simplifies the sample space on which the approximate process lives.

\section*{Acknowledgments}
This work was supported by the National Science Foundation (grant EF-1038677) and also in part by the Mathematical Biosciences Institute and the NSF (grant DMS 0931642), by a grant from the Simons Foundation (\#259837 to PJT), and by the Council for the International Exchange of Scholars (CIES).  The authors thank R.~\galan and N.~Schmandt for bringing the stochastic shielding approximation to their attention, and for helpful discussions.  Thanks also to H.~Chiel for critical comments on the manuscript.

\bibliographystyle{plain}
%\bibliography{math,Snyder,networks,stoch_chem,Dicty,neuroscience,PJT}
\bibliography{Schmidt-Thomas-JMN}

\appendix

\section{Stochastic Shielding Construction of Schmandt and \galan}
\label{app:stoch-shield-GS}

In \cite{SchmandtGalan2012PRL}, Schmandt and \galan considered discrete time simulations approximating a continuous time, finite state Markov chain
\begin{eqnarray}
N_i(t)&=&N_i(0)+\sum_{j\not=i}\left(\tilde{N}_{ji}(t)-\tilde{N}_{ij}(t)\right)
\end{eqnarray}
where $N_i(t)$ is the number of individuals in a population (of size $N_{\text{tot}}$) in state $i$ at time $t$, and $\tilde{N}_{ij}(t)$ counts the number of $i\to j$ transitions that have occurred as of time $t$. The transition counts $\tilde{N}_{ij}(t)$ may be written using the random time change representation \cite{AndersonKurtz2011chapter} as
\begin{eqnarray}
\tilde{N}_{ij}(t)&=&Y_{ij}\left[\int_{s=0}^tN_i(s)\alpha_{ij}(s)\,ds\right].
\end{eqnarray}
By convention we take $N_{ii}(t)\equiv 0$ and $\alpha_{ii}(t)\equiv 0$.  The $Y_{ij}$ are independent unit rate Poisson processes driving the different state-to-state transitions.  The transition from state $i$ to state $j$ occurs with \textit{per capita} rate $\alpha_{ij}$.    In a conductance based model, such as a discrete stochastic version of the Hodgkin Huxley equations, the vector $(N_1(t),\cdots,N_K(t))$ would represent the number of ion channels in each of $K$ distinct states, and the transition rates could vary with time, e.g.~through dependence on membrane potential or second messenger concentration.   Although Schmandt and \galan consider both the stationary and time varying case, we restrict attention to the stationary case, which corresponds experimentally to a voltage clamped preparation.

One may (approximately) simulate trajectories of the Markov chain using a discrete time step approach.  Following \cite{SchmandtGalan2012PRL}, we fix a time step $h>0$ and define $N_{ij}$ as $N_{ij}(t)=\tilde{N}_{ij}(t+h) -\tilde{N}_{ij}(t)$, that is,  the number of $i\to j$ transitions occurring in the interval $(t,t+h]$.  The net increments in the state occupancy numbers $N_i$ are then given by 
\begin{equation}
\Delta_i(t)\equiv N_i(t+h)-N_i(t)=\sum_{j\not= i} N_{ji}(t)-N_{ij}(t).
\end{equation}
To obtain a practical algorithm, Schmandt and \galan set $N_{ij}(t)\sim\text{Binom}[N_i(t),\alpha_{ij}(t)h]$.  Since there is then a finite probability that $N_i(t+h)<0$, one must include an iterative resampling scheme to force $N_i(t+h)\ge 0$.  As an alternative, we consider instead a multinomial representation of the destinations of all $N_i(t)$ individuals beginning the time step at node $i$.  That is,  for each $i$, $1\le i \le K$, we set
\begin{equation}
(N_{i1},\cdots,N_{ii},\cdots, N_{iK})\sim\text{Multi}\left[N_i(t),\left(\alpha_{i1}h,\cdots,\left(1-\sum_{j\not=i}\alpha_{ij}h\right),\cdots,\alpha_{iK}h\right)\right].
\end{equation}
The multinomial distribution produces an integer-valued random vector with mean and marginal distributions the same as that given by the binomial distribution; the only difference is that transitions emanating from a common node are not assumed to be independent.  

The first and second moments arising from the multinomial transition distribution are
\begin{eqnarray}
E[N_{ij}|\vec{N}(t)]&=&N_i(t)\alpha_{ij}h,\text{ for }i\not=j\\
E[N_{ii}|\vec{N}(t)]&=&N_i(t)\left(1-\sum_{j\not=i}\alpha_{ij}h\right)=N_i(t)-\sum_{j\not=i}E[N_{ij}]\\
V[N_{ij}|\vec{N}(t)]&=&N_i(t)\alpha_{ij}h(1-\alpha_{ij}h),\text{ for }i\not=j\\
V[N_{ii}|\vec{N}(t)]&=&N_i(t)\left(\sum_{j\not=i}\alpha_{ij}h \right)\left(1-\sum_{j\not=i}\alpha_{ij}h\right)\\
\cov[N_{ij},N_{ij'}|\vec{N}(t)]&=&-N_i(t)\alpha_{ij}\alpha_{ij'}h^2,\text{ for }j\not=j', j\not=i, j'\not=i\\
\cov[N_{ij},N_{ii}|\vec{N}(t)]&=&-N_i(t)\alpha_{ij}h\left(1-\sum_{j'\not=i}\alpha_{ij'}h\right),\text{ for } j\not=i.
\end{eqnarray}
Here all expectations are conditioned on the current state of the system, 
$$\vec{N}(t)=(N_1(t),\cdots,N_i(t),\cdots,N_K(t)).$$

The mean increment given the current distribution of the population, $\bar{\Delta}_i(t)\equiv E[\Delta_i(t)|\vec{N}(t)]$, is written in terms of the mean transitions as
\begin{equation}
\bar{\Delta}_i(t)
=\sum_{j\not=i}\left(E[N_{ji}(t)|N_j(t)] - E[N_{ij}(t)|N_i(t)]\right)
=\sum_{j\not=i}\left(N_j(t)\alpha_{ji}h-N_i(t)\alpha_{ij}h\right).
\end{equation}
The deviation of the actual number of $i\to j$ transitions from the expected number is 
\begin{equation}
\delta\Delta_i(t)\equiv \Delta_i(t)-\bar{\Delta}_i(t)
=\sum_{j\not=i}\left((N_{ji}(t)-N_j(t)\alpha_{ji}h)-(N_{ij}(t)-N_i(t)\alpha_{ij}h)\right)
=\sum_{j\not=i}(\delta N_{ji}(t)-\delta N_{ij}(t))
\end{equation}
where $\delta N_{ij}(t)=N_{ij}(t)-E[N_{ij}(t)|\vec{N}(t)]$ is the deviation of the number of $i\to j$ transitions from the expected number.  
The mean of $\delta N_{ij}(t)$ is zero for all $i,j,$ and all $t$, by construction.  
The stochastic shielding approximation amounts to setting $\delta N_{ij}(t)$ to zero for selected $i\to j$ transitions, namely for those transitions between ``unobservable states", or (equivalently) between any two states with the same value of the measurement observable, i.e.~the conductance.   Since $E[\delta N_{ij}(t)|\vec{N}(t)]=0$ already, the only error introduced by suppressing the fluctuations associated with the $i\to j$ transition comes from the propagation of the fluctuations through the network to the observable states.  But the fluctuations in the transitions, $N_{ij}$, are only weakly correlated with the fluctuations in the occupancy numbers of observable states, $N_k(t)$, when $i$ and $j$ have the same conductance.  To introduce this shielding effect, Schmandt and \galan calculate the second moments for the population increments $\delta \Delta_i(t)$.  As an example, in the three node case, for the multinomial transition model, the variances are given by
\begin{eqnarray}
E[\delta\Delta^2_1(t)|\vec{N}(t)]&=& V[N_{12}|\vec{N}(t)]+V[N_{21}|\vec{N}(t)]\\ \nonumber 
&=&N_1(t)\alpha_{12}h(1-\alpha_{12}h)+N_2(t)\alpha_{21}h(1-\alpha_{21}h)\\  
E[\delta\Delta^2_2(t)|\vec{N}(t)]
&=& V[N_{12}|\vec{N}(t)]+V[N_{21}|\vec{N}(t)]+ V[N_{23}|\vec{N}(t)]+V[N_{32}|\vec{N}(t)]+\\  \nonumber 
&&2 \cov[N_{21},N_{23}|\vec{N}(t)]\\  
&=&N_1(t)\alpha_{12}h(1-\alpha_{12}h)+N_2(t)\alpha_{21}h(1-\alpha_{21}h)+\\  \nonumber 
&&N_2(t)\alpha_{23}h(1-\alpha_{23}h)+N_3(t)\alpha_{32}h(1-\alpha_{32}h)\\  \nonumber 
&&-2N_2(t)\alpha_{21}\alpha_{23}h^2
\\
E[\delta\Delta^2_3(t)|\vec{N}(t)]&=&  V[N_{23}|\vec{N}(t)]+V[N_{32}|\vec{N}(t)]\\ \nonumber 
&=&N_2(t)\alpha_{23}h(1-\alpha_{23}h)+N_3(t)\alpha_{32}h(1-\alpha_{32}h)
\end{eqnarray}
and the covariances are given by
\begin{eqnarray}  
E[\delta\Delta_1(t)\delta\Delta_2(t)|\vec{N}(t)]&=& -V[N_{12}|\vec{N}(t)]-V[N_{21}|\vec{N}(t)] - \cov[N_{21},N_{23}|\vec{N}(t)]\\  \nonumber 
&=&-N_1\alpha_{12}h(1-\alpha_{12}h)-N_2\alpha_{21}h(1-\alpha_{21}h)+N_2\alpha_{21}\alpha_{23}h^2\\
E[\delta\Delta_1(t)\delta\Delta_3(t)|\vec{N}(t)]&=& \cov[N_{21},N_{23}|\vec{N}(t)]
=-N_2(t)\alpha_{21}\alpha_{23}h^2\\
E[\delta\Delta_2(t)\delta\Delta_3(t)|\vec{N}(t)]&=&  \nonumber 
 -V[N_{23}|\vec{N}(t)]-V[N_{32}|\vec{N}(t)] - \cov[N_{21},N_{23}|\vec{N}(t)]\\
&=&-N_2\alpha_{23}h(1-\alpha_{23}h)-N_3\alpha_{32}h(1-\alpha_{32}h)+N_2\alpha_{21}\alpha_{23}h^2
\end{eqnarray} 
Schmandt and \galan obtain similar expressions that agree up to order $O(h)$; the difference between the binomial and multinomial expressions only appears in the $O(h^2)$ terms.  For example, they assert that $E[\delta\Delta_1(t)\delta\Delta_3(t)|\vec{N}(t)]\equiv0$, while under the multinomial model this covariance is equal to $-N_2(t)\alpha_{21}\alpha_{23}h^2$. Fortunately, this difference does not undermine the main argument.

From this point, Schmandt and \galan obtain an expression for the stationary covariance matrix of the reduced process (compare equation (8) in \cite{SchmandtGalan2012PRL}   with our Lemma 1) and decompose the covariance into a sum over direct and indirect connections to 
a single conducting or observable state.   This situation corresponds, in our analysis, to the case where the measurement vector $M$ contains a single nonzero entry.   Schmandt and \galan
argue that suppressing the fluctuations associated with transitions not directly affecting the observable state decrease their contribution to the variance of the observable state occupancy, while increasing the contribution of the direct transitions to the same variance.  In addition, they show through numerical comparisons that Hodgkin-Huxley equations with a full Markov process and the reduced process are practically indistinguishable both under voltage clamp (stationary transition rates) and current clamp (time varying transition rates) conditions.

\section{Derivation of Tau-Leaping for an Arbitrary Finite Graph}
\label{app:tau-leaping}
\subsection{Tau-Leaping: General Case}

We will use standard tau-leaping arguments \cite{Gillespie2001JChemPhys,PetzoldGillespie2003JChemPhys} to derive the multidimensional Ornstein-Uhlenbeck process in \S \ref{ssec:OUP} (Equation \ref{eq:OUPforX}).  Given a symmetric directed graph $\mathcal{G}(\mathcal{V},\mathcal{E})$ with $n$ nodes, let $N(t)\in\N^n$ be the population process (Markov jump process) representing the number of individuals in each of $n$ states at time $t$.  Let $\ntot \ge 1$ be the total number of individuals in the system.  Recall the random time-change representation in terms of Poisson processes given in Equation \ref{eq:population-process}: 
\begin{equation}\label{eq:population-process2}
N(t)=N(0)+\sum_{k\in\mathcal{E}}\zeta_k Y_k \left(\int_{0}^t\alpha_k N_{i(k)}(s)\,ds  \right).
\end{equation}
Each $Y_k$ is an independent unit rate Poisson process counting the occurrence of reaction $k$ (transition from state $i(k)$ to $j(k)$); $\alpha_k$ is the \textit{per capita} transition rate of reaction $k$; $N_{i(k)}(s)$ is the number of individuals at state $i(k)$ at time $s$, and $\zeta_k$ is the stoichiometry vector for reaction $k$.  For simplicity, we will suppress ``$k$'' in our notation so that state $i$ means state $i(k)$.  

In the case $\ntot=1$, let $p_i(t) = P(X(t)=i)$ be the probability that a single random walker occupies state $i$ at time $t$.  Clearly, $\sum_{i=1}^n p_i(t)=1$ for each $t$.  The time evolution of the probability vector $p(t)=[p_1(t), \dots, p_n(t)]$ is given by the following master equation
\begin{equation}
\frac{dp}{dt} = p L
\end{equation}
where
\begin{equation}\label{eq:Ldecomposition}
L = -\sum_{k\in\mathcal{E^*}}\alpha_k \zeta_k\zeta_k^\intercal
\end{equation}
is the graph Laplacian which can be represented as the sum over all undirected edges (denoted by the set $\mathcal{E^*}$) given in Equation \ref{eq:Ldecomposition}.  

Let $\pi$ represent the steady state distribution, i.e.~the row vector satisfying $\pi L=0$ with entries such that $\sum_{i=1}^n\pi_i=1$.  
Suppose we represent $N(t)$ as the deviation from its mean, $\bar{N} = \pi \ntot$, so that $N(t) = \bar{N} + X(t)$, where $X(t)$ is a mean zero stochastic process.  Then
\begin{eqnarray}\label{eq:Xderivation}
X(t) &=& N(t) - \bar{N} \\
&=& N(0) - \bar{N} + \sum_{k\in\mathcal{E}} \zeta_k Y_k \left( \int_0^t \alpha_k N_i(s)ds \right)\\
&=& X(0) + \sum_{k\in\mathcal{E}} \zeta_k Y_k \left( \int_0^t \alpha_k [\bar{N_i} + X_i(s)]ds \right)\\
&=& X(0) + \sum_{k\in\mathcal{E}} \zeta_k Y_k \left( t \alpha_k \bar{N_i} + \int_0^t \alpha_k X_i(s)ds \right)
\end{eqnarray}
since $N_i(s) = \bar{N_i}(s)+X_i(s)$ and $\alpha_k$ and $\bar{N_i}$ are constants.

Now following standard tau-leaping results \cite{Gillespie2001JChemPhys,PetzoldGillespie2003JChemPhys}, 
\begin{eqnarray}
X(t+\tau)-X(t) &=& \sum_{k\in\mathcal{E}} \zeta_k \left[ Y_k \left( (t+\tau) \alpha_k \bar{N_i} + \int_0^{t+\tau} \alpha_k X_i(s)ds \right) - Y_k \left( t \alpha_k \bar{N_i} + \int_0^t \alpha_k X_i(s)ds \right) \right]\label{eq:Xderivation-tau-leaping1}\\
&\approx& \sum_{k\in\mathcal{E}} \zeta_k \tilde{Y}_k \left( \tau \alpha_k \bar{N_i} + \tau \alpha_k X_i(t) \right)\label{eq:Xderivation-tau-leaping2}\\
&=& \sum_{k\in\mathcal{E}} \zeta_k \tilde{Y}_k \left( \tau \alpha_k [\bar{N_i} + X_i(t)] \right).
\end{eqnarray}
which says that we can approximate Equation \ref{eq:Xderivation-tau-leaping1} using an almost equivalent set of Poisson processes $\tilde{Y}_k$ where each $\tilde{Y}_k$ at time $t$ is approximately Gaussian distributed with mean and variance $\tau\alpha_k [\bar{N_i}+X_i(t)]$.  Since we are assuming that $|X_i(t)|<<\bar{N}_i$ (uniformly in time) and since we want the noise amplitude to be independent of $X$, we further approximate
\begin{equation} 
\tilde{Y}_k(\tau\alpha_k [\bar{N_i} + X_i(t)]) \approx \mathcal{N}(\tau\alpha_k [\bar{N_i} + X_i(t)], \tau\alpha_k\bar{N_i})
\end{equation}
by dropping the dependence of the variance on $X$. 

Dividing by $\tau$ and taking the limit as $\tau\to\infty$ yields the SDE
\begin{equation}\label{eq:dXtau}
dX = \sum_{k\in\mathcal{E}} \zeta_k \left( [\bar{N}_i+X_i]\alpha_k dt + \sqrt{\bar{N}_i \alpha_k} dW_k \right).
\end{equation}

Recalling that the $k^{th}$ reaction is from node $i(k)$ to $j(k)$, then the $k^{th}$ reaction in the first term in the RHS of Equation \ref{eq:dXtau} can be written as
\begin{equation}
\zeta_k (\bar{N}_l+X_l)\alpha_k dt = \begin{cases}
	-(\bar{N}_l+X_l)\alpha_k dt & \text{if component $l=i(k)$}\\
	(\bar{N}_i+X_i)\alpha_k dt & \text{if component $l=j(k)$}\\
	0 & \text{otherwise}
	\end{cases}
\end{equation}
Keeping track of components, we sum over the source and destination nodes for each reaction.  Then for the $l^{th}$ component of $X$ we have
\begin{equation}
dX_l = \sum_i(\bar{N}_i+X_i)\alpha_{il}dt - (\bar{N}_l+X_l)\sum_j \alpha_{lj}dt
\end{equation}
which yields
\begin{equation}
dX = (\bar{N}+X)Q dt \quad\mbox{ where }\quad 
(Q)_{ij} = \begin{cases}
	\alpha_{ij} & \text{if $i \ne j$}\\
	-\sum_{j \ne i}\alpha_{ij} & \text{if $i=j$}
	\end{cases}																									
\end{equation}
where $Q$ is the generator matrix.  Note that we changed notation slightly to illustrate that $\alpha_{ij}$ is the transition rate from state $i$ to $j$ rather than indexing by reaction $k$.  The graph Laplacian we consider in Equation \ref{eq:OUPforX} is actually $L = Q^\intercal$ so we have
$dX = L(\bar{N}+X)dt$.  Since $\bar{N}$ is proportional to the stationary distribution $\pi$, we have that $L\bar{N} = 0$, and hence the first term in the SDE is $dX = LXdt$. 

Now the second term in the RHS of Equation \ref{eq:dXtau} can be written as
\begin{equation}
\zeta_k \sqrt{\bar{N}_l\alpha_k} dW_k = \begin{cases}
	-\sqrt{\bar{N}_l \alpha_k} dW_k & \text{if component $l=i(k)$}\\
	\sqrt{\bar{N}_i \alpha_k} dW_k & \text{if component $l=j(k)$}\\
	0 & \text{otherwise}
	\end{cases}
\end{equation}
Keeping track of components, here we sum over all $m$ reactions to find
\begin{eqnarray}\label{eq:BdW}
dX &=& \matrix{c}{\sqrt{\bar{N}_{l(1)}\alpha_1}\zeta_1, ~ \sqrt{\bar{N}_{l(2)}\alpha_2}\zeta_2, ~\dots~, \sqrt{\bar{N}_{l(m)}\alpha_m}\zeta_m} \matrix{c}{dW_1\\ dW_2\\ \vdots\\ dW_m} \\
&=& B dW
\end{eqnarray}
where $\sigma_k=\sqrt{\bar{N}_{i(k)}\alpha_k}$ in the definition of matrix $B$.

Therefore, putting the first and second terms together, we have derived the OU process $dX = LXdt + BdW$ given in Equation \ref{eq:OUPforX}.

\subsection{Tau-Leaping: 3-State Example}

Here we will explicitly derive the OU process from the population process given in \S \ref{ssec:PopulationProcess} by using the tau-leaping argument above for the 3-state example in \S \ref{ssec:3state}.  We have $N(t)\in\N^3$ and by Equation \ref{eq:population-process2}, 
\begin{eqnarray}
N_1(t) &=& N_1(0) -Y_1 \left[\int_0^t N_1(s)\alpha_1 ds \right] +Y_2 \left[\int_0^t N_2(s)\alpha_2 ds \right] \\
N_2(t) &=& N_2(0) +Y_1 \left[\int_0^t N_1(s)\alpha_1 ds \right] -Y_2 \left[\int_0^t N_2(s)\alpha_2 ds \right]\nonumber \\ 
& & -Y_3 \left[\int_0^t N_2(s)\alpha_3 ds \right] +Y_4 \left[\int_0^t N_3(s)\alpha_4 ds \right] \\
N_3(t) &=& N_3(0) +Y_3 \left[\int_0^t N_2(s)\alpha_3 ds \right] -Y_4 \left[\int_0^t N_3(s)\alpha_4 ds \right] 
\end{eqnarray}
following the notation given in \S \ref{ssec:3state}, specifically the labeling of reactions given in Table \ref{tab:3-node-stoich}.  Note that $\alpha_k$ could be time dependent $\alpha_k(t)$.
%Here we use the standard notation $\alpha_{ij}$ to denote the instantaneous per capita transition rate from state $i$ to $j$

The tau-leaping approximation above gives
{\small
\begin{eqnarray}
X_1(t) &=& X_1(0) -\int_0^t X_1(s)\alpha_1 ds -\int_0^t \sqrt{X_1(s)\alpha_1} dW_1(s) +\int_0^t X_2(s)\alpha_2 ds +\int_0^t \sqrt{X_2(s)\alpha_2} dW_2(s)\\
X_2(t) &=& X_2(0) +\int_0^t X_1(s)\alpha_1 ds +\int_0^t \sqrt{X_1(s)\alpha_1} dW_1(s) -\int_0^t X_2(s)\alpha_2 ds -\int_0^t \sqrt{X_2(s)\alpha_2} dW_2(s)\nonumber \\ 
& & -\int_0^t X_2(s)\alpha_3 ds -\int_0^t \sqrt{X_2(s)\alpha_3} dW_3(s) +\int_0^t X_3(s)\alpha_4 ds +\int_0^t \sqrt{X_3(s)\alpha_4} dW_4(s)\\
X_3(t) &=& X_3(0) +\int_0^t X_2(s)\alpha_3 ds +\int_0^t \sqrt{X_2(s)\alpha_3} dW_3(s) -\int_0^t X_3(s)\alpha_4 ds -\int_0^t \sqrt{X_3(s)\alpha_4} dW_4(s). 
\end{eqnarray}}
Equivalently, we could write these integral equations in differential form  
\begin{eqnarray}
dX_1 &=& -X_1\alpha_1 dt -\sqrt{X_1\alpha_1} dW_1 +X_2\alpha_2 dt +\sqrt{X_2\alpha_2} dW_2\label{eq:dXdifferential}\\
dX_2 &=& X_1\alpha_1 dt +\sqrt{X_1\alpha_1} dW_1 -X_2\alpha_2 dt -\sqrt{X_2\alpha_2} dW_2\nonumber\\ 
& & -X_2\alpha_3 dt -\sqrt{X_2\alpha_3} dW_3 +X_3\alpha_4 dt +\sqrt{X_3\alpha_4} dW_4\\
dX_3 &=& X_2\alpha_3 dt +\sqrt{X_2\alpha_3} dW_3 -X_3\alpha_4 dt -\sqrt{X_3\alpha_4} dW_4\label{eq:dXdifferential2}. 
\end{eqnarray}
%$N_1(t)+N_2(t)+N_3(t) = N_{tot}$, constant
These equations are nonlinear since the noise intensity depends on $X_i$.  Note that for any $t$, $X_1(t)+X_2(t)+X_3(t) = N_{tot}$ so that the total population is constant.  The mean $\bar{X}$ satisfies
\begin{eqnarray}
\frac{d\bar{X}}{dt} = \bar{X} \matrix{ccc}{-\alpha_1 & \alpha_1 & 0\\ \alpha_2 & -(\alpha_2+\alpha_3) & \alpha_3\\ 0 & \alpha_4 & -\alpha_3}
\end{eqnarray}
where the matrix above is the generator $Q$, or our $L^\intercal$.  In the case where $Q$ is fixed, $\bar{X}$ is proportional to the null left eigenvector of $Q$; biologically, this is the voltage clamp case.  Let $(\bar{X}_1, \dots, \bar{X}_n)$ be the corresponding stationary vector.  Now we linearize Equations \ref{eq:dXdifferential}-\ref{eq:dXdifferential2} around the stationary vector.

Let $V = X-\bar{X}$ and assume that $\frac{|V|}{\bar{X}}<<1$.  Then since $\sqrt{X_i \alpha_k} = \sqrt{(\bar{X}_i + V_i)\alpha_k} = \sqrt{\bar{X}_i\alpha_k} + O\left(\frac{V_i}{\bar{X}_i}\right)$, we have
\begin{eqnarray}
dV_1 &=& (-V_1\alpha_1 + V_2\alpha_2)dt -\sqrt{X_1\alpha_1}dW_1 +\sqrt{X_2\alpha_2}dW_2 + O\left(\frac{|V|}{N_{tot}}\right)\\
dV_2 &=& (V_1\alpha_1 - V_2\alpha_2 - V_2\alpha_3 + V_3\alpha_4)dt +\sqrt{V_1\alpha_1}dW_1 -\sqrt{V_2\alpha_2}dW_2\nonumber\\ 
& & -\sqrt{V_2\alpha_3}dW_3 +\sqrt{V_3\alpha_4}dW_4 + O\left(\frac{|V|}{N_{tot}}\right)\\
dV_3 &=& (V_2\alpha_3 - V_3\alpha_4)dt +\sqrt{V_2\alpha_3}dW_3  -\sqrt{V_3\alpha_4}dW_4 + O\left(\frac{|V|}{N_{tot}}\right). 
\end{eqnarray}
Neglecting the $O\left(\frac{|V|}{N_{tot}}\right)$ terms gives us the multidimensional Ornstein-Uhlenbeck process of Equation \ref{eq:OUPforX} for the 3-state example.

\section{Proofs and Calculations}
\label{app:proofs}

\subsection{Stationary Covariance of a Multidimensional OU Process}
\label{app:OUP-covariance}

The SDE for $X(t)$ in Equation \ref{eq:OUPforX} has explicit solution (see \cite{Gardiner2009book}, chapter 4.5)
\begin{equation}
X(t) = \exp(Lt) X(0) + \int_0^t{\exp(L(t-t')) B dW(t')}.
\end{equation}
Assuming the initial condition is either deterministic or Gaussian, then $X(t)$ is Gaussian with mean 
\begin{equation}
E[X(t)] = \exp(Lt) E[X(0)]
\end{equation}
and correlation function 
\begin{eqnarray}\label{eq:correlationX}
\text{Cov}[X(t),X^\intercal (s)] &=& \exp(Lt)E[X(0),X^\intercal (0)]\exp(Ls) \nonumber\\
& & +\int_0^{t\wedge s}\exp[L(t-t')]BB^\intercal \exp[L^\intercal (s-t')]\,dt',
\end{eqnarray}
where $t\wedge s$ means the minimum of $t$ and $s$. 
%%
%If we take the initial condition $X(0)$ from the equilibrium distribution, then $E[X(0),X^\intercal(0)] = 0$ and Equation \ref{eq:correlationX} simplifies to\footnote{QQQ Check -- this doesn't look quite right!}
%\begin{equation}\label{eq:correlationXsimple}
%\cov[X(t),X^\intercal(s)] = \int_0^{t\wedge s}\exp[-L(t-t')]BB^\intercal\exp[-L^\intercal (s-t')]\,dt'.
%\end{equation}
%
Setting $s=t$ and taking the limit as $t\to\infty$, we obtain the stationary covariance function
\begin{equation}\label{eq:stationarycovX}
\cov[X(t),X^\intercal (t)] = \lim_{t\to\infty}\int_0^t\exp[L(t-t')]BB^\intercal \exp[L^\intercal (t-t')]\,dt'.
\end{equation}
We exploit the fact that not only does $B$ decompose into the sum $B=\sum_{k=1}^m B_k$, but in the case of a first-order reaction process, $BB^\intercal$ also decomposes into the following sum 
\begin{equation}\label{eq:decompBBt}
BB^\intercal=\sum_{k=1}^m B_k B_k^\intercal, 
\end{equation}
and futher, $B_k B_k^\intercal = \sigma_k^2 \zeta_k \zeta_k^\intercal$ for each edge (reaction) $k\in\mathcal{E}$.  
Therefore, the stationary covariance of the full process decomposes into a sum of the contributions from the $m$ different reaction processes:
\begin{equation}\label{eq:covariance-of-X}
\cov[X(t),X^\intercal (t)] = \lim_{t\to\infty}\int_0^t\sum_{k=1}^m \sigma_k^2\exp[L(t-t')]\zeta_k\zeta_k^\intercal \exp[L^\intercal (t-t')]\,dt'.
\end{equation}
We note that the eigenvector corresponding to the leading (0) eigenvalue of $L$ has constant components, therefore it lies in the kernel of the matrix $B_kB_k^\intercal$ for each $k$, which guarantees finite covariance in Equation \ref{eq:covariance-of-X}.

\subsection{Computation of Edge Importance $R_k$ and Proof of Lemma \ref{lem:decompose-U}}
\label{subsec:proof-lemma1}
Using the spectral properties of the graph Laplacian $L$, we can rewrite the stationary covariance of $X(t)$ (Equation \ref{eq:stationarycovX}) by replacing each expression involving a matrix exponential by the sum over the orthogonal eigendecomposition of $L$.  Let $v_i$ be the $i^{th}$ eigenvector of $L$ (written as a column vector), with eigenvalue $\lambda_i$, i.e. $L v_i = \lambda_i v_i$.  Summing over each eigenvalue, we can write $L = \sum\limits_{i=1}^n \lambda_i v_i v_i^\intercal$.  Note that this decomposition is only valid when $L$ is symmetric; the non-symmetric case is discussed in \S \ref{sec:HH}.  Then we have the following expression from Equation \ref{eq:stationarycovX}
\begin{eqnarray}\label{eq:QBBQ-madness}
&&\exp[L(t-t')]BB^\intercal \exp[L^\intercal (t-t')]\\
&=& \left(\sum_{i=1}^n e^{\lambda_i(t-t')} v_i v_i^\intercal \right) BB^\intercal \left(\sum_{j=1}^n e^{\lambda_j(t-t')} v_j v_j^\intercal \right)\\
&=& \sum_{i,j=1}^n e^{(\lambda_i+\lambda_j)(t-t')}(v_i v_i^\intercal ) (BB^\intercal ) (v_j v_j^\intercal). 
\end{eqnarray}

Using the decomposition of matrix $B$ (Equations \ref{eq:decompB} and \ref{eq:decompBBt}), it follows that
\begin{equation}
BB^\intercal=\sum_{k=1}^m B_kB_k^\intercal = \sum_{k=1}^m \sigma_k^2 \zeta_k \zeta_k^\intercal.
\end{equation}
The covariance of the full process $X$ is therefore given by 
\begin{eqnarray} \label{eq:covarianceXsimple}
\cov[X(t),X^\intercal(t)] &=& \int_0^t \sum_{i,j=1}^n e^{(\lambda_i+\lambda_j)(t-t')}(v_i v_i^\intercal) (BB^\intercal) (v_j v_j^\intercal)\,dt'\\
&=& \sum_{k=1}^m \sigma_k^2 \sum_{i=2}^n\sum_{j=2}^n
\left(\frac{1-e^{(\lambda_i+\lambda_j)t}}{-(\lambda_i+\lambda_j)}\right) (v_i v_i^\intercal) (\zeta_k\zeta_k^\intercal) (v_j v_j^\intercal). 
\label{eq:covarianceXsimple2}
\end{eqnarray}
By construction of the graph Laplacian, its leading eigenvalue $\lambda_1\equiv 0$.  The corresponding (right) eigenvector has constant components, $v_1=(1,\cdots,1)^\intercal/\sqrt{n}$.  Therefore, for each stoichiometry vector we have $\zeta_k^\intercal v_1\equiv 0$.  Consequently the terms in the inner summation (\ref{eq:covarianceXsimple2}) with index $i=1$ or $j=1$ vanish, and may be omitted without changing the result.
Taking the limit as $t\to\infty$ of the covariance function gives us the stationary covariance
\begin{eqnarray} \label{eq:stationarycovarianceX}
\cov[X(t),X^\intercal(t)] &=& \lim_{t\to\infty} \sum_{k=1}^m \sigma_k^2 \sum_{i=2}^n\sum_{j=2}^n
\left(\frac{1-e^{(\lambda_i+\lambda_j)t}}{-(\lambda_i+\lambda_j)}\right) (v_i v_i^\intercal) (\zeta_k\zeta_k^\intercal) (v_j v_j^\intercal)\\
&=& \sum_{k=1}^m\sigma_k^2 \sum_{i=2}^n\sum_{j=2}^n
\left(\frac{-1}{\lambda_i+\lambda_j}\right) (v_i v_i^\intercal) (\zeta_k\zeta_k^\intercal) (v_j v_j^\intercal). 
\end{eqnarray}
Recall that we are interested in the linear measurement functional $M\in\R^n$ projected onto $X(t)$, i.e. the projection $Y(t) = M^\intercal X(t)$. For edges $k\in\mathcal{E}'$ neglected in the approximation $\tilde{Y}=M^\intercal \tilde{X}(t)$, we take the limit as $t\to\infty$ of the mean squared error of $\tilde{Y}(t)-Y(t) = M^\intercal U(t)$ to get
\begin{eqnarray} \label{eq:errorU}
R[\mathcal{E}']&=&\lim_{t\to\infty}E\left[||(\tilde{Y}(t)-Y(t))||^2_2 \right]\\
&=&\lim_{t\to\infty}E\left[||M^\intercal U(t)||^2_2 \right]\\
&=&\lim_{t\to\infty}\left(M^\intercal\cov[U(t),U^\intercal(t)]M\right)\\
&=&\sum_{k\in\mathcal{E}'}\sigma_k^2 \sum_{i=2}^n\sum_{j=2}^n
\left(\frac{-1}{\lambda_i+\lambda_j}\right) (M^\intercal v_i) (v_i^\intercal \zeta_k) (\zeta_k^\intercal v_j) (v_j^\intercal M)\label{eq:errorU2}\\
&=&\sum_{k\in\mathcal{E}'} R_k.
\end{eqnarray}

%Note that the sum $\sum_{i=2}^n\sum_{j=2}^n\left(\frac{1}{\lambda_i+\lambda_j}\right)$ is negative because all non-zero eigenvalues are either negative or have negative real part, making the rest of the sum in Equation \ref{eq:errorU2} negative. In defining $R_k$, we choose to make it positive by multiplying $S$ by $-1$.  Therefore, $R_k = -R_k^*$ is given by
%\begin{equation}
%R_k := \sigma_k^2\sum_{i=2}^n\sum_{j=2}^n \left(\frac{-1}{\lambda_i+\lambda_j}\right) (M^\intercal v_i) (v_i^\intercal \zeta_k) (\zeta_k^\intercal v_j) (v_j^\intercal M).
%\end{equation}

\subsection{Proof of Lemma \ref{lem:random}}
\label{ssec:proof_lemma2}

Suppose that assumptions A0-A5 given in \S \ref{ssec:assumptions} hold.  Recall that $M \in\{0,1\}^n$ is an arbitrary a measurement vector.  Suppose there are $n_1>0$ ones and $n_0>0$ zeros such that $n_1+n_0=n$.  We assume $n_1 = O(1)$, that is, we exclude the case where $n_1$ grows without bound as $n$ grows.  If we look at the corresponding measurement value of the $l_-^{th}$ and $l_+^{th}$ components of $\zeta_k$ (see Equation \ref{eq:zeta_k_2}), we have three possible cases:
\begin{enumerate}\label{list:cases} 
\item $l_{\pm}\in 1_{M}$, i.e. $M(l_-)=M(l_+)=1$
\item $l_{\pm}\notin 1_{M}$, i.e. $M(l_-)=M(l_+)=0$
\item $l_-\in 1_{M}$ and $l_+ \notin 1_{M}$, i.e. $M(l_-)=1$ and $M(l_+)=0$ (respectively, $M(l_-)=0$ and $M(l_+)=1$, equivalent up to a sign change) 
\end{enumerate}
For each part of Lemma \ref{lem:random}, we will prove the result for these three cases.  If we let $n_1^*$ denote the number of terms in the set $1_M\backslash\{l_\pm\}$, then 
\begin{equation}
  n_1^*=\begin{cases}
    n_1-2, & \text{if $l_{\pm}\in 1_{M}$ (Case 1)}\\
    n_1, & \text{if $l_{\pm}\notin 1_{M}$ (Case 2)}\\
    n_1-1, & \text{if $l_-\in 1_{M}$ and $l_+\notin 1_{M}$ (Case 3)}\\
  \end{cases}
\end{equation}
and we can consider all three cases at once using this notation where now $n_1^{*}=O(1)$. 

Let $a=v_i(l_-), b=v_i(l_+)$, and $c=\sum_{l\in 1_M\backslash\{l_{\pm}\}}v_i(l)$.  By assumption A2, we have that $E[a]=E[b]=E[c]=0$ and $E[a^2]=E[b^2]=n^{-1}$ from the normalization of the eigenvectors, and it follows from A3b that $E[c^2]=(n_1^*)n^{-1}+O(n^{-3})$, as $n\to\infty$.  Assumption A3 gives conditions on second order terms. Assumptions A4 and A5 give conditions on fourth order moments and fourth order products of $a, b$ and $c$.

%%%%%%%%%%%%%%%%%%%%%%%%%%%%%%%%
% Proof of Part A
%%%%%%%%%%%%%%%%%%%%%%%%%%%%%%%%
\subsubsection{Proof of Part A}

We will show that, as $n\to\infty$,
\begin{equation}
E[M^\intercal v_i v_i^\intercal \zeta_k] = \frac{1}{n}(M^\intercal\zeta_k) + O\left(\frac{1}{n^2}\right).
\end{equation}
By definition
\begin{equation}\label{eq:definition}
E[M^\intercal v_i v_i^\intercal \zeta_k] = E\left[\sum_{l\in 1_M} v_i(l) (v_i(l_+)-v_i(l_-))\right]
\end{equation}
since $M^\intercal v_i = \sum_{l\in 1_M} v_i(l)$ and $v_i^\intercal \zeta_k = v_i(l_+)-v_i(l_-)$. We compute this expectation for the three cases listed at the beginning of \S \ref{list:cases}.

Using the notation introduced above, we note that this expectation has the form:
\begin{eqnarray}
& &E[(a+b+c)(b-a)] \quad\text{for Case 1}\\
& &E[c(b-a)] \quad\text{for Case 2}\\
& &E[(a+c)(b-a)] \quad\text{for Case 3}.
\end{eqnarray}

\emph{Case 1:} $l_{\pm}\in 1_M$

Expanding the expected value yields
\begin{eqnarray}
E[(a+b+c)(b-a)] &=& E[b^2 - a^2 + bc -ac]\\
&=& E[b^2] - E[a^2] + E[bc] - E[ac]\\
&=& \frac{1}{n} - \frac{1}{n} + E[bc] - E[ac] 
\end{eqnarray}
since $E[a^2]=E[b^2]=n^{-1}$ by assumption A2 (eigenvector normalization).  Note that $E[ac]=E[bc]$, and each contains $n_1^{*}$ terms with the following expectation as $n\to\infty$:
\begin{eqnarray}\label{eq:Eac}
E[ac] &=& E[v_i(l_-) \sum_{l\in 1_M\backslash\{l_{\pm}\}}v_i(l)] = \sum_{l\in 1_M\backslash\{l_{\pm}\}} E[v_i(l_-) v_i(l)] = n_1^{*} O(n^{-2})\\ &=& O(n^{-2}).
\end{eqnarray}
This follows from the assumptions that, as $n\to\infty$, $E[v_i(l)v_i(l')]=O(n^{-2})$ for $l\ne l'$ (A3b) and $n_1^{*}=O(1)$.  Thus, since $M^\intercal \zeta_k = -1+1 = 0$ in this case, as $n\to\infty$, 
\begin{equation}
E[(a+b+c)(b-a)] = \frac{1}{n}(M^\intercal \zeta_k) + O(n^{-2}).
\end{equation}

\emph{Case 2:} $l_{\pm}\notin 1_M$

Expanding the expected value yields
\begin{equation}
E[c(b-a)] = E[bc] - E[ac] = O(n^{-2})
\end{equation}
as $n\to\infty$, by Equation \ref{eq:Eac} in Case 1 above.  Thus, since $M^\intercal \zeta_k = 0$ in this case, as $n\to\infty$, 
\begin{equation}
E[c(b-a)] = \frac{1}{n}(M^\intercal \zeta_k) + O(n^{-2}).
\end{equation}

\emph{Case 3:} $l_-\in 1_M$ and $l_+\notin 1_M$

Expanding the expected value yields
\begin{eqnarray}
E[(a+c)(b-a)] &=& E[-a^2 + ab + bc - ac]\\
&=& -E[a^2] + E[ab] + E[bc] - E[ac]\\
&=& -\frac{1}{n} + O(n^{-2})
\end{eqnarray}
as $n\to\infty$, which follows by Equation \ref{eq:Eac} from Case 1 and by the assumptions that $E[v_i(l)v_i(l')]=O(n^{-2})$ for $l\ne l'$ (A3b) and $n_1^{*}=O(1)$, as $n\to\infty$.  Since $M^\intercal \zeta_k = -1$ in this case, then as $n\to\infty$, 
\begin{equation}
E[(a+c)(b-a)] = \frac{1}{n}(M^\intercal \zeta_k) + O(n^{-2}).
\end{equation}

Similarly, the alternate Case 3 where $l_+\in 1_M$ and $l_-\notin 1_M$ gives, as $n\to\infty$, 
\begin{eqnarray}
E[(b+c)(b-a)] &=& E[-a^2 + ab + bc - ac]\\
&=& E[b^2] - E[ab] + E[bc] - E[ac]\\
&=& \frac{1}{n} + O(n^{-2}) 
\end{eqnarray}
and since $M^\intercal \zeta_k = 1$ in this case, we have as $n\to\infty$
\begin{equation}
E[(b+c)(b-a)] = \frac{1}{n}(M^\intercal \zeta_k) + O(n^{-2}).
\end{equation}

%%%%%%%%%%%%%%%%%%%%%%%%%%%%%%%%
% Proof of Part B
%%%%%%%%%%%%%%%%%%%%%%%%%%%%%%%%
\subsubsection{Proof of Part B}

We will show that, as $n\to\infty$,
\begin{equation}\label{eq:partBresult}
E[M^\intercal v_i v_i^\intercal \zeta_k]^2 = \frac{1}{n^2}|M^\intercal\zeta_k| + O\left(\frac{1}{n^4}\right)
\end{equation}
where now we take the absolute value of the term $M^\intercal\zeta_k$.  By definition (see Equation \ref{eq:definition}), we have that 
\begin{equation}
E[M^\intercal v_i v_i^\intercal \zeta_k]^2 = E\left[\sum_{l\in 1_M} v_i(l) (v_i(l_+)-v_i(l_-))\right]^2.
\end{equation}

Using the notation introduced above, this expectation has the following structure in each case:
\begin{eqnarray}
& &E[(a+b+c)(b-a)]^2 \quad\text{for Case 1}\\
& &E[c(b-a)]^2 \quad\text{for Case 2}\\
& &E[(a+c)(b-a)]^2 \quad\text{for Case 3}.
\end{eqnarray}

By Lemma \ref{lem:random} part A, we have that, as $n\to\infty$,
\begin{eqnarray}
E[(a+b+c)(b-a)] &=& 0 + O(n^{-2})\\
E[c(b-a)] &=& 0 + O(n^{-2})\\
E[(a+c)(b-a)] &=& -\frac{1}{n} + O(n^{-2})\\
E[(b+c)(b-a)] &=& \frac{1}{n} + O(n^{-2})  
\end{eqnarray}
where the last two equations fall under Case 3.  Squaring these terms yields, as $n\to\infty$,
\begin{eqnarray}
E[(a+b+c)(b-a)]^2 &=& 0 + O(n^{-4})\\
E[c(b-a)]^2 &=& 0 + O(n^{-4})\\
E[(a+c)(b-a)]^2 &=& \frac{1}{n^2} + O(n^{-4})\\
E[(b+c)(b-a)]^2 &=& \frac{1}{n^2} + O(n^{-4}).  
\end{eqnarray}
In this case, both versions of Case 3 are positive so we multiply $1/n^2$ by $|M^\intercal \zeta_k|$ which gives us the desired result in Equation \ref{eq:partBresult}.

%%%%%%%%%%%%%%%%%%%%%%%%%%%%%%%%
% Proof of Part C
%%%%%%%%%%%%%%%%%%%%%%%%%%%%%%%%
\subsubsection{Proof of Part C}

We will show that, as $n\to\infty$,
\begin{equation}
E[(M^\intercal v_i v_i^\intercal \zeta_k)^2] =O(n^{-q}) \text{ for some } q>1.
\end{equation}

It follows by definition that
\begin{equation}
E[(M^\intercal v_i v_i^\intercal\zeta_k)^2] = E[(M^\intercal v_i) (v_i^\intercal \zeta_k\zeta_k^\intercal v_i) (v_i^\intercal M)] = E\left[\left(\sum_{l\in 1_M} v_i(l) \right)^2 (v_i(l_+)-v_i(l_-))^2 \right ],
\end{equation}
since $M^\intercal v_i = v_i^\intercal M = \sum_{l\in 1_M} v_i(l)$ and $v_i^\intercal \zeta_k\zeta_k^\intercal v_i = (v_i(l_+)-v_i(l_-))^2$.

Note that this term has the following structure in each case:
\begin{eqnarray}
& &E[(a+b+c)^2(b-a)^2] \quad\text{for Case 1}\label{eq:case1sum2}\\
& &E[c^2(b-a)^2] \quad\text{for Case 2}\\
& &E[(a+c)^2(b-a)^2] \quad\text{for Case 3}\label{eq:case3sum2}.
\end{eqnarray}
Expanding the sums above (Equations \ref{eq:case1sum2}-\ref{eq:case3sum2}), we see that all but one term for Cases 2 and 3 also appear in Case 1, and that term $E[a^3b]$ is of smaller order of magnitude than $E[a^3c]$ which appears in Case 1.  Thus, it suffices to consider only Case 1.  Expanding the sum (Equation \ref{eq:case1sum2}) gives
\begin{eqnarray}
E[(a+b+c)^2(b-a)^2]&=&E[(a^2+b^2+c^2+2ab+2ac+2bc)(a^2-2ab+b^2)]\\
&=& E[a^4-2a^2b^2+b^4 + a^2c^2-2abc^2+b^2c^2\\
&&-2ab^2c - 2a^2bc + 2a^3c + 2b^3c]\\
&=& E[a^4] + E[b^4] + O(n^{-2}) \text{ as } n\to\infty.
\end{eqnarray}

The leading order terms are $E[a^4]=E[b^4]=O(n^{-q})$ as $n\to\infty$ for some $q>1$ by assumption A4a and the term $E[a^2 b^2]=O(n^{-2})$ as $n\to\infty$ by assumption A4b.  Note that all terms involving powers of $c$ carry an extra factor of $n_1^*$ (or $(n_1^*)^2$), but this doesn't change the order of magnitude since we're assuming $n_1^*=O(1)$.  Therefore, the terms $E[a^2c^2]$ and $E[b^2c^2]$ are also $O(n^{-2})$, as shown below.  As $n\to\infty$
\begin{eqnarray}
E[a^2c^2] &=& E\left[ a^2 \left( \sum_{l\in1_M\backslash\{l_\pm\}} v_i(l) \right)^2 \right]\\
&=& E\left[ a^2 v_i(l_1)^2 +\dots+ a^2 v_i(l_{n_1^*})^2 + \sum_{l_j,l_k \in1_M\backslash\{l_\pm\},j\ne k} a^2 v_i(l_j)v_i(l_k)\right]\\
&=& n_1^* E\left[ a^2 v_i(l_1)^2\right] + {n_1^* \choose 2} E\left[a^2 v_i(l_1)v_i(l_2)\right]\\
&=& O(n_1^* n^{-2}) + O((n_1^*)^2 n^{-3}) \quad\text{by Assumptions A4b and A5}\\
&=& O(n^{-2}).
\end{eqnarray}

The same holds for $E[b^2c^2]$ since $E[a^2c^2]=E[b^2c^2]$.  We can do a similar calculation for $E[abc^2]$, replacing $a^2$ with $ab$, and noting that assumption A5 holds for terms of the form $a b v_i(l_1)$ and $a b v_i(l_1)v_i(l_2)$ with distinct eigenvector components. Hence, $E[abc^2]= O(n^{-3})$ as $n\to\infty$.

All other cross terms ($E[ab^2c]$, $E[a^2bc]$, $E[a^3c]$, $E[b^3c]$) are of order $O(n_1^* n^{-3}) = O(n^{-3})$ as $n\to\infty$ by assumption A5.  Therefore, since the leading order terms are $O(n^{-q})$, it follows that 
\begin{equation}
E[(a+b+c)^2(b-a)^2]=O(n^{-q}) \text{ as } n\to\infty \text{ for some } q>1.
\end{equation}

\end{document}